\documentclass[twoside]{amsart}
\usepackage{float}
\usepackage[bookmarksnumbered,plainpages,hypertex]{hyperref}
\usepackage{graphicx}
\usepackage{amssymb}
\usepackage[usenames]{color}

 \input{xy}
 \xyoption{all}
\newtheorem{theorem}{\sc Theorem}[section]
\newtheorem{proposition}[theorem]{\sc Proposition}
\newtheorem{lemma}[theorem]{\sc Lemma}
\newtheorem{corollary}[theorem]{\sc Corollary}
\theoremstyle{definition}
\newtheorem{definition}[theorem]{\sc Definition}

\newtheorem{example}[theorem]{\sc Example}

\theoremstyle{remark}
\newtheorem{remark}[theorem]{\sc Remark}

\newtheorem{claim}[theorem]{\sc}

\def\Hom{\mathrm{Hom}}

\def\ot{\otimes}
\def\can{\mathrm{can}}

\def\cM{\mathfrak M}
\def\cC{\mathcal C}
\def\cD{\mathcal D}
\def\hH{\mathcal H}
\def\stac#1{\raise-.2cm\hbox{$\stackrel{\,\displaystyle\otimes\,}{\scriptscriptstyle{#1}}$}}
\def\sstac#1{\otimes_{#1}}

\numberwithin{equation}{section}
\begin{document}
\title{A Schneider type Theorem for Hopf Algebroids}
\author{A. Ardizzoni}
\address{University of Ferrara, Department of Mathematics, Via Machiavelli
35, Ferrara, I-44100, Italy} \email{rdzlsn@unife.it}
\urladdr{http://www.unife.it/utenti/alessandro.ardizzoni}
\author{G. B\"{o}hm}
\address{Research Institute for Particle and Nuclear Physics, Budapest, H-1525
  Budapest 114, P.O.B.49\\
\indent Hungary } \email{G.Bohm@rmki.kfki.hu}
\urladdr{http://www.rmki.kfki.hu/$\sim$bgabr}
\author{C. Menini}
\address{University of Ferrara, Department of Mathematics, Via Machiavelli
35, Ferrara, I-44100, Italy} \email{men@unife.it}
\urladdr{http://www.unife.it/utenti/claudia.menini}
\date{December 2006, revised May 2007}
\begin{abstract}
Comodule algebras of a Hopf algebroid ${\mathcal H}$ with a
bijective antipode, i.e. algebra extensions {$B\subseteq A$} by
$\hH$, are studied.
{Assuming that a lifted canonical map is a split epimorphism of
modules
  of the (non-commutative) base algebra of $\hH$, {\em relative injectivity} of
  the $\hH$-comodule algebra $A$ is related to the {\em Galois property} of
  the extension $B\subseteq A$ and also to the {\em equivalence} of the
  category of relative Hopf modules to the category of $B$-modules.}
This extends a classical theorem by H.-J. Schneider on Galois extensions by a
Hopf algebra.
 Our main tool is an observation that relative injectivity of a
comodule algebra is equivalent to {\em relative separability} of a forgetful
functor, a notion introduced and analyzed hereby.

{\color{blue}
In the first version of this
submission, we heavily used the statement that two constituent
bialgebroids in a Hopf algebroid possess isomorphic comodule categories.
This statement was based on \cite[Theorem 2.6]{Brz:Corext}, whose proof turned
out to contain an unjustified step. In the revised version we return to an
earlier definition of a comodule of a Hopf algebroid, that distinguishes
between comodules of the two constituent bialgebroids, and modify the
statements and proofs in the paper accordingly.
}
\end{abstract}

 \maketitle
\tableofcontents

\pagestyle{headings}

\section{Introduction}

Galois extensions of non-commutative algebras by a Hopf algebra
generalize Galois extensions of commutative rings by groups and
are known as the algebraic (dual) versions of (non-commutative)
principal bundles. By a Hopf Galois extension the following
structure is meant. Comodules over a Hopf algebra $H$ form a
monoidal category ${\mathfrak M}^H$, whose monoids are called {\em
comodule algebras}. This means an algebra and $H$-comodule $A$,
such that the coaction $\rho^A:A\to A\otimes H$ is an algebra map
(with respect to the tensor product algebra structure of the
codomain). It can be looked at as {a} notion dual to the action of
a group on a manifold. Dualizing the notion of invariant points,
{\em coinvariants} of $A$ are defined as those elements on which
coaction is trivial, i.e. the elements of the subalgebra
$$
B: = \{\ b\in A\ \vert\ \rho^A(b)=b\otimes 1_H\ \}.
$$
In this situation the algebra $A$ is called an {\em extension of
$B$ by $H$}. The algebra extension $B\subseteq A$ is said to be
{\em $H$-Galois} if in addition the so called {\em canonical map}
\begin{equation}\label{eq:can.Hopf}
\can : A\stac B A \to A\otimes H, \qquad a\stac B a' \mapsto a\rho^A(a')
\end{equation}
is bijective (hence an isomorphism of left $A$-modules and right
$H$-comodules). This is a dual formulation of the condition that a group
action on a manifold is free.

(Right-right) {\em relative Hopf modules} are (right) modules for an
$H$-comodule algebra $A$ and (right) comodules for the Hopf algebra $H$,
satisfying a compatibility condition with the $H$-coaction in $A$. In the case
of an $H$-Galois extension $B\subseteq A$, relative Hopf modules are
canonically identified with descent data for the extension $B\subseteq A$.
Hence if $A$ is faithfully flat as a left $B$-module, it follows by the
Faithfully Flat Descent Theorem that the category ${\mathfrak M}^H_A$ of
right-right relative Hopf modules is equivalent to the category $\cM_B$ of
right $B$-modules.

In the study of Hopf Galois extensions, important tools are provided by
theorems,
stating that in appropriate situations surjectivity of the canonical map
\eqref{eq:can.Hopf}
implies its bijectivity. One group of such results (e.g. \cite[Theorem
1.7]{KreTak:HopfGal}, \cite[Corollary 2.4.8 1]{Scha:HGal&biGal}
\cite[Theorem 3.1]{SchSch:genGal}, \cite[Theorem 4.2 and Corollary
4.3]{Bohm:hgdGal}) can be called {\em
`Kreimer-Takeuchi type'} theorems (as their first representative was proven
in \cite[Theorem 1.7]{KreTak:HopfGal}). In this group of theorems
projectivity of the regular comodule of the
coacting Hopf algebra is assumed. The other group involves {\em
`Schneider type'} theorems (after \cite[Theorem I]{Schn:PriHomS}, see
e.g. \cite[Theorem 4.9]{SchSch:genGal}, \cite[Theorem 4.6]{Brz:Galcom},
\cite[Theorem 3.15]{MenMil:intYD}, \cite[Theorem 3.9]{MenMil:affDK}). Here
relative injectivity of the Hopf comodule algebra in question is assumed.

The starting point of our work is an observation that the proofs of all
above theorems share a common philosophy. Related to a comodule algebra $A$ of
a Hopf algebra $H$ over a commutative ring $k$, there are forgetful functors
\begin{equation}\label{fig:forgfHa}
\xymatrix{
{\cM}^H_A\ar[r]^{\mathbb R} &{\cM}^H\ar[r]^{\mathbb U}&
{\cM}_k.
}
\end{equation}
If $H$ is a projective $k$-module then the codomain of the canonical map
\eqref{eq:can.Hopf} is a
projective $A$-module. Then it
follows from the surjectivity of the canonical map that its lifted version
\begin{equation}\label{eq:can.lifted}
{\widetilde {\can}}: A\stac k A \to A \stac k H,\qquad a\stac k a'
\mapsto
  a\rho^A(a')
\end{equation}
{has a $k$-linear right inverse, i.e. it} is a retraction of
$k$-modules. The various Schneider type theorems give sufficient
conditions for the forgetful functor ${\mathbb U}$ to reflect
(certain) retractions. Then bijectivity of the canonical map
\eqref{eq:can.Hopf} follows by a result of Schauenburg
\cite[Corollary 2.4.8]{Scha:HGal&biGal} stating that -- under the
additional assumption that all $H$-coinvariants of the obvious
right $H$-comodule $A\otimes_k A$ are elements of $A\otimes_k B$
-- the canonical map is bijective, provided that its lifted
version \eqref{eq:can.lifted} is a retraction of $H$-comodules.

In the present paper we introduce the notion of {\em separability}
of a functor ${\mathbb U}:{\mathfrak A}\to {\mathfrak B}$, {\em
relative} to a functor ${\mathbb R}:{\mathfrak R}\to {\mathfrak
A}$ (not to be mixed with separability of the second kind in
\cite{CaeMil:SecK}). An ${\mathbb R}$-relative separable functor
${\mathbb U}$ reflects retractions in the sense that, for a
morphism $f$ in ${\mathfrak R}$ such that ${\mathbb U}{\mathbb
R}(f)$ is a retraction, ${\mathbb R}(f)$ is a retraction. As it
turns out, the conditions of all known Schneider type theorems
imply the separability of the forgetful functor ${\mathbb U}$ in
\eqref{fig:forgfHa}, relative to ${\mathbb R}$.

Our strategy, of tracing back Schneider type theorems to
properties of a forgetful functor, can be compared to that of Caenepeel, Ion,
Militaru and Zhu, when in \cite{CaeIonMilZhu:SepDoiH} they explained all known
Maschke type theorems by the separability of a forgetful functor.

The motivation of our work comes from a wish to prove a Schneider
type theorem for more general algebra extensions by a Hopf
algebroid, replacing the Hopf algebra $H$ above. A
Kreimer-Takeuchi type theorem was proven in \cite{Bohm:hgdGal}. In
that paper similar methods have been used as in
\cite{SchSch:genGal}: the entwining structure (over a
non-commutative base), determined by a comodule algebra of a Hopf
algebroid, has been studied. It turns out that this framework is
not sufficient to obtain a Schneider type theorem for extensions
by Hopf algebroids. Recall that a Hopf algebroid $\hH$ consists of
two related coring (and bialgebroid) structures, over two
different base algebras $L$ and $R$.
{\color{blue} The proper definition of an $\hH$-comodule consists
of a compatible pair of comodules, one for
each constituent coring. This results in a monoidal
category $\cM^\hH$ of $\hH$-comodules. By definition, a right
$\hH$-comodule algebra is an algebra in $\cM^\hH$. As in the Hopf
algebra case, right $A$-modules {\em in} $\cM^\hH$ are called {\em relative
Hopf modules}. Their category $\cM^\hH_A$ admits forgetful functors}
\begin{equation}\label{eq:R&U}
\xymatrix{
{\cM}^\hH_A\ar[r]^{\mathbb R}& {\cM}^\hH\ar[r]^{\mathbb U} &
{\cM}_L.
}
\end{equation}
{\color{blue}
The fruitful  approach to a Schneider type theorem for Hopf algebroids turns
out to be a study of these forgetful functors.}
\bigskip

The paper is organized as follows.
In Section \ref{sec:relsep} the notion of a
{\em separable} functor $\mathbb{U}$, {\em relative} to functors
$\mathbb{L}:{\mathfrak L}\to \mathfrak{A}$ and
$\mathbb{R}:{\mathfrak R}\to \mathfrak{A}$, is introduced and
investigated. Section \ref{sec:entw.str} concerns
relative separability of a forgetful functor $\cM^\cD\to \cM_L$,
associated to an entwining structure $(A,\cD,\psi)$ over an
algebra $L$. If $\cD$ possesses a grouplike element, relative
separability of the forgetful functor is shown to imply relative
injectivity of $A$ as a $\cD$-comodule and, in the case when in
addition the entwining map is bijective, also relative injectivity
of $A$ as an entwined module (see Theorem \ref{thm:rsentw} and
Proposition \ref{prop:cdrinj}).
In Section \ref{sec:hgd} separability of the forgetful
functor ${\mathbb U}:\cM^\hH\to \cM_L$, relative to the forgetful functor
${\mathbb R}$ from the category of relative Hopf modules to the category of
$\hH$-comodules, is studied,
for a Hopf algebroid $\hH$ and its comodule algebra $A$, cf. \eqref{eq:R&U}.  
In the case when the antipode of
$\hH$ is bijective, it is shown to be {\em equivalent} to relative
injectivity of the $\hH$-comodule $A$ {(see Theorem
\ref{thm:4.1.s})}. This result enables us to answer a question
posed in \cite{Bohm:hgdGal}. That is, {in Proposition
\ref{cor:fgp}} we prove that, in a Galois extension $B\subseteq A$
by a finitely generated and projective Hopf algebroid $\hH$ {with
a bijective antipode}, $A$ is faithfully flat as a left $B$-module
if and only if it is faithfully flat as a right $B$-module. The
main result is a Schneider type theorem in Section
\ref{sec:Schneider}.
Recall that Schneider's classical Theorem I in
\cite{Schn:PriHomS}
  deals with an algebra extension $B\subseteq A$ by a $k$-Hopf algebra $H$
  with a bijective antipode. It is assumed that {$H$} is a projective
  $k$-module
  and the canonical map \eqref{eq:can.Hopf} is surjective. Clearly, in this
  case the lifted canonical map \eqref{eq:can.lifted} is a split epimorphism
  of $k$-modules. As a proper generalization to an algebra extension
  $B\subseteq A$ by a Hopf algebroid $\hH$, in Theorem \ref{thm:schhgd} we
  assume that some lifted canonical map is a split epimorphism of modules for
  the (non-commutative) base algebra $L$ of $\hH$. This assumption is related
  to surjectivity of the canonical map and some projectivity conditions in
  Remark \ref{rem:tcan.split}. Under the assumption that, for an algebra
  extension $B\subseteq A$ by a Hopf algebroid $\hH$ with a bijective
  antipode, the lifted canonical map is a split epimorphism of $L$-modules,
the Galois property of the extension is related to relative
injectivity of the $\hH$-comodule $A$ and to the equivalence of
the category {$\cM^\hH_A$} of
relative Hopf modules to the category of
$B$-modules.  Section \ref{sec:eqproj} is devoted to a study of
(relative)
  equivariant injectivity and projectivity properties.
Preliminary results about entwining structures (over arbitrary non-commutative
{algebras}), coring extensions (in the sense of \cite{Brz:Corext})
and Hopf algebroids are collected in Appendix \ref{sec:prelims}.

{\color{blue}
As an experiment, using informality of the arXiv, in this submission
corrections (with respect to the first version) are written in blue. We hope
it would be helpful to the readers of the original submission.}

Throughout this paper the term {\em algebra} is used for an associative and
unital but not necessarily commutative algebra over a fixed commutative ring
$k$. Multiplication is denoted by juxtaposition and the unit element is
denoted by $1$.
For an algebra $A$, the opposite algebra is denoted by $A^{op}$.
The category of right (respectively, left) modules for an algebra $A$ is
denoted by $\cM_A$ (respectively, ${}_A \cM$). The set of morphisms between
two $A$-modules $M$ and $M'$ is denoted by $\Hom_A(M,M')$ (respectively,
${}_A\Hom(M,M')$). The category of $A$-$A$ bimodules is denoted by ${}_A\cM_A$
and its Hom sets by ${}_A \Hom_A(M,M')$.

For the coproduct in a coring $\cC$ over an algebra $A$, we use a Sweedler
type index notation $c\mapsto c^{(1)}\sstac A c^{(2)}$, for $c\in \cC$, where
implicit summation is understood. Similarly, for a right $\cC$-coaction
we use an index notation of the form $\varrho^M(m)=m^{[0]}\sstac A m^{[1]}$,
for $m\in M$. The category of right $\cC$-comodules is denoted by $\cM^\cC$
and its Hom sets by $\Hom^\cC(M,M')$. Symmetrical notations are used for left
$\cC$-comodules. The coaction is denoted by ${}^M\varrho(m)=m^{[-1]}\sstac A
m^{[0]}$, for a left $\cC$-comodule $M$ and $m\in M$. The category of left
$\cC$-comodules is denoted by ${}^\cC\cM$ and its Hom sets by
${}^\cC\Hom(M,M')$.

\section{Relative separable functors} \label{sec:relsep}

We start by recalling some material about separable functors. For
more information we refer to \cite[Chap. IX, page 307-312]{HS},
\cite[Chap. 8, page 279-281]{Weibel} and \cite{CMZ}. Throughout
the paper we use the following terminology. A morphism $f:C_1\to C_2$ in a
category ${\mathfrak C}$ is said to be a {\em split monomorphism}
or {\em section} if it is {\em
  cosplit} by  some morphism $h:C_2\to C_1$ in ${\mathfrak C}$, i.e. $h\circ
f=C_1$.
Dually, $f$ is called a {\em split epimorphism} or {\em retraction} provided
that it is {\em split} by some morphism $g:C_2\to C_1$ in ${\mathfrak C}$,
i.e. $f\circ g=C_2$.

\begin{definition}
\label{def 4.2.32} Let $ \mathfrak{C}$ be a
category and let $\mathcal{S}$ be a
class of morphisms in $\mathfrak{C}$.
For a morphism $f:C_{1}\rightarrow C_{2}$ in $\mathfrak{C}$,
an object $P\in \mathfrak{C}$ is called $f$-\emph{projective} if the map
$\Hom_\mathfrak{C}(P,f):\Hom_\mathfrak{C} (P,C_{1})\rightarrow
\Hom_\mathfrak{C}(P,C_{2})$ is surjective. $P$ is
$\mathcal{S}$\textbf{-}\emph{projective} if it is
$f$-projective for every $f\in \mathcal{S}$.

Dually, an object $I\in \mathfrak{C}$ is called $f$-\emph{injective} if it is
$f^{op}$-projective in the opposite category $\mathfrak{C}^{op}$, where
$f^{op}:C_{2}\rightarrow C_{1}$ is considered to be a morphism
$\mathfrak{C}^{op}$. $I$ is called $ \mathcal{S}$\textbf{-}\emph{injective }if
it is $f$-injective for every $ f\in \mathcal{S}$.
\end{definition}
All results below about projective objects can be dualized to get their
analogues for injective objects.
\begin{theorem}\cite{Ar2}\label{teo 4.2.33}
Let $\mathbb{H}:\mathfrak{B}\rightarrow
\mathfrak{A}$ be a covariant functor and consider a class of morphisms
\begin{equation}
\mathcal{E}_{\mathbb{H}}:=\{g\in \mathfrak{B}\mid
\mathbb{H}(g)\text{ is a split epimorphism in }\mathfrak{A}\}.
\label{proj class}
\end{equation}
Assume that $\mathbb{T}:\mathfrak{A}\rightarrow \mathfrak{B}$ is a left
adjoint of $\mathbb{H}$ and denote by $\varepsilon :\mathbb{TH}\rightarrow
{\mathfrak{B}}$ the counit of the adjunction.
Then, for an object $P\in \mathfrak{B}$, the following assertions are
equivalent.
\begin{itemize}
\item[$(a)$] $P$ is $\mathcal{E}_{\mathbb{H}}$-projective.
\item[$(b)$] $\varepsilon _{P}:\mathbb{TH}(P)\rightarrow P$ is a
split epimorphism.
\item[$(c)$] There is a split epimorphism $\pi :\mathbb{T}(X)\rightarrow P$,
for a suitable object $X\in \mathfrak{A}$.
\end{itemize}

In particular, all objects of the form $\mathbb{T}(X)$, for $X\in
\mathfrak{A}$, are $\mathcal{E}_{\mathbb{H}}$-projective.
\end{theorem}

In \cite{Ar2} also a dual version of Theorem \ref{teo 4.2.33} can be found.
It deals with ${\mathcal I}_{\mathbb T}$-injective objects in a category
${\mathfrak A}$, for a left adjoint functor ${\mathbb T}:{\mathfrak A}\to
{\mathfrak B}$, and
$$
{\mathcal I}_{\mathbb T}=\{ f\in {\mathfrak A} \mid
\mathbb{T}(f)\text{ is a split monomorphism in }\mathfrak{B}\}.
$$
Using the current terminology, relative injective right comodules
of an $A$-coring $\cC$, discussed in {Section} \ref{preli:coring},
can be characterized as $\mathcal{I}_{\mathbb{U}}$-injective
objects, where ${\mathbb U}:\cM^\cC\to \cM_A$ denotes the
forgetful functor.
  As recalled in {Section} \ref{preli:coring}, the forgetful functor
  ${\mathbb U}$   possesses a right adjoint, the functor
  $\bullet \sstac A \cC$. The unit of the adjunction is given by the
  $\cC$-coaction. Therefore,
 the dual version of Theorem \ref{teo 4.2.33} (a)$\Leftrightarrow$ (b)
includes the claim, recalled in {Section} \ref{preli:coring}, that
  a right $\cC$-comodule $M$ is relative injective if and only if the coaction
  $\varrho^M$ in it is a split monomorphism in $\cM^\cC$.

Since any covariant functor preserves split epimorphisms and split
monomorphisms, we immediately have that,
for any two functors
${\mathbb F}:\mathfrak{A}\rightarrow \mathfrak{B}$ and
${\mathbb G}:\mathfrak{B}\rightarrow \mathfrak{C}$,
\begin{equation}\label{eq:non.so}
\mathcal{E}_{\mathbb F}\subseteq \mathcal{E}_{{\mathbb G}{\mathbb F}}\qquad \
\text{and}\qquad
\mathcal{I}_{\mathbb F}\subseteq \mathcal{I}_{{\mathbb G}{\mathbb F}}.
\end{equation}

As explained in the Introduction, in the area of Schneider type
  theorems one often faces the following problem. Consider an entwining
  structure $(A,\cD,\psi)$ over an algebra $L$. Assume that some
  map {in  ${\cM}^{\cD}_A(\psi)$} (practically the canonical map) is a
  retraction in $\cM_L$. Under what
  assumptions is it a retraction also in $\cM^\cD$? Putting the question in
  a more functorial way, we can ask in which cases is
  ${\mathcal E}_{\mathbb F}={\mathcal E}_{\mathbb{GF}}$, for the
  forgetful functors ${\mathbb F}: {\cM}^{\cD}_A(\psi)\to \cM^\cD$ and
  ${\mathbb   G}: \cM^\cD \to \cM_L$.
  For these particular functors ${\mathbb F}$ and ${\mathbb G}$, property 1)
  in Proposition \ref{prop:Es.equal}  below reduces to a similar (but somewhat
  weaker) assumption as in a Schneider type theorem \cite[Theorem
  5.9]{SchSch:genGal} (see also \cite[Theorem 4.6]{Brz:Galcom}).
  Properties like in part 2) of Proposition \ref{prop:Es.equal}
  are assumed e.g. in \cite[Corollary 4.8]{SchSch:genGal}.

\begin{proposition}\label{prop:Es.equal}
For two functors $\mathbb{F}:\mathfrak{A}\to \mathfrak{B}$ and
$\mathbb{G}:\mathfrak{B}\to \mathfrak{C}$, {
$\mathcal{E}_{\mathbb{F}}= \mathcal{E}_{\mathbb {G}\mathbb{F}}$}
whenever any of the following properties hold.
\begin{itemize}
\item[$1)$] $\mathbb{F}\left( A\right) $ is $\mathcal{E}_{\mathbb
  G}$-projective, for every object $A\in \mathfrak{A}$.
\item[$2)$] $\mathfrak{A,B}$ and $\mathfrak{C}$ are abelian
categories, ${\mathbb G}$ is left exact and reflects epimorphisms,
$\mathbb{F}$ is left
  exact and $\mathbb{F} \left( A\right) $ is
$\mathcal{I}_{\mathbb G}$-injective, for every object $A\in \mathfrak{A}$.
\end{itemize}
Dually, { $\mathcal{I}_{\mathbb{F}}= \mathcal{I}_{{\mathbb
G}\mathbb{F}}$} whenever any of the following properties hold.
\begin{itemize}
\item[$1^{op})$] $\mathbb{F}\left( A\right) $ is
$\mathcal{I}_{\mathbb G}$-injective, for every object $A\in
\mathfrak{A}$. \item[$2^{op})$] $\mathfrak{A,B}$ and
$\mathfrak{C}$ are abelian categories, $\mathbb{G}$ is right exact
and reflects monomorphisms, $\mathbb{F}$ is right exact and
$\mathbb{F}\left( A\right) $ is $\mathcal{E}_{\mathbb
G}$-projective, for every object $A\in\mathfrak{A}$.
\end{itemize}
\end{proposition}
\begin{proof}
$1)$ Let $f:A_1\rightarrow A_2$ be a morphism in
  $\mathcal{E}_{\mathbb{G} \mathbb{F}}$.
Then $\mathbb{F}\left( f\right):\mathbb{F}\left( A_1\right)
\rightarrow \mathbb{F}\left( A_2\right)$ belongs to
$\mathcal{E}_{\mathbb G}$ and hence, by hypothesis, it is a split
  epimorphism. Thus $f\in \mathcal{E}_{\mathbb{F}}.$

$2)$ For $f\in \mathcal{E}_{\mathbb{G}\mathbb{F}}$,
consider the exact sequence (kernel diagram)
\begin{equation*}
0\rightarrow K\overset{i}{\longrightarrow }A_1\overset{f}{\longrightarrow }A_2
\end{equation*}
in ${\mathfrak A}$.
The left exact functor $\mathbb{F}$ takes it to the exact sequence
$$
\xymatrix{
0\ar[r]&
\mathbb{F}\left( K\right)\ar[r]^{\mathbb{F}\left( i\right) }&
\mathbb{F}\left( A_1\right)\ar[r]^{\mathbb{F}\left( f\right) }&
\mathbb{F}\left( A_2\right)
}
$$
in ${\mathfrak B}$.
Since $f$ is an element of $\mathcal{E}_{\mathbb{G}\mathbb{F}}$, the morphism
$\mathbb{G}\mathbb{F}\left(f\right) $ is a split epimorphism.
Since $\mathbb{G}$ is left exact and $\mathfrak{C}$ is an abelian
category, the sequence
$$
\xymatrix{
0\ar[r]&
\mathbb{G}\mathbb{F}\left( K\right)\ar[r]^{\mathbb{G}\mathbb{F}\left( i\right)
}&
\mathbb{G}\mathbb{F}\left( A_1\right)\ar[r]^{\mathbb{G}\mathbb{F}
\left( f\right) }&
\mathbb{G}\mathbb{F}\left( A_2\right)\ar[r]&
0
}
$$
in $\mathfrak{C}$ is split exact. Thus we deduce that $i\in
\mathcal{I}_{\mathbb{G}\mathbb{F}}$. Moreover, since $\mathbb{G}$
reflects epimorphisms, $\mathbb{F}\left(f\right) $ is an
epimorphism. So the sequence
$$
\xymatrix{
0\ar[r]&
\mathbb{F}\left( K\right)\ar[r]^{\mathbb{F}\left( i\right) }&
\mathbb{F}\left( A_1 \right)\ar[r]^{\mathbb{F}\left( f\right) }&
\mathbb{F}\left( A_2 \right)\ar[r]&
0
}
$$
in ${\mathfrak B}$ is exact too. Since $i$ is an element of
$\mathcal{I}_{\mathbb{G}\mathbb{F}}$, its image ${\mathbb F}(i)$ is in
$\mathcal{I}_\mathbb{G}$. By assumption $\mathbb{F}\left( K\right) $ is
$\mathcal{I}_{\mathbb{G}}$-injective hence
the monomorphism $\mathbb{F}\left( i\right) $ is split. Since
$\mathfrak{B}$ is an abelian category, we conclude that
$\mathbb{F}\left( f\right) $ is a split epimorphism, i.e. that
$f\in \mathcal{E}_{\mathbb{F}}$.

Claims $1^{op})$ and $2^{op})$ follow by duality.
\end{proof}

The {most important} notions of this section are introduced in
the following definition.

\begin{definition}\label{def:rsep}
Consider the following diagram of functors.
\begin{equation*}\xymatrix@C=0.5cm{&&\mathfrak{B}   \\
   \mathfrak{L} \ar[rr]^{\mathbb{L}} && \mathfrak{A}\ar[u]^{\mathbb{U}}  &&
   \mathfrak{R}\ar[ll]_{\mathbb{R}} }
\end{equation*}They give rise
to two functors
\begin{equation*}
\Hom_{\mathfrak{A}}(\mathbb{L}\left( \bullet \right), \mathbb{R}\left(
\bullet \right) )\quad \textrm{and}\quad
\Hom_{\mathfrak{B}}({\mathbb{U}}\mathbb{L}(\bullet ),{
\mathbb{U}}\mathbb{R}(\bullet )):\mathfrak{L}^{op}\times \mathfrak{R}
\rightarrow \underline{\underline{\mathfrak{Sets}}}
\end{equation*}
and a natural transformation between them
\begin{equation}\label{eq:nat.tr.U}
\Phi(\mathbb{U},\mathbb{L},\mathbb{R}):
\Hom_{\mathfrak{A}}(\mathbb{L}\left(\bullet \right) ,\mathbb{R}
\left( \bullet \right) )\rightarrow
\Hom_{\mathfrak{B}}({\mathbb{U}}\mathbb{L}(\bullet ),{\mathbb{U}}
\mathbb{R}(\bullet)),\qquad
\Phi(\mathbb{U},\mathbb{L},\mathbb{R})_{{L},{R}}({f}):={\mathbb{U}}({f}),
\end{equation}
for all objects ${L}\in \mathfrak{L},{R}\in \mathfrak{R}$ and for
every morphism $f:\mathbb{L}(L)\to \mathbb{R}(R)$. We say that
\begin{itemize}
\item[1)] ${\mathbb{U}}$ is \emph{$\left(
\mathbb{L},\mathbb{R}\right)
  $-faithful} if
$\Phi(\mathbb{U},\mathbb{L},\mathbb{R})_{L,R}$ is injective, for
every objects ${L}\in \mathfrak{L}$ and ${R}\in \mathfrak{R}$.
\item[2)] ${\mathbb{U}}$ is \emph{$\left(
\mathbb{L},\mathbb{R}\right)
  $-full} if
$\Phi(\mathbb{U},\mathbb{L},\mathbb{R})_{L,R}$ is surjective, for
every objects ${L}\in \mathfrak{L}$ and ${R}\in  \mathfrak{R}$.
\item[3)] ${\mathbb{U}}$ is \emph{$\left(\mathbb{L},\mathbb{R}\right)
  $-separable} if
$\Phi(\mathbb{U},\mathbb{L},\mathbb{R})$
{is a split natural monomorphism}. \item[4)] ${\mathbb{U}}$ is
\emph{$\left( \mathbb{L},\mathbb{R}\right) $-coseparable} if
$\Phi(\mathbb{U},\mathbb{L},\mathbb{R})$
is a split natural epimorphism.
\end{itemize}
\end{definition}
When both $\mathbb{L}$ and $\mathbb{R}$ are identity functors, we
recover the classical definitions of a faithful, full, separable
and naturally full (here called coseparable) functor. We are
particularly interested in the case when either $\mathbb{L}$ or
$\mathbb{R}$ is the identity functor. Anyway, some of our results
can be stated for the general case.

\begin{remark}
\label{claim: Rafael}
Following \cite[page 1446]{Rafael}, one can
prove that Definition \ref{def:rsep} 3) can be reformulated (in the spirit of
a characterization of separable functors in \cite{NdO}) as follows.
A functor $\mathbb{U}:\mathfrak{A}\rightarrow \mathfrak{B} $ is
$\left({\mathbb{L}},{\mathbb{R}}\right) $-separable, for some
functors $\mathbb{L}: \mathfrak{L}\rightarrow \mathfrak{A}$ and
$\mathbb{R}:\mathfrak{R} \rightarrow \mathfrak{A}$,
if and only if there is a map
\begin{equation*}
{\widetilde {\Phi}}(\mathbb{U},\mathbb{L},\mathbb{R})_{L,R}:
\Hom_{\mathfrak{B}}({\mathbb{UL}}(L),{\mathbb{UR}}(R))\rightarrow
\Hom_{\mathfrak{A}}({\mathbb{L}}(L),{\mathbb{R}}(R)),
\end{equation*}
for all objects ${L}\in \mathfrak{L}$ and ${R}\in \mathfrak{R}$, satisfying the
following identities.
\begin{itemize}
\item[$S1)$] $
\widetilde{\Phi}(\mathbb{U},\mathbb{L},\mathbb{R})_{L,R}
({\mathbb{U}}({f}))={f}$, for any ${f}\in
  \Hom_{\mathfrak{A}}({\mathbb{L}}(L),{\mathbb{R}}(R))$.
\item[$S2)$] $
\widetilde{\Phi}(\mathbb{U},\mathbb{L},\mathbb{R})_{L',R'}(h')\circ
{\mathbb{L}}({l})={\mathbb{R}}({r})\circ
\widetilde{\Phi}(\mathbb{U},\mathbb{L},\mathbb{R})_{L,R}
(h)$, for every commutative diagram in $\mathfrak{B}$ of the following form.
\begin{equation*}
\xymatrix@R=15pt@C=30pt{
  \mathbb{UL}(L) \ar[d]_{\mathbb{UL}(l)} \ar[r]^{h} & \mathbb{UR}(R)
  \ar[d]^{\mathbb{UR}(r)} \\
  \mathbb{UL}(L') \ar[r]_{h'} & \mathbb{UR}(R')  }
\end{equation*}
\end{itemize}
\end{remark}
\begin{remark}
\label{rem faithful} Recall that faithful functors reflect mono, and
  epimorphisms. Analogously, for an $\left( {\mathbb{L}},{\mathbb{R}} \right)
  $-faithful functor $\mathbb{U}$ the following hold true.
\begin{itemize}
\item[$1)$] Assume that ${\mathbb{R}}$ is surjective on the objects
and let $f:A\rightarrow {\mathbb{L}}(L)$ be a morphism in
$\mathfrak{A}$. Then $f$ is an epimorphism whenever $\mathbb{U}(f)$ is.
\item[$2)$] Assume that ${\mathbb{L}}$ is surjective on the objects and let
  $f:{\mathbb{R}}(R)\rightarrow A$ be a morphism in $\mathfrak{A}$.
Then $f$ is a monomorphism whenever $\mathbb{U}(f)$ is.
\end{itemize}
\end{remark}

In the rest of the section we extend some standard results about
  separable functors to relative separable functors in Definition
  \ref{def:rsep} 3). Analogous results can be obtained for coseparable
  functors by a careful dualization.
\begin{theorem}\label{teo compos of separable}
Consider the following diagram of functors.
\begin{equation*}\xymatrix@R=15pt@C=30pt{&&\mathfrak{C}   \\
&&\mathfrak{B}\ar[u]^{\mathbb{V}}    \\
  \mathfrak{L'} \ar[r]^{\mathbb{L'}} & \mathfrak{L} \ar[r]^{\mathbb{L}} &
  \mathfrak{A}\ar[u]^{\mathbb{U}}  & \mathfrak{R}\ar[l]_{\mathbb{R}}  &
  \mathfrak{R'}\ar[l]_{\mathbb{R'}}}
\end{equation*}
The following assertions hold true.
\begin{itemize}
\item[$1)$] If $\mathbb{U}$ is $(\mathbb{L},\mathbb{R})$-separable
then
  $\mathbb{U}$ is $(\mathbb{LL}^{\prime },\mathbb{RR}^{\prime })$-separable.
\item[$2)$] If $\mathbb{U}$ is $(\mathbb{L},\mathbb{R})$-separable
and
  $\mathbb{V}$ is $(\mathbb{U}\mathbb{L},\mathbb{U}\mathbb{R})$-separable
  then $\mathbb{V}\mathbb{U}$ is $(\mathbb{L},\mathbb{R})$-separable.
\item[$3)$] If $\mathbb{V}\mathbb{U}$ is
$(\mathbb{L},\mathbb{R})$-separable then $\mathbb{U}$ is
$(\mathbb{L},\mathbb{R})$-separable.
\end{itemize}
\end{theorem}

\begin{proof}
The proof is similar to \cite[I.3 Proposition 46 and Corollary 9]{CMZ}.

$1)$ Since $\mathbb{U}$ is $(\mathbb{L},\mathbb{R})$-separable,
there exists a natural retraction
$\widetilde{\Phi}(\mathbb{U},\mathbb{L},\mathbb{R})$
of the natural transformation \eqref{eq:nat.tr.U}.
For any objects $L^{\prime }\in \mathfrak{L}^{\prime}$ and $R^{\prime
}\in \mathfrak{R}^{\prime }$, the maps
\begin{equation*}
\widetilde{\Phi}(\mathbb{U},\mathbb{L},\mathbb{R})
_{\mathbb{L}^{\prime }(L^{\prime }),\mathbb{R}^{\prime}(R^{\prime })}
:\Hom_{\mathfrak{B}}(\mathbb{U}{\mathbb{L}}\mathbb{L}^{\prime
}(L^{\prime }),\mathbb{U}{\mathbb{R}}\mathbb{R}^{\prime }(R^{\prime
}))\rightarrow
\Hom_{\mathfrak{A}}({\mathbb{L}}\mathbb{L}^{\prime }(L^{\prime }),{\mathbb{R}}
\mathbb{R}^{\prime }(R^{\prime }))
\end{equation*}
define a natural transformation which is a retraction of
$\Phi(\mathbb{U},\mathbb{L}\mathbb{L}',\mathbb{R}\mathbb{R}')$, defined
analogously to \eqref{eq:nat.tr.U}.

$2)$ The natural transformation
\begin{equation}\label{eq:nat.tr.VU}
\Phi(\mathbb{V}\mathbb{U},\mathbb{L},\mathbb{R}):
\Hom_{\mathfrak{A}}({\mathbb{L}}(\bullet ),{\mathbb{R}}(\bullet
))\rightarrow \Hom_{\mathfrak{C}}(\mathbb{V}\mathbb{U}{\mathbb{L}}(\bullet
),\mathbb{V}\mathbb{U}{\mathbb{R}}(\bullet )),\quad f\mapsto
\mathbb{V}\mathbb{U}\left( f\right)
\end{equation}
is a composite of the split natural monomorphisms $\Phi(\mathbb{U},
\mathbb{L},\mathbb{R})$, corresponding via \eqref{eq:nat.tr.U} to the
$(\mathbb{L},\mathbb{R})$-separable functor $\mathbb{U}$, and
$\Phi(\mathbb{V},\mathbb{U}\mathbb{L},\mathbb{U}\mathbb{R})$, corresponding to
the $(\mathbb{U}\mathbb{L},\mathbb{U}\mathbb{R})$-separable
functor $\mathbb{V}$. Hence \eqref{eq:nat.tr.VU} is a split natural
monomorphism too,
proving $(\mathbb{L},\mathbb{R})$-separability of $\mathbb{V}\mathbb{U}$.

$3)$ Since the functor $\mathbb{V}\mathbb{U}$ is
$(\mathbb{L},\mathbb{R})$-separable, the corresponding natural
transformation \eqref{eq:nat.tr.U} possesses a retraction
$\widetilde{\Phi}(\mathbb{V} \mathbb{U},\mathbb{L},\mathbb{R})$.
The composite $\widetilde{\Phi}(\mathbb{V}
\mathbb{U},\mathbb{L},\mathbb{R})\circ
\Phi(\mathbb{V},\mathbb{U}\mathbb{L}, \mathbb{U}\mathbb{R})$ is a
natural retraction of $\Phi(\mathbb{U},\mathbb{L},\mathbb{R})$ in
\eqref{eq:nat.tr.U}.
\end{proof}

\begin{theorem}[Maschke type Theorem]\label{thm:rmaschke}
Let ${\mathbb{U}}:\mathfrak{A}\rightarrow \mathfrak{B}$,
$\mathbb{L}:\mathfrak{L}\rightarrow \mathfrak{A}$ and
$\mathbb{R}:\mathfrak{R}\rightarrow \mathfrak{A}$ be functors.
\begin{itemize}
\item[1)] If ${\mathbb{U}}$ is
  $\left({\mathfrak{A}},\mathbb{R}\right)$-separable then,
for any objects $R\in \mathfrak{R}$ and $A\in \mathfrak{A}$, a morphism
$f:\mathbb{R}\left( R\right) \rightarrow A$ is a split
monomorphism whenever ${\mathbb{U}}\left( f\right) $ is a split
monomorphism. Moreover, in this case
$\mathcal{I}_{\mathbb{R}}=\mathcal{I}_{{\mathbb{U}}\mathbb{R}}$
and
$\mathcal{E}_{\mathbb{R}}=\mathcal{E}_{{\mathbb{U}}\mathbb{R}}$.
\item[2)] If ${\mathbb{U}}$ is $\left(
\mathbb{L},{\mathfrak{A}}\right)$-separable then, for any objects
$L\in \mathfrak{L}$ and $A\in \mathfrak{A}$, a morphism
$f:A\rightarrow \mathbb{L}\left( L\right) $ is a split epimorphism
whenever ${\mathbb{U}}\left( f\right) $ is a split epimorphism.
Moreover, in this case
$\mathcal{I}_{\mathbb{L}}=\mathcal{I}_{{\mathbb{U}}\mathbb{L}}$
and
$\mathcal{E}_{\mathbb{L}}=\mathcal{E}_{{\mathbb{U}}\mathbb{L}}.$
\end{itemize}
\end{theorem}

\begin{proof}
Let $A$, $R$ and $f$ be as in part 1).
Let $\widetilde{\Phi}(\mathbb{U},{ \mathfrak{A}},\mathbb{R})$ be a
natural retraction of ${\Phi}(\mathbb{U},{
\mathfrak{A}},\mathbb{R})$ in \eqref{eq:nat.tr.U}. In view of
$S2)$ in Remark \ref{claim: Rafael}, any retraction $\pi$
of ${\mathbb{U}}\left( f\right)$ satisfies
\begin{equation*}
\widetilde{\Phi}(\mathbb{U},\mathfrak{A},\mathbb{R})_{A,R}
\left( \pi \right) \circ f={\mathbb{R}\left( R\right) }.
\end{equation*}
That is, $f$ is a split monomorphism. In particular, $f:=\mathbb{R}(g)$
is a split monomorphism, for any $g\in
\mathcal{I}_{\mathbb{U}\mathbb{R}}$. Together with \eqref{eq:non.so}
this proves $\mathcal{I}_{{\mathbb{U}}\mathbb{R}}=\mathcal{I}_{\mathbb{R}}.$
Next take a morphism $g:R\rightarrow R^{\prime}$ in
$\mathcal{E}_{{\mathbb{U}}\mathbb{R}}$, and a section $\sigma$
of ${\mathbb{U}} \mathbb{R}(g)$.
Then, by naturality of
$\widetilde{\Phi}(\mathbb{U},\mathfrak{A},\mathbb{R})$,
\begin{equation*}
{\mathbb{R}\left( R^{\prime }\right) }
=\widetilde{\Phi}(\mathbb{U},\mathfrak{A},\mathbb{R})_{\mathbb{R}\left(
  R^{\prime }\right), R^{\prime }}\left( {\mathbb{U}\mathbb{R}\left( R^{\prime
  }\right) }\right)
=\widetilde{\Phi}(\mathbb{U},\mathfrak{A},\mathbb{R})_{\mathbb{R}\left(
  R^{\prime  }\right),R^{\prime }}\big( {\mathbb{U}}\mathbb{R}\left( g\right)
  \circ {\mathbb{\sigma }}\big)
=\mathbb{R}\left( g\right) \circ
  \widetilde{\Phi}(\mathbb{U},\mathfrak{A},\mathbb{R}) _{\mathbb{R}\left(
    R^{\prime }\right) ,R}\left( {\mathbb{\sigma }}\right).
\end{equation*}
This implies that $\mathbb{R}\left( g\right) $ is a split epimorphism, i.e.
$g\in \mathcal{E}_{\mathbb{R}}.$
In view of \eqref{eq:non.so}, we have
$\mathcal{E}_{\mathbb{R}}=\mathcal{E}_{{\mathbb{U}}\mathbb{R}}$
proven.

Part 2) is proven by dual reasoning.
\end{proof}

\begin{corollary}\label{coro Rafael}
Let $(\mathbb{T},\mathbb{H})$ be an adjunction of functors
$\mathbb{T}:\mathfrak{A}\rightarrow \mathfrak{B}$ and
$\mathbb{H}:\mathfrak{B} \rightarrow \mathfrak{A}$. For any
functors $\mathbb{L}:\mathfrak{L}\rightarrow \mathfrak{B}$ and
$\mathbb{R}:\mathfrak{R}\rightarrow \mathfrak{A}$, the following
hold.
\begin{itemize}
\item[$1)$] If the functor $\mathbb{H}$ is
$(\mathbb{L},{\mathfrak{B}})$-separable then $\mathbb{L}\left(
L\right)$ is
  $\mathcal{E}_{\mathbb{H}}$-projective for every $L\in \mathfrak{L}$.
\item[$2)$] If the functor $\mathbb{T}$ is
$({\mathfrak{A}},\mathbb{R})$-separable then $\mathbb{R}\left(
R\right)$ is
  $\mathcal{I}_{\mathbb{T}}$-injective for every $R\in \mathfrak{R}$.
\end{itemize}
\end{corollary}

\begin{proof} Let $\eta:{\mathfrak A}\to \mathbb{H}\mathbb{T}$ be the unit
 and $\varepsilon:\mathbb{T}\mathbb{H}\to {\mathfrak B}$ be the counit
 of the adjunction $(\mathbb{T},\mathbb{H})$.

 1) For any object $L\in \mathfrak{L}$, the epimorphism
  $\mathbb{H}(\varepsilon _{\mathbb{L}\left(L\right) })$ is split by
$\eta _{\mathbb{HL}\left( L\right) }$. Hence,
by Theorem \ref{thm:rmaschke} 2), $\varepsilon _{\mathbb{L}\left( L\right) }$
{is a split epimorphism} in $\mathfrak{B}$. By Theorem \ref{teo
4.2.33} $(b)\Rightarrow (a)$, $\mathbb{L}\left( L\right) $ is
$\mathcal{E}_{\mathbb{H}}$-projective.

2) For any object $R\in {\mathfrak R}$, the monomorphism
$\mathbb{T}(\eta _{\mathbb{R}\left( R\right)})$ is split by
$\varepsilon_{\mathbb{TR}\left(R\right)}$. Hence the claim follows
analogously to part 1), by Theorem \ref{thm:rmaschke} 1) and
a dual form of Theorem \ref{teo 4.2.33}.
\end{proof}
In the following theorem functors
preserving and reflecting relative projective (resp. injective) objects are
studied.
\begin{theorem}
\label{teo F and P-project} Let $(\mathbb{T},\mathbb{H})$ and $(\mathbb{
T^{\prime }},\mathbb{H^{\prime }})$ be adjunctions and consider
the following (not necessarily commutative) diagrams of functors.
\begin{equation*}
\xymatrix@R=15pt@C=50pt{
  \mathfrak{A} \ar[d]_{\mathbb{T}} \ar[r]^{\mathbb{F}'} & \mathfrak{A'}
  \ar[d]^{\mathbb{T}'} \\
  \mathfrak{B} \ar[r]_{\mathbb{F}} & \mathfrak{B'}   }
\text{ \qquad } \xymatrix@R=15pt@C=50pt{
  \mathfrak{A} \ar[r]^{\mathbb{F}'} & \mathfrak{A'}  \\
  \mathfrak{B} \ar[u]_{\mathbb{H}} \ar[r]_{\mathbb{F}} & \mathfrak{B'}
  \ar[u]^{\mathbb{H}'}  }
\end{equation*}
If ${\mathbb{T}^{\prime }}\mathbb{F}^{\prime }$ and
$\mathbb{F}{\mathbb{T}}$ are naturally equivalent, then the following
hold.
\begin{itemize}
\item[$1)$] If an object $P$ in $\mathfrak{B}$ is
  $\mathcal{E}_{\mathbb{H}}$-projective then $\mathbb{F}(P)$ is
  $\mathcal{E}_{\mathbb{H}^{\prime }}$-projective.
\item[$2^{op}$)] Assume that $\mathbb{F}^{\prime }$ is
  $({\mathfrak{A}},\mathbb{R})$-separable for some functor
  $\mathbb{R}:\mathfrak{R} \rightarrow \mathfrak{A}$.
  If, for an object $R\in \mathfrak{R}$,
  the object $\mathbb{F}^{\prime }{\mathbb{R}}\left( R\right)$ is
  $\mathcal{I}_{\mathbb{T}^{\prime }}$-injective, then ${\mathbb{R}}\left(
  R\right)$ is $\mathcal{I}_{\mathbb{T}}$-injective.
\end{itemize}
If $\mathbb{F}^{\prime }{\mathbb{H}}$ and ${\mathbb{H}^{\prime
}}\mathbb{F}$ are naturally equivalent, then the following hold.
\begin{itemize}
\item[$1^{op}$)] If an object $I$ in $\mathfrak{A}$ is
  $\mathcal{I}_{\mathbb{T}}$-injective, then $\mathbb{F}^{\prime }(I)$ is
  $\mathcal{I}_{\mathbb{T}^{\prime}}$-injective.
\item[$2)$] Assume that $\mathbb{F}$ is
  $(\mathbb{L},{\mathfrak{B}})$-separable for some functor
  $\mathbb{L}:\mathfrak{L}\rightarrow \mathfrak{B}$.
  If, for an object $L\in \mathfrak{L}$, the object
  $\mathbb{F}{\mathbb{L}}\left( L\right) $
  is $\mathcal{E}_{\mathbb{H}^{\prime }}$-projective then ${\mathbb{L}}\left(
  L\right)$ is   $\mathcal{E}_{\mathbb{H}}$-projective.
\end{itemize}
\end{theorem}
\begin{proof}
Denote by $\eta:{\mathfrak{A}}\rightarrow \mathbb{H}\mathbb{T}$
the unit and by $\varepsilon :\mathbb{T}\mathbb{H}\rightarrow {\mathfrak{B}}$
the counit of the adjunction $(\mathbb{T},\mathbb{H})$.

1) By Theorem \ref{teo 4.2.33} (a) $\Rightarrow$ (b),
$\mathcal{E}_{\mathbb{H}}$-projectivity of $P$ implies that
$\varepsilon _{P}:\mathbb{T}\mathbb{H}(P)\rightarrow P$
{is a split epimorphism}. Hence also $\mathbb{F}(\varepsilon
_{P}):\mathbb{T}^{\prime} \mathbb{F}^{\prime}
\mathbb{H}(P)\sim\mathbb{F}\mathbb{T}\mathbb{H}(P) \rightarrow
\mathbb{F}(P)$ is a split epimorphism. Application of Theorem
\ref{teo 4.2.33} (c) $\Rightarrow$ (a) to the adjunction
$(\mathbb{T}^{\prime},\mathbb{H}^{\prime })$ completes the proof
of $\mathcal{E}_{\mathbb{H}^{\prime }}$-projectivity of
$\mathbb{F}(P)$.

2) For any object $P$ in $\mathfrak{B}$, the monomorphism $\eta
_{\mathbb{H}(P)}$ is split by $\mathbb{H}(\varepsilon_{P})$. Hence
$\mathbb{F}^{\prime }\mathbb{H} (\varepsilon _{P})$ {is a split
epimorphism}.
By the natural equivalence
$\mathbb{F}^{\prime }{\mathbb{H}}\sim{\mathbb{H}^{\prime
}}\mathbb{F}$, also $\mathbb{H}^{\prime }\mathbb{F}(\varepsilon
_{P})$
{is a split epimorphism}, yielding that $\mathbb{F}(\varepsilon
_{P}):\mathbb{F}\mathbb{T}\mathbb{H}(P)\rightarrow \mathbb{F}(P)$
belongs to $\mathcal{E}_{\mathbb{H}^{\prime }}$.
In the case when
$\mathbb{F}(P)$ is $\mathcal{E}_{\mathbb{H}^{\prime
}}$-projective, we conclude that $\mathbb{F}(\varepsilon_{P})$ is a split
epimorphism in $\mathfrak{B}'$. Now put
$P={\mathbb{L}}\left( L\right)$, such that ${\mathbb{FL}}\left(
L\right)$ is $\mathcal{E}_{\mathbb{H}^{\prime }}$-projective as in
the claim. Then, by Theorem \ref{thm:rmaschke}$~ 2), \varepsilon
_{{\mathbb{L}}\left( L\right)}$
{is a split epimorphism} in $\mathfrak{B}$ and hence
${\mathbb{L}}\left( L\right)$ is $\mathcal{
E}_{\mathbb{H}}$-projective by Theorem \ref{teo 4.2.33} (b)
$\Rightarrow$ (a).

The remaining claims $1^{op}$) and $2^{op}$) follow by dual
reasoning.
\end{proof}
Let $(\mathbb{T},\mathbb{H})$ be an adjunction of functors ${\mathbb{T}}:
\mathfrak{A}\rightarrow \mathfrak{B}$ and ${\mathbb{H}}:
\mathfrak{B}\rightarrow \mathfrak{A}$. Denote by $\varepsilon :{\mathbb{T}}
\mathbb{H}\rightarrow {\mathfrak{B}}$ and $\eta:{\mathfrak{A}}\rightarrow
\mathbb{H}{\mathbb{T}}$ the counit and the unit of the adjunction,
respectively. Consider the canonical isomorphism
\begin{equation}\label{eq:isophi}
\phi _{A,B}:\Hom_{\mathfrak{B}}(\mathbb{T}\left( A\right) ,B)\rightarrow
\Hom_{\mathfrak{A}}(A,\mathbb{H}\left( B\right) ),\qquad\phi _{A,B}\left(
f\right)
=\mathbb{H}\left( f\right) \circ \eta _{A}
\end{equation}
with inverse
\begin{equation*}
\phi _{A,B}^{-1}
:\Hom_{\mathfrak{A}}(A,\mathbb{H}\left( B\right) )\rightarrow
\Hom_{\mathfrak{B}}(\mathbb{T}\left( A\right) ,B),\qquad
\phi _{A,B}^{-1}\left(g\right) =\varepsilon _{B}\circ
\mathbb{T}\left(g\right).
\end{equation*}
In terms of the natural
transformations \eqref{eq:nat.tr.U} and \eqref{eq:isophi}, for any functors
$\mathbb{L}:\mathfrak{L}\rightarrow \mathfrak{A}$ and
$\mathbb{R}:\mathfrak{R}\rightarrow \mathfrak{A}$, define a natural
transformation
$$
\Omega  :=\phi \circ\Phi(\mathbb{T},\mathbb{L},\mathbb{R}):
\Hom_{\mathfrak{A}}(\mathbb{L}\left(
\bullet \right) ,\mathbb{R}\left( \bullet \right) )\rightarrow \Hom_{
\mathfrak{A}}(\mathbb{L}(\bullet ),\mathbb{H{T}R}(\bullet )).
$$
Then, for every morphism $f:\mathbb{L}\left( L\right) \rightarrow
\mathbb{R}\left(
R\right) ,$ one has
\begin{equation}\label{eq:Omega}
\Omega _{L,R}\left( f\right)  =\mathbb{H}{\mathbb{T}}\left(
f\right) \circ \eta _{\mathbb{L}\left( L\right) }=\eta
_{\mathbb{R}\left( R\right) }\circ f.
\end{equation}
Dually, for functors $\mathbb{L}:\mathfrak{L}\rightarrow
\mathfrak{B}$ and $\mathbb{R}:\mathfrak{R}\rightarrow
\mathfrak{B}$, there is a natural transformation
$$
\mho  :=\phi^{-1} \circ\Phi(\mathbb{H},\mathbb{L},\mathbb{R}):
\Hom_{\mathfrak{A}}(\mathbb{L}\left( \bullet
\right) ,\mathbb{R}\left( \bullet \right) )\rightarrow \Hom_{\mathfrak{A}}({
\mathbb{T}}\mathbb{HL}(\bullet ),\mathbb{R}(\bullet )),
$$
mapping a morphism $f:\mathbb{L}(L)\to \mathbb{R}(R)$ to
$$
\mho _{L,R}\left( f\right)  =\varepsilon _{\mathbb{R}\left(
R\right)
}\circ {\mathbb{T}}\mathbb{H}\left( f\right) =f\circ \varepsilon _{\mathbb{L}
\left( L\right) }.
$$
\begin{lemma}\label{lem:ff}
On the category of functors and natural transformations consider the following
endofunctor $\alpha$. It maps a functor ${\mathbb F}:{\mathfrak A}\to
{\mathfrak B}$ to the functor $\mathrm{Hom}_{\mathfrak B}(\bullet, {\mathbb
  F}(\bullet)):{\mathfrak B}^{op}\times {\mathfrak A}\to
\underline{\underline{\mathfrak{Sets}}}$,
and it maps a natural transformation $\sigma\in
\underline{\underline{\mathrm{Nat}}}({\mathbb F},{\mathbb G})$
to $\mathrm{Hom}_{\mathfrak B}(\bullet,\sigma_\bullet)$, i.e.
$$
\alpha(\sigma)_{B,A}:
\mathrm{Hom}_{\mathfrak B}(B,{\mathbb F}(A))\to
\mathrm{Hom}_{\mathfrak B}(B,{\mathbb G}(A)),\qquad
g\ \mapsto \ \sigma_A\circ g.
$$
The functor $\alpha$ is fully faithful.
\end{lemma}
\begin{proof}
The bijectivity of the maps
$$
\alpha_{{\mathbb F},{\mathbb G}}:
\underline{\underline{\mathrm{Nat}}}({\mathbb F},{\mathbb G})\to
\underline{\underline{\mathrm{Nat}}}\big(\mathrm{Hom}_{\mathfrak B}
(\bullet, {\mathbb F}(\bullet)),
\mathrm{Hom}_{\mathfrak B}(\bullet,{\mathbb G}(\bullet))\big)\qquad
\sigma\mapsto \alpha(\sigma),
$$
for any functors ${\mathbb F},{\mathbb G}:{\mathfrak A}\to {\mathfrak B}$,
is proven by constructing the inverse
$(\alpha_{{\mathbb F},{\mathbb G}})^{-1}({\mathcal P})_A:=
{\mathcal P}_{{\mathbb F}(A),A}({\mathbb F}(A))$,
for ${\mathcal P}\in
\underline{\underline{\mathrm{Nat}}}\big(\mathrm{Hom}_{\mathfrak B}
(\bullet, {\mathbb F}(\bullet)), \mathrm{Hom}_{\mathfrak B}
(\bullet, {\mathbb G}(\bullet))\big)$ and $A\in {\mathfrak{A}}$.
It is straightforward to check that the naturality of ${\mathcal P}$ (i.e. the
identity ${\mathbb G}(a)\circ {\mathcal P}_{B,A}(g)\circ b=
{\mathcal P}_{B',A'}\big({\mathbb F}(a)\circ g\circ b\big)$, for $a\in
\mathrm{Hom}_{\mathfrak A}(A,A')$, $b\in \mathrm{Hom}_{\mathfrak B}(B',B)$
and $g\in \mathrm{Hom}_{\mathfrak B}(B,{\mathbb F}(A))$) implies the
naturality of $(\alpha_{{\mathbb F},{\mathbb G}})^{-1}({\mathcal P})$.
Furthermore, (keeping the notation),
$$
\alpha_{{\mathbb F},{\mathbb G}}
\big((\alpha_{{\mathbb F},{\mathbb G}})^{-1} ({\mathcal P})\big)_{B,A}(g)=
(\alpha_{{\mathbb F},{\mathbb G}})^{-1}({\mathcal P})_A\circ g=
{\mathcal P}_{{\mathbb F}(A),A}({\mathbb F}(A))\circ g={\mathcal P}_{B,A}(g),
$$
where the last equality follows by the naturality of ${\mathcal P}$. Also,
$$
(\alpha_{{\mathbb F},{\mathbb G}})^{-1}
\big(\alpha_{{\mathbb F},{\mathbb G}}(\sigma)\big)_A=
\alpha_{{\mathbb F},{\mathbb G}}(\sigma)_{{\mathbb F}(A),A}({\mathbb F}(A))=
\sigma_A,
$$
{what completes the proof.}
\end{proof}
\begin{theorem}\label{thm:rrafael}
Let $(\mathbb{T},\mathbb{H})$ be an adjunction of functors
${\mathbb{T}}:\mathfrak{A}\rightarrow \mathfrak{B}$ and
${\mathbb{H}}:\mathfrak{B}\rightarrow \mathfrak{A}$, with unit
$\eta$ and counit $\varepsilon$. Consider any functor
$\mathbb{R}:\mathfrak{R}\rightarrow\mathfrak{A}$ and a functor
$\mathbb{L}:\mathfrak{L}\rightarrow \mathfrak{A}$ which is
surjective on the objects (e.g. the identity functor $\mathbb{L}=
{\mathfrak{A}}$). Then the following assertions hold.
\begin{itemize}
\item[$1)$] ${\mathbb{T}}$ is $\left( \mathbb{L},\mathbb{R}\right)
$-faithful if and only if it is
$\left({\mathfrak{A}},\mathbb{R}\right)$-faithful and if and only
if $\eta
  _{\mathbb{R}\left( R\right) }$ is a monomorphism, for every object $R\in
  \mathfrak{R}$.
\item[$2)$] ${\mathbb{T}}$ is $\left(
\mathbb{L},\mathbb{R}\right)$-full
  if and only if it is $\left({\mathfrak{A}},\mathbb{R}\right) $-full and if
  and only if $\eta
  _{\mathbb{R}\left( R\right) }$
  is a split epimorphism, for every object $R\in \mathfrak{R}$.
\item[$3)$] (Rafael type Theorem) ${\mathbb{T}}$ is
$\left({\mathfrak{A }},\mathbb{R}\right)$-separable if and only if
$\eta _{\mathbb{R}\left( \bullet \right) }$
is a split natural monomorphism.
\item[$4)$] (Dual Rafael type Theorem) ${\mathbb{T}}$ is
 $\left({\mathfrak{A}},\mathbb{R}\right)$-coseparable if and only if
 $\eta_{\mathbb{R}\left( \bullet \right) }$
is a split natural epimorphism.
\end{itemize}
\end{theorem}

\begin{proof}
Recall that the natural transformation $\phi $ in \eqref{eq:isophi} is an
isomorphism.

$1)$ $\left( \mathbb{L},\mathbb{R}\right) $-faithfulness of
$\mathbb{T}$, i.e. injectivity of the natural transformation
$\Phi(\mathbb{T},\mathbb{L},\mathbb{R})_{L,R}$ in
\eqref{eq:nat.tr.U}, for every object $L\in \mathfrak{L}$ and
$R\in\mathfrak{R}$, is equivalent to injectivity of $\Omega
_{L,R}$, for every $L\in \mathfrak{L}$ and $R\in \mathfrak{R}$.
Since $\mathbb{L}$ is surjective on the objects, in light of
\eqref{eq:Omega} this is equivalent to saying that $\eta
_{\mathbb{R}\left( R\right) }$ is a monomorphism for every $R\in
\mathfrak{R}$.

$2)$ $\left( \mathbb{L},\mathbb{R}\right) $-fullness of
$\mathbb{T}$, i.e. surjectivity of the natural transformation
$\Phi(\mathbb{T},\mathbb{L},\mathbb{R})_{L,R}$ in
\eqref{eq:nat.tr.U}, for every object $L\in \mathfrak{L}$ and
$R\in \mathfrak{R}$, is equivalent to surjectivity  of $\Omega
_{L,R}$, for every $L\in \mathfrak{L}$ and $R\in \mathfrak{R}$.
Let us prove that this is equivalent to saying that $\eta
_{\mathbb{R}\left( R\right) }$ is a split epimorphism for every
$R\in \mathfrak{R.}$ In fact, since $\mathbb{L}$ is surjective on
the objects, for every $R\in \mathfrak{R}$ there exists an object
$L\in \mathfrak{L}$ such that
$\mathbb{H}{\mathbb{T}}\mathbb{R}(R)=\mathbb{L}(L)$. Thus if
$\Omega _{L,R}$ is surjective then, by
 ${\mathbb{H}{\mathbb{T}}\mathbb{R}(R)}\in
\Hom_{\mathfrak{A}}(\mathbb{H}{\mathbb{T}}\mathbb{R}(R),\mathbb{H}{\mathbb{T}}
\mathbb{R}(R))=\Hom_{\mathfrak{A}}(\mathbb{L}(L),\mathbb{H{T}R}(R))$,
there exists $\sigma \in \Hom_{\mathfrak{A}}(
\mathbb{L}(L),\mathbb{R}\left( R\right) )$ such that $\eta _{\mathbb{R}
\left( R\right) }\circ \sigma ={\mathbb{H{T}R}(R)}.$
Conversely, let $g$ be any morphism in
$\Hom_{\mathfrak{A}}(\mathbb{L}(L),\mathbb{H{T}R}(R))$, for some ${L}\in
\mathfrak{L}$ and ${R}\in \mathfrak{R}$.  Let $\sigma$ be a
section of $\eta _{\mathbb{R}\left( R\right)}$.
Define $f\in \Hom_{\mathfrak{A}}(\mathbb{L}
\left( L\right) ,\mathbb{R}\left( R\right) )$ by $f:=\sigma \circ g.$ Then $
\Omega_{L,R} \left( f\right) =\eta _{\mathbb{R}\left( R\right) }\circ
f=g.$

$3)$ $\left({\mathfrak A},\mathbb{R}\right)$-separability of
$\mathbb{T}$, i.e. natural cosplitting of
$\Phi(\mathbb{T},\mathfrak{A},\mathbb{R})$, is equivalent to
natural cosplitting of $\Omega$. Note that $\Omega $ is the image
of the natural transformation $\eta _{{\mathbb{R} }(\bullet )}$
under the functor $\alpha$ in Lemma \ref{lem:ff}. Hence the claim
follows by Lemma \ref{lem:ff}, as a fully faithful functor
preserves and reflects split monomorphisms.

$4)$ $\left( \mathfrak{A},\mathbb{R}\right) $-coseparability of
$\mathbb{T}$, i.e. natural splitting of
$\Phi(\mathbb{T},\mathfrak{A},\mathbb{R})$, is equivalent to
natural splitting of $\Omega$. Hence this claim follows by the
same argument as $3)$ does, as a fully faithful functor preserves
and reflects split epimorphisms as well.
\end{proof}
Dually, one proves the following result.
\begin{theorem}
\label{thm:Dualrrafael} Let $(\mathbb{T},\mathbb{H})$ be an
adjunction of functors ${\mathbb{H}}:\mathfrak{A}\rightarrow \mathfrak{B}$
and ${\mathbb{T}}:\mathfrak{B}\rightarrow \mathfrak{A}$, with unit $\eta$ and
counit $\varepsilon$.
Consider any functor $\mathbb{L}:\mathfrak{L}\rightarrow \mathfrak{A}$ and a
functor $\mathbb{R}:\mathfrak{R}\rightarrow \mathfrak{A}$ which is
surjective on the objects (e.g. the identity functor
$\mathbb{R}={\mathfrak{A}}$).
Then the following  assertions hold.
\begin{itemize}
\item[$1)$] ${\mathbb{H}}$ is $\left( \mathbb{L},\mathbb{R}\right)
$-faithful if and only if it is
$\left(\mathbb{L},{\mathfrak{A}}\right)$-faithful and if and only
if $\varepsilon    _{\mathbb{L}\left( L\right) }$
  is an epimorphism for every object $L\in \mathfrak{L}$.
\item[$2)$] ${\mathbb{H}}$ is $\left( \mathbb{L},\mathbb{R}\right)
$-full if and only if it is $\left(
\mathbb{L},{\mathfrak{A}}\right)$-full and if and only if
$\varepsilon _{\mathbb{L}\left( L\right) }$ is a split
monomorphism for every object $L\in \mathfrak{L}$.
\item[$3)$] (Rafael type Theorem) ${\mathbb{H}}$ is $\left(
  \mathbb{L},\mathfrak{A}\right) $-separable if and only if $\varepsilon
  _{\mathbb{L}\left( \bullet \right) }$
is a split natural epimorphism. \
\item[$4)$] (Dual Rafael type Theorem) ${\mathbb{H}}$ is $\left(
  \mathbb{L},\mathfrak{A}\right) $ -coseparable if and only if $\varepsilon
  _{\mathbb{L}\left(\bullet \right) }$
is a split natural monomorphism.
\end{itemize}
\end{theorem}

A notion somewhat reminiscent to our {\em relative separability}
of a functor was introduced in \cite{CaeMil:SecK} under the name
of {\em separability of the
  second kind}. Our next task is to find a relation between the two notions.

\begin{definition}
Let $\mathbb{R}:\mathfrak{A}\rightarrow \mathfrak{A}^{\prime }$
and $\mathbb{T}:\mathfrak{A}\rightarrow \mathfrak{B}$ be covariant
functors. Following \cite[Definition 2.1]{CaeMil:SecK} and using
the notation introduced in \eqref{eq:nat.tr.U}, ${\mathbb{T}}$ is
called ${\mathbb R}$\emph{-separable of the second kind} if the
natural transformation $\Phi
({\mathbb{R}},{\mathfrak{A}},{\mathfrak{A}})$ factors through
$\Phi ({\mathbb{T}},{\mathfrak{A}},{\mathfrak{A}})$.
\end{definition}

\begin{proposition}
\label{pro:seck}Let $(\mathbb{T},\mathbb{H})$ and $(\mathbb{
T^{\prime }},\mathbb{H^{\prime }})$ be adjunctions with respective units
$\eta$ and $\eta'$. Consider the following diagrams of functors.
\begin{equation*}
\xymatrix@R=15pt@C=50pt{
  \mathfrak{A} \ar[d]_{\mathbb{T}} \ar[r]^{\mathbb{R}} & \mathfrak{A'}
  \ar[d]^{\mathbb{T}'} \\
  \mathfrak{B}  & \mathfrak{B'}   }
\text{ \qquad } \xymatrix@R=15pt@C=50pt{
  \mathfrak{A} \ar[r]^{\mathbb{R}} & \mathfrak{A'}  \\
  \mathfrak{B} \ar[u]_{\mathbb{H}}  & \mathfrak{B'}
  \ar[u]^{\mathbb{H}'}  }
\end{equation*}
The following assertions are equivalent.
\begin{itemize}
\item[(a)] ${\mathbb{T}}$ is ${\mathbb{R}}$-separable of the second
  kind.
\item[(b)] There exists a natural transformation $\nu :{\mathbb{R}}
  {\mathbb{H}}{\mathbb{T}}\rightarrow {\mathbb{R}}$,
satisfying $\nu _{A}\circ {\mathbb{R}}(\eta _{A})={\mathbb{R}}(A)$, for any
  $A\in {\mathfrak{A}}$.
\end{itemize}
Assume that there exists a natural equivalence $\xi
:\mathbb{H^{\prime}T^{\prime}R}\rightarrow \mathbb{RHT}$ such that
\begin{equation}\label{eq:rel_uns}
\xi \circ \eta _{{\mathbb{R}}(\bullet )}^{\prime }=
{\mathbb{R}}(\eta _{\bullet}).
\end{equation}
Then the following assertion is also equivalent to the foregoing
ones.
\begin{itemize}
    \item[(c)]${\mathbb{T}}^{\prime }$ is $({\mathfrak{A}}^{\prime
    },{\mathbb{R}})$-separable.
\end{itemize}
\end{proposition}

\begin{proof}
$\left( a\right) \Leftrightarrow \left( b\right) $ This
equivalence {was proven in} \cite[Theorem 2.7]{CaeMil:SecK}.

$\left( b\right) \Leftrightarrow \left( c\right) $ This equivalence follows by
  Theorem \ref{thm:rrafael} 3), in view of \eqref{eq:rel_uns}.
\end{proof}

\section{Application to entwining structures}\label{sec:entw.str}

As it is recalled in Section \ref{preli:cor.ext}, a coring $\cD$ over an
algebra $L$ is said to be a right extension of a coring $\cC$ over an algebra
$A$ provided that $\cC$ is a $\cC$-$\cD$ bicomodule, via the left regular
coaction.
{\color{blue}
Under the additional assumption that the coring extension is {\em pure}
(cf. Section \ref{preli:cor.ext}), there exists
}
a $k$-linear functor
${\mathbb R}:\cM^\cC\to \cM^\cD$, making the following diagram,
involving four forgetful functors, commutative.
\begin{equation}\label{fig:Brzfunc}
\xymatrix
{\cM^\cC\ar[rrrr]^{\mathbb R}\ar[rd]^{{\mathbb U}^\cC} &&&&\cM^\cD
\ar[ld]_{{\mathbb U}^\cD}\\
&\cM_A\ar[rd]^{{\mathbb F}^A} && \cM_L\ar[ld]_{{\mathbb F}^L}&\\
&&\cM_k&&
}
\end{equation}
The functor $\mathbb{R}$ was explicitly constructed in \cite{Brz:Corext},
cf. Section \ref{preli:cor.ext}.
In this section we study
{\color{blue}
pure}
coring extensions, especially those ones which arise
from entwining structures,
{\color{blue} cf. Section \ref{ex:entw.pure}}.
 We focus on the problem of $(\cM^\cD,\mathbb{R})$
-separability of the functor $\mathbb{U}^\cD$ in Figure \eqref{fig:Brzfunc}.

The following first result is an easy generalization of \cite[Corollary
3.6]{Brz:str} to 
{\color{blue} pure}
coring extensions.
\begin{proposition}\label{prop:rscorext}
Consider an $L$-coring $\cD$ which is
{\color{blue}
a pure}
right extension of an $A$-coring $\cC$,
and the corresponding functors in Figure \eqref{fig:Brzfunc}. The forgetful
functor ${\mathbb U}^\cD$ is $(\cM^\cD,{\mathbb R})$-separable if and only if
the right $\cD$-coaction $\tau_\cC$ in $\cC$ is a
split monomorphism of left $\cC$-comodules and right
$\cD$-comodules.
\end{proposition}
\begin{proof}
The functor ${\mathbb U}^\cD$ possesses a right adjoint, the
functor $\bullet \otimes_L \cD:\cM_L\to \cM^\cD$ (cf.
\cite[18.13]{BrzWis:cor}). The unit of the adjunction is given by
the $\cD$-coaction $\tau$. Hence, by Theorem \ref{thm:rrafael} 3),
${\mathbb
  U}^\cD$ is $(\cM^\cD,{\mathbb R})$-separable if and only if there exists
a natural retraction $\nu$ of $\tau_{{\mathbb R}(\bullet)}$. Therefore if
${\mathbb U}^\cD$ is $(\cM^\cD,{\mathbb R})$-separable then in particular
$\tau_\cC$ possesses a right $\cD$-colinear retraction $\nu_\cC$. We claim
that $\nu_\cC$ is
also left $\cC$-colinear. Indeed, for any right $A$-module $N$ and $n\in N$,
the map $\cC\to N\otimes_A \cC$, $c\mapsto n\otimes_A c$ is right
$\cC$-colinear. Hence by the naturality of $\nu$,
$$
\nu_{N\sstac{A}\cC}(n\stac{A} c\stac{L} d)=n\stac{A} \nu_\cC (c\stac{L} d),
$$
for $n\in N$, $c\in \cC$ and $d\in \cD$. In particular, taking
$N=A$, we conclude on the left $A$-linearity of $\nu_\cC$.
Furthermore, a right $\cC$-coaction {$\varrho^M:m\mapsto
m^{[0]}\sstac A m^{[1]}$} (being coassociative) is $\cC$-colinear,
hence the naturality of $\nu$ implies
$$
\rho^M\big(\nu_M(m\stac{L} d)\big)=\nu_{M\sstac{A} \cC}\big(\rho^M(m)\stac{L}
d\big),
$$
for any $\cC$-comodule $M$, $m\in M$ and $d\in \cD$.
Therefore
$$
\nu_M(m\stac{L} d)^{[0]}\stac{A} \nu_M(m\stac{L} d)^{[1]}=m^{[0]}\stac{A}
\nu_\cC (m^{[1]}\stac{L} d).
$$
Taking $M=\cC$ we have the left $\cC$-colinearity of $\nu_\cC$ proven.

Conversely, let ${\tilde \nu}$ be a left $\cC$-colinear right $\cD$-colinear
retraction of $\tau_\cC$. The natural transformation $\nu$ is constructed as
follows. For any right $\cC$-comodule $M$, put
\begin{equation}\label{eq:nu_tilde}
\nu_M: M\stac{L} \cD\to M, \qquad m\stac{L} d\mapsto m^{[0]}\epsilon_\cC\circ
   {\tilde \nu}(m^{[1]}\stac{L} d).
\end{equation}
Its naturality is obvious. It follows by the $\cD$-colinearity of
a $\cC$-coaction $\rho^M$ that $\nu_M$ is a retraction of
$\tau_M{:  m\mapsto m_{[0]}\sstac L m_{[1]}}$. Indeed,
$$
\nu_M\circ \tau_M(m)={m_{[0]}}^{[0]}\epsilon_\cC\circ
{\tilde \nu}({m_{[0]}}^{[1]}\stac{L} m_{[1]})=m^{[0]}
\epsilon_\cC\circ {\tilde \nu}\circ \tau_\cC(m^{[1]})=m.
$$
It remains to check the $\cD$-colinearity of $\nu_M$. For $m\sstac{L} d\in
M\sstac{L} \cD$,
\begin{eqnarray*}
\tau_M\circ\nu_M(m\stac{L} d)
&=& \big( m^{[0]}\epsilon_\cC\circ {\tilde \nu}(m^{[1]}\stac{L} d)\big)_{[0]}
\stac{L}
\big( m^{[0]}\epsilon_\cC\circ {\tilde \nu}(m^{[1]}\stac{L} d)\big)_{[1]} \\
&=& \big( m^{[0]}\epsilon_\cC\circ {\tilde \nu}(m^{[1]}\stac{L} d)\big)^{[0]}
\epsilon_\cC\big({\big(m^{[0]}\epsilon_\cC\circ {\tilde \nu}(m^{[1]}\stac{L}
d)\big)^{[1]}}_{[0]}\big)  \stac{L}
{\big( m^{[0]}\epsilon_\cC\circ {\tilde \nu}(m^{[1]}\stac{L}
  d)\big)^{[1]}}_{[1]} \\
&=& m^{[0]} \epsilon_\cC\big(
\big( m^{[1]}\epsilon_\cC\circ {\tilde \nu}(m^{[2]}\stac{L} d)\big)_{[0]}
\big)\stac{L}
\big( m^{[1]}\epsilon_\cC\circ {\tilde \nu}(m^{[2]}\stac{L} d)\big)_{[1]}\\
&=& m^{[0]} \epsilon_\cC\big( {\tilde \nu}(m^{[1]}\stac{L} d)_{[0]}\big)
\stac{L} {\tilde \nu}(m^{[1]}\stac{L} d)_{[1]}\\
&=& m^{[0]} \epsilon_\cC\circ {\tilde \nu}(m^{[1]}\stac{L} d_{(1)})
\stac{L} d_{(2)}
= (\nu_M\stac{L} \cD)\circ (M\stac{L} \Delta_\cD)(m\stac{L} d),
\end{eqnarray*}
where the second equality follows by the explicit form
of the functor ${\mathbb R}$, relating $\tau_M$ to $\varrho^M$,
cf. \eqref{eq:D_coac}, the third one follows by the right
$A$-linearity of a $\cC$-coaction, and the fourth and fifth
equalities follow by the left $\cC$-colinearity and the right
$\cD$-colinearity of ${\tilde \nu}$, respectively.
\end{proof}

If the two corings $\cC$ and $\cD$ are equal and ${\mathbb R}$ is the identity
functor, then Proposition \ref{prop:rscorext} reduces to \cite[Corollary
  3.6]{Brz:str}.
More generally, if $\cC$ and $\cD$ are corings over the same base
  algebra $A$
{\color{blue}
and the right $A$-actions of the $A$-coring $C$ and the right $\cD$-comodule
${\mathcal C}$ coincide},
then $\cD$ is a right extension of $\cC$ if and only if there
  exists a homomorphism of $A$-corings $\kappa:\mathcal{C}\to \mathcal{D}$ (in
  terms of  which the $\mathcal{D}$-coaction on $\mathcal{C}$ is given by
  $\tau_{\mathcal{C}}:= (\mathcal{C}\otimes_A \kappa)\circ
  \Delta_{\mathcal{C}}$), cf. \cite[Corrigendum]{BohmBrz:pre_tor}. In this case,
  using the same methods in \cite[Corollary 3.6]{Brz:str},
  the map $\tau_{\mathcal C}$ is checked to be a split monomorphism of left
  $\mathcal{C}$-comodules and right $\mathcal{D}$-comodules (i.e. the functor
  ${\mathbb U}^\cD$ in Figure \eqref{fig:Brzfunc} is checked to be
  $(\cM^\cD,{\mathbb R})$-separable) if and only if
  there exists an $A$-$A$ bimodule map $\widehat{\nu}:\mathcal{C} \otimes_A
  \mathcal{D}\to A$, such that
$$
\widehat{\nu} \circ (\mathcal{C}\stac A \kappa)\circ \Delta_{\mathcal{C}}
=\epsilon_{\mathcal{C}} \qquad \textrm{and} \qquad
\kappa\circ (\mathcal{C}\stac A \widehat{\nu} )\circ (\Delta_{\mathcal{C}}
\stac A \mathcal{D}) = (\widehat{\nu}\stac A\mathcal{D})\circ
(\mathcal{C} \stac A \Delta_{\mathcal{D}}).
$$
This extends \cite[Theorem 3.5]{Brz:str}. On the other hand, for an
arbitrary
{\color{blue}
pure}
coring extension $\cD$ of $\cC$, \cite[Corollary
3.6]{Brz:str} together with Theorem \ref{teo compos of separable}
1) implies that if $\cD$ is a
coseparable coring then the functor ${\mathbb U}^\cD$ in Figure
\eqref{fig:Brzfunc} is $(\cM^\cD,{\mathbb R})$-separable. This fact follows
alternatively by Proposition \ref{prop:rscorext}: if $\zeta$ is a
$\cD$-$\cD$ bicolinear retraction of $\Delta_\cD$, then $(\cC\sstac{L}
\epsilon_\cD\circ \zeta)\circ (\tau_\cC\sstac{L}\cD)$ is a $\cC$-$\cD$
bicolinear retraction of $\tau_\cC$.

Note that, {by Corollary \ref{coro Rafael} 2),} for any
{\color{blue}
pure}
coring extension $\cD$ of $\cC$, $(\cM^\cD,{\mathbb R})$-separability of
${\mathbb U}^\cD$ implies in particular that every right
$\cC$-comodule is relative injective as a right $\cD$-comodule. In
what follows we turn to analyzing more consequences of
$(\cM^\cD,\mathbb{R})$-separability of $\mathbb{U}^\cD$, for
coring extensions arising from entwining structures $(A,\cD,\psi)$
over an algebra $L$. As the main results of the section, {Theorem
\ref{thm:rsentw} and Proposition \ref{prop:cdrinj}} show that if
$\psi$ is bijective and there exists a grouplike element in $\cD$,
then $(\cM^\cD,\mathbb{R})$-separability of $\mathbb{U}^\cD$
implies that $A$ is relative injective {\em also as an entwined
module}. A key notion of our study is the following generalization
of Doi's total integral in \cite{Doi:totint}.
\begin{definition}\label{def:totint}
Let $(A,\cD,\psi)$ be an entwining structure over an algebra $L$.
Assume that $\cD$ possesses a grouplike element $e$ so that $A$ is
a right $\cD$-comodule with coaction $a\mapsto \psi(e\sstac L a)$, cf.
\eqref{eq:A.entw.mod}. A
right $\cD$-comodule map $j:\cD\to A$, satisfying the
normalization condition $j(e)=1_A$, is called a {\em right total
integral}.

In a bijective  entwining structure $(A,\cD,\psi)$ over an algebra
$L$, such that $\cD$ possesses a grouplike element $e$, a {\em
left total integral} is defined as a right total integral in the
$L^{op}$-entwining structure $(A^{op},\cD_{cop},\psi^{-1})$. This
is the same as a left $\cD$-comodule map $j:\cD\to A$, with
respect to the coaction $a\mapsto \psi^{-1}(a\sstac L e)$, cf.
\eqref{eq:entw.mod.A}, satisfying the normalization condition $j(e)=1_A$.
\end{definition}

Consider an entwining structure $(A,\cD,\psi)$ over an algebra
$L$, and denote by $\cC$ the associated $A$-coring $A\sstac L \cD$
(cf. Section \ref{rem:entwcorext}). Consider the following diagrams of
functors \begin{equation}\label{eq:fig_seck}
\xymatrix@R=15pt@C=50pt{
  {\cM^\cC\cong \cM^\cD_A(\psi)} \ar[d]_{\mathbb{T}={\mathbb U}^\cC}
  \ar[r]^{\mathbb{R}} & {\cM^\cD }
  \ar[d]^{{\mathbb{T}'={\mathbb U}^\cD}} \\
  \cM_A
& {\cM_L}   }
\text{ \qquad } \xymatrix@R=15pt@C=50pt{
  {\cM^\cC\cong \cM^\cD_A(\psi)} \ar[r]^{\mathbb{R}} & \cM^\cD  \\
  \cM_A
\ar[u]_{\mathbb{H}=\bullet\sstac{L} \cD}  & \cM_L
  \ar[u]^{\mathbb{H}'=\bullet\sstac{L} \cD}  }
\end{equation}
where $\mathbb{T}={\mathbb U}^\cC$,$\mathbb{T}'={\mathbb U}^\cD$
and $\mathbb{R}$ are forgetful functors (cf. Figure
\eqref{fig:Brzfunc}). Note that $(\mathbb{T},\mathbb{H})$ and
$(\mathbb{T^{\prime}},\mathbb{H^{\prime }})$ are adjunctions and
the respective units $\eta$ and $\eta'$ are given by the right
$\cD$-coaction, in both cases (cf. {Section} \ref{preli:coring}).
Hence they satisfy ${\mathbb R}(\eta_\bullet)= \eta'_{{\mathbb
R}(\bullet)}$.

\begin{theorem}\label{thm:rsentw}
Let $(A,\cD,\psi)$ be an entwining structure over an algebra $L$.
Consider the functors in Figure \eqref{eq:fig_seck}. The following assertions
are equivalent.

\begin{itemize}
\item[$(a)$] ${\mathbb{T}}={\mathbb U}^\cC$ is
${\mathbb{R}}$-separable of the second kind.

\item[$(b)$]  There exists a natural transformation $\nu :{\mathbb{R}
}{\mathbb{H}}{\mathbb{T}}\rightarrow {\mathbb{R}}$,
satisfying $\nu _{M}\circ {\mathbb{R}}(\eta _{M})={\mathbb{R}}
(M)$, for any $M\in {\cM^\cD_A(\psi)}$.

\item[$(c)$]  ${\mathbb{T}}^{\prime }={\mathbb U}^\cD$ is $(\cM^\cD,{
\mathbb{R}})$-separable.

\item[$(d)$]  There exists a morphism $\theta \in
  {}_L\mathrm{Hom}_L(\cD\sstac{L}\cD,A)$ satisfying, for all $d,d'\in \cD$,
  \begin{equation}\label{eq:CM.tot.int}
\theta(d\stac{L} d'_{(1)})\stac{L}
d'_{(2)}=\psi\big(d_{(1)}\stac{L}\theta(d_{(2)} \stac{L}
d')\big)\quad \textrm{and}\quad \theta(d_{(1)}\stac{L}
d_{(2)})=\eta_A\circ
\epsilon_\cD(d).
\end{equation}\end{itemize}
If these equivalent conditions hold, and in addition there exists
a grouplike element in $\cD$, then there exists a right total
integral in the $L$-entwining structure $(A,\cD,\psi)$.
\end{theorem}

\begin{proof}
The equivalence of assertions $(a),(b)$ and $(c)$ is a consequence of
  Proposition  \ref{pro:seck}.

The equivalence of assertions $(a)$ and $(d)$
is proven by an easy extension to non-commutative base of
arguments in \cite[Proposition 4.12]{CaeMil:SecK}, about entwining structures
over commutative rings.

Assume that there exists a grouplike element $e$ in $\cD$, hence
$A$ is a right $\cD$-comodule with coaction \eqref{eq:A.entw.mod}.
In this situation the map
\begin{equation*}
j:\cD\to A\qquad d\mapsto \theta(e\stac L d)
\end{equation*}
  is right $\cD$-colinear and satisfies the normalization condition $j(e)=
  1_A$. That is, $j$ is a right total integral in the sense of Definition
  \ref{def:totint}.
\end{proof}

Note that, following the proof of \cite[Proposition
4.12]{CaeMil:SecK}, a bijective correspondence can be obtained
between maps $\theta$ as in \eqref{eq:CM.tot.int} and left
$\cC=(A\sstac{L} \cD)$-colinear right $\cD$-colinear retractions
of the $\cD$-coaction $A\sstac{L} \Delta_\cD$. The explicit
relation is given by the same formulae as in \cite{CaeMil:SecK},
in the paragraph preceding Proposition 4.12. Since in view of
Proposition \ref{prop:rscorext} the existence of a left
$\cC=(A\sstac{L} \cD)$-colinear right $\cD$-colinear retraction of
the $\cD$-coaction $A\sstac{L} \Delta_\cD$ is equivalent to
assertion (c) in Theorem \ref{thm:rsentw}, in \cite[Proposition
4.12]{CaeMil:SecK} implicitly also the equivalence of assertions
(a) and (c) in Theorem \ref{thm:rsentw} is proven.

In contrast to \cite{CaeMil:SecK}, in the current paper
  the term {\em total integral} is used only in the more restricted sense of
  Definition   \ref{def:totint}.

The following proposition extends \cite[Proposition
4.2]{BohmBrz:relchern}. It clarifies the role of total integrals
in bijective entwining structures with a grouplike element.
For the notion of coinvariants, with respect to a grouplike element in
  a coring, consult Section \ref{preli:coring}.

\begin{proposition}\label{prop:cdrinj}
Consider a {\em bijective} entwining structure $(A,\mathcal{D},\psi )$ over an
algebra $L$, such that there exists a
grouplike element $e$ in $\mathcal{D}$. Let $\mathcal{C}:
=A\otimes _{L}\mathcal{D}$ be the associated $A$-coring.
The following assertions are equivalent.
\begin{itemize}
\item[$\left( a\right) $] $A$ is a relative injective right (resp. left)
  $\mathcal{C}$-comodule.
\item[$\left( b\right) $] $A$ is a relative injective right (resp. left)
  $\mathcal{D}$-comodule.
\item[$(c)$] There exists a {right }(resp. {left}) total
  integral in the entwining structure $(A,\mathcal{D},\psi )$.
\end{itemize}
If these equivalent conditions hold then $B:=A^{co\mathcal{C}}=
A^{co\mathcal{D}}$ is a direct summand of $A$ as a right (resp. left)
$B$-module.
\end{proposition}

\begin{proof}
$(a)\Rightarrow (b)$ For a relative injective right $\mathcal{C}$-comodule $
M $, the right $\mathcal{D}$-coaction has a right $A$-linear right
$\mathcal{D}$-colinear retraction. Hence it is a relative injective right
$\mathcal{D}$ -comodule.

$(b)\Rightarrow (a)$ Assume that $A$ is a relative injective right $\mathcal{
D}$-comodule. Similarly to the proof of \cite[Lemma 4.1]{SchSch:genGal}, in
terms of a right $\mathcal{D}$-colinear retraction $\nu _{A}$
of the $\mathcal{D}$-coaction
\eqref{eq:A.entw.mod} in $A$, a right $\mathcal{C}$-colinear retraction is
given by
\begin{equation*}
\mu _{A}\circ \lbrack \nu _{A}(1_{A}\stac {L}\bullet )\stac  {L}A]\circ \psi
  ^{-1}:A\stac {L}\mathcal{D}\rightarrow A.
\end{equation*}

The equivalence $(b)\Leftrightarrow (c)$ was proven in
  \cite[Proposition 4.2]{BohmBrz:relchern}, as follows. To a right
  $\cD$-colinear
  retraction $\nu_A$ of the $\cD$-coaction \eqref{eq:A.entw.mod} in $A$, one
  associates a right total
  integral $j:d\mapsto \nu_A(1_A\sstac L d)$. Conversely, in terms of a right
  total integral $j$, a right $\cD$-colinear retraction of the $\cD$-coaction
  \eqref{eq:A.entw.mod} in $A$ is constructed as $\nu_A:= \mu_A\circ(j\sstac L
  A)\circ \psi^{-1}$.

It remains to prove the last statement. By property $\left(
a\right)$, the right $\mathcal{C}$-coaction in $A$ is a split
monomorphism in $\mathfrak{M}^{\mathcal{C}}$. Taking the
$\mathcal{C}$-coinvariants part (with respect to the grouplike
element $1_{A}\otimes _{L} e$) of its retraction, we obtain right
$B$-linear retraction of the inclusion $B\rightarrow A$.

In order to prove the claim about the left comodule structures,
the same
arguments can be applied to the entwining structure $(A^{op},\mathcal{C}
_{cop},\psi ^{-1})$ over the algebra $L^{op}$.
\end{proof}

The following Lemma is a simple generalization of \cite[Lemma
4.1]{SchSch:genGal}.

\begin{lemma}\label{lem:coinv}
Let $\cC$ be an $A$-coring possessing a grouplike element $g$.
Assume that $A$ is a relative injective left $\cC$-comodule via
the coaction $a\mapsto ag$, determined by $g$. Denote by $B\colon\!\!\!
=A^{co\cC}$ the coinvariants of $A$ with respect to $g$.
Then the unit of the adjunction $(\bullet\sstac{B}
A,(\bullet)^{co\cC})$, i.e. the natural transformation \eqref{eq:co_un},
is an isomorphism.
\end{lemma}

\begin{proof}
Let $M$ be a relative injective left $\cC$-comodule
and $\nu_M$ be a left $\cC$-colinear retraction of the coaction ${}^M\rho$.
Introduce a further map $\xi^M:M\to \cC\sstac{A} M$, $m\mapsto
g\sstac{A} m$. We claim that
$$
\xymatrix { M\ar@<1.3ex>[rrr]^{\hspace{-.3cm}
{}^M\rho}\ar@<-1.3ex>[rrr]_{\hspace{-.3cm}
  \xi^M}
&&& \cC\stac{A} M\ar[lll]|{\hspace{.1cm} {\nu_M}} }
$$
is a contractible pair in ${}_B\cM$. Clearly, all morphisms ${}^M\rho$,
$\xi^M$ and $\nu_M$ are left $B$-linear. By definition $\nu_M\circ
{}^M\rho=M$. Hence we conclude after observing that, for $m\in M$,
$$
{}^M\rho\circ \nu_M\circ \xi^M(m)=\nu_M(g\stac{A}
m)^{[-1]}\stac{A} \nu_M(g\stac{A}m)^{[0]}=g\stac{A}
\nu_M(g\stac{A} m)= \xi^M\circ \nu_M\circ \xi^M(m),
$$
where in the second equality the left $\cC$-colinearity of $\nu_M$
has been used. {In particular we deduce that} the equalizer
of ${}^A\rho:A\to \cC$, $a\mapsto ag$ and $\xi^A:A\to \cC$,
$a\mapsto ga$ cosplits in ${}_B\cM$. Hence it is preserved by the
functor $N \sstac{B} \bullet: {}_B\cM\to \cM_k$, for any right
$B$-module $N$. Recall that $A$ is a right $\cC$-comodule with
coaction $\xi^A$ and $N\sstac{B}A$ is a right $\cC$-comodule with
coaction $N\sstac{B} \xi^A$. Therefore
$$
(N\stac{B} A)^{co\cC}=\mathrm{Ker}(N\stac{B} {}^A\rho-N\stac{B}
\xi^A)= N\stac{B} \mathrm{Ker}({}^A\rho-\xi^A) = N \stac{B}
B\cong N.
$$
{This proves that \eqref{eq:co_un} is a natural isomorphism, as
  stated.}
\end{proof}

The following proposition formulates a functorial criterion
for a coring with a grouplike element to be a Galois coring.

\begin{proposition} \label{pro:fullyfaith} Let
$\mathcal{C}$ be an $A$-coring possessing a
grouplike element $g$. Denote by $B:=A^{co\mathcal{C}}$ the coinvariants of $
A$ with respect to $g$. Consider the adjunction $(\bullet \otimes
_{B}A,(\bullet)^{co\mathcal{C}})$ in {Section} \ref{preli:coring}
and the canonical map
$\can:A\sstac B A \to \cC$
in \eqref{eq:cor_can}. The following statements hold.
\begin{enumerate}
\item[1)] $\mathrm{can}$ is an epimorphism if and only if the
functor
  $(\bullet )^{co\mathcal{C}}$ is $\left( \bullet \otimes
    _{A}\mathcal{C},\mathfrak{M}^{\mathcal{C}}\right) $-faithful.

\item[2)] $\mathrm{can}$ is a split monomorphism if and only if
the functor $(\bullet)^{co\mathcal{C}}$ is $\left( \bullet \otimes
  _{A}\mathcal{C},\mathfrak{M}^{\mathcal{C}}\right) $-full.
\end{enumerate}
In particular, $\mathcal{C}$ is a Galois coring if and only
if the functor $(\bullet)^{co\mathcal{C}}$ is $\left( \bullet
\otimes
  _{A}\mathcal{C},\mathfrak{M}^{\mathcal{C}}\right) $-fully faithful.
\end{proposition}
\begin{proof}
Denote the counit of the coring $\cC$ by $\epsilon$. For any right
$A$-module $M$, the counit $n$ of the adjunction $(\bullet \otimes
_{B}A,(\bullet)^{co\mathcal{C}})$ (cf. \eqref{eq:co_coun}) is
subject to the equality of maps $(M\sstac A \cC)^{co\cC}\sstac B
A\to M\sstac A \cC$,
\begin{equation*}
( M\stac {A}\mathrm{can}) \circ (M\stac A \epsilon\stac {B}A) =
n_{M\sstac {A}\mathcal{C}}.
\end{equation*}
Since the restriction of $M\sstac A \epsilon$ is an isomorphism
$(M\sstac A \cC)^{co\cC}\to M\sstac A A$, the claims follow by
Theorem \ref{thm:Dualrrafael} 1) and 2), respectively.
\end{proof}

\section{Comodule algebras of Hopf algebroids}\label{sec:hgd}

Consider a Hopf algebroid $\hH$, with constituent left bialgebroid
  $\hH_L$ over the base algebra $L$ and right bialgebroid
  $\hH_R$ over $R$ (cf. Section \ref{preli:hgd}),
and a right $\hH$-comodule algebra $A$ (cf. Appendix \ref{def:Hopf_mod}). 
Recall (from {Section} \ref{preli:hgd.comodule}) that latter means a right
$\hH_R$-comodule algebra
{\color{blue}
and}
right $\hH_L$-comodule algebra $A$, with coactions $a\mapsto
  a^{[0]}\sstac R a^{[1]}$ and $a\mapsto a_{[0]}\sstac L a_{[1]}$,
  respectively, related as in
{\color{blue}
\eqref{eq:hgd_coac}}.
Recall (from Appendix \ref{preli:left-right}) that if the antipode of $\hH$ is
bijective, then the right $\hH$-comodule algebra structure of $A$ is
equivalent to a left $\hH$-comodule algebra structure of $A^{op}$.
{\color{blue}
Consider the forgetful functors
\begin{equation}\label{eq:U}
\xymatrix{
\cM^\hH_A \ar[rr]^-{{\mathbb R}}&&
\cM^\hH \ar[rr]^-{{\mathbb V}}&&
\cM^{\hH_L} \ar[rr]^-{{\mathbb U}}&&
\cM_L.
}
\end{equation}
In this section we study
relative separability of ${\mathbb U}$ with respect to ${\mathbb V}
{\mathbb R}$ and relative separability of ${\mathbb U} {\mathbb V} $ with
respect to ${\mathbb R}$. }

{\color{blue}
\begin{theorem}\label{thm:4.1.w}
Consider a Hopf algebroid $\hH$, with constituent
left bialgebroid $\hH_L$, right bialgebroid $\hH_R$ and antipode $S$.
For a right $\hH$-comodule algebra $A$, the following assertions are equivalent.
\begin{itemize}
\item[{(a)}] There exists a right total integral in the (bijective)
  $L$-entwining structure \eqref{eq:hgdentwL} with grouplike element
$1_H$, i.e. a morphism $j\in {\rm Hom}^{\mathcal H_L}(H,A)$, normalized as
$j(1_H)=1_A$.
\item[{(b)}] $A\in \cM^{\hH_L}$
is ${\mathcal I}_{{\mathbb U}}$-injective
  (i.e. $A$ is a relative injective right $\hH_L$-comodule).
\item[{(c)}] Any object in the image of ${\mathbb V} {\mathbb R}$ is
  ${\mathcal I}_{{\mathbb U}}$-injective (i.e. injective {\em with respect to}
  ${\mathbb U}$). 
\item[{(d)}] The functor ${\mathbb U}$ is $(\cM^{\hH_L},{\mathbb V}
  {\mathbb R})$-separable.
\end{itemize}
If the antipode of $\hH$ is bijective then the following statements are also
equivalent to the foregoing ones.
\begin{itemize}
\item[{(e)}] There exists a {left total integral} in the
bijective $R$-entwining structure \eqref{eq:hgdentwR} with grouplike element
$1_H$, i.e. a left $\hH_R$-colinear map 
$j^{op}_{cop}:H\to A$, normalized as $j^{op}_{cop}(1_H)=1_A$.
\item[{(f)}] $A$ is a relative injective left $\hH_R$-comodule.
\item[{(g)}] Any object of ${}^\hH\cM_A$ is a relative injective left
  $\hH_R$-comodule.
\item[{(h)}] The forgetful functor ${}^{\hH_R}\cM\to {}_R\cM$ is
  $({}^{\hH_R}\cM, ({\mathbb V}{\mathbb R})^{op}_{cop})$-relative separable, where
  $({\mathbb V}{\mathbb R})^{op}_{cop}$ denotes the forgetful functor
  ${}^\hH\cM_A \to {}^{\hH_R}\cM$.
\end{itemize}
If the antipode of $\hH$ is bijective then the following statements are also
equivalent to each other (but not necessarily to the foregoing ones).
\begin{itemize}
\item[{(i)}] There exists a {right total integral} in the $R$-entwining
structure \eqref{eq:hgdentwR} with grouplike element $1_H$, i.e. a morphism
$j^{op}\in {\rm Hom}^{\mathcal H_R}(H,A)$, normalized as $j^{op}(1_H)=1_A$.
\item[{(j)}] $A$ is a relative injective right $\hH_R$-comodule.
\item[{(k)}] Any object of ${}_A \cM^\hH$ is a relative injective right
  $\hH_R$-comodule. 
\item[{(l)}] The forgetful functor $\cM^{\hH_R}\to \cM_R$ is
  $(\cM^{\hH_R}, ({\mathbb V}{\mathbb R})^{op})$-relative separable, where
  $({\mathbb V}{\mathbb R})^{op}$ denotes the forgetful functor ${}_A\cM^\hH \to
  \cM^{\hH_R}$.
\item[{(m)}] There exists a left total integral in the bijective
$L$-entwining structure \eqref{eq:hgdentwL} with grouplike element $1_H$,
i.e. a left $\hH_L$-colinear map 
$j_{cop}:H\to A$, normalized as $j_{cop}(1_H)=1_A$.
\item[{(n)}] $A$ is a relative injective left $\hH_L$-comodule.
\item[{(o)}] Any object of ${}^\hH_A\cM$ is a relative injective left
  $\hH_L$-comodule. 
\item[{(p)}] The forgetful functor ${}^{\hH_L}\cM\to {}_L\cM$ is
  $({}^{\hH_L}\cM, ({\mathbb V}{\mathbb R})_{cop})$-relative separable, where
  $({\mathbb V}{\mathbb R})_{cop}$ denotes the forgetful functor ${}^\hH_A\cM\to
  {}^{\hH_L}\cM$.
\end{itemize}
\end{theorem}}

\begin{proof}
$(a)\Leftrightarrow (b)$ This equivalence follows by Proposition
\ref{prop:cdrinj} $(b)\Leftrightarrow (c)$, since the entwining
map \eqref{eq:hgdentwL} is bijective.

$(a)\Rightarrow (d)$
{\color{blue}
In light of Theorem \ref{thm:rrafael}~3), we need to construct a
  right $\hH_L$-colinear
natural retraction $\nu_M$ of the $\hH_L$-coaction, for any $M\in
\cM^\hH_A$. In terms of the map $j$ in part (a), it is given by the well
defined maps
\begin{equation}\label{eq:nu_M}
\nu_M:M\ot_L H \to M, \qquad m\ot_L h \mapsto m^{[0]}j\big(S(m^{[1]}) h \big),
\end{equation}
where the Sweedler type index notation $m\mapsto m^{[0]}\ot_R m^{[1]}$ is used
for the $\hH_R$-coaction on $M$.}

$(d)\Rightarrow (c)$ The forgetful functor
$\cM^{\hH_L}\to \cM_L$
has a right adjoint, the functor $\bullet\sstac{L} H$,
cf. Section \ref{preli:coring}.
Hence the claim follows by Corollary \ref{coro Rafael} 2).

$(c)\Rightarrow (b)$ This implication is trivial as $A$ itself is an object in
$\cM^{\mathcal H}_A$.

{\color{blue}
If the antipode is bijective then implications $(e)\Leftrightarrow (f)
\Leftrightarrow (g) \Leftrightarrow (h)$ follow by
applying $(a)\Leftrightarrow (b)\Leftrightarrow (c) \Leftrightarrow
(d)$ to the opposite-coopposite Hopf algebroid ${\mathcal
H}^{op}_{cop}$ and its right comodule algebra $A$ (with $({\mathcal
H}_R)^{op}_{cop}$-coaction $a\mapsto  a_{[0]} \sstac{R^{op}} S(a_{[1]})$
and $({\mathcal H}_L)^{op}_{cop}$-coaction $a\mapsto a^{[0]} \sstac{L^{op}}
S(a^{[1]}) $).
}

{\color{blue}
$(a)\Leftrightarrow (e)$ In terms of $j$ in part (a), a normalized left
  $\hH_R$-comodule map is given by $j\circ S^{-1}$. Clearly, it yields a
  bijective correspondence.

The other sequence of equivalences (i)-(p) follows by applying the proven
result to the opposite (or coopposite) Hopf algebroid and its right comodule
algebra $A^{op}$.
}
\end{proof}

{\color{blue}
\begin{theorem}\label{thm:4.1.s}
Let ${\mathcal H}$ be a Hopf algebroid with constituent left
bialgebroid $(H,L,s_L,t_L,\gamma_L,\pi_L)$, right bialgebroid
$(H,R,s_R,t_R,\gamma_R, \pi_R)$, and antipode $S$.
For a right $\hH$-comodule algebra $A$, the following assertions are
equivalent.
\begin{itemize}
\item[{(a)}] There exists a normalized right $\hH$-comodule map $j:H \to A$.
\item[{(b)}] $A \in \cM^\hH$ is ${\mathcal I}_{{\mathbb U} {\mathbb
    V}}$-injective. 
\item[{(c)}] Any object in the image of ${\mathbb R}$ is ${\mathcal
  I}_{{\mathbb U} {\mathbb V}}$-injective.
\item[{(d)}] The functor ${\mathbb U} {\mathbb V}$ is $(\cM^{\hH},{\mathbb
  R})$-separable.
\end{itemize}
If the antipode of $\hH$ is bijective, then these equivalent statements are
equivalent also to the existence of a normalized left $\hH$-comodule map
$H \to A$, hence the symmetrical counterparts of (b)-(d).
\end{theorem}

\begin{proof}
(d)$\Rightarrow$(c) follows by statement 1) in Appendix \ref{prop:ff_adj} and
  Corollary 2.9~2). 

(c)$\Rightarrow$(b) is obvious.

(b)$\Rightarrow$(a) Denote by $\eta:R\to A$ the unit of the $R$-ring $A$.
Since $\eta \circ \pi_R \circ t_L:L \to A$ and $t_L:L \to H$ are
$\hH$-comodule maps and $t_L$ is a split monomorphism of right $L$-modules,
using ${\mathcal I}_{{\mathbb U} {\mathbb V}}$-injectivity of $A$, $j$ is
constructed as the unique $\hH$-comodule map for which $j\circ t_L = \eta
\circ \pi_R \circ t_L$.

(a)$\Rightarrow$(d) We need to construct a natural retraction $\nu_M$ of the
$\hH_L$-coaction, for any object $M$ in $\cM^\hH_A$. In terms of the
map $j$ in part (a), it is given by the same formula \eqref{eq:nu_M}. Since
$j$ is an $\hH$-comodule map, so in $\nu_M$.

If the antipode is bijective then any (normalized) right $\hH$-comodule map
$j:H\to A$ determines a (normalized) left $\hH$-comodule map $j\circ S^{-1}:H
\to A$. This correspondence is clearly bijective.
\end{proof}

Obviously, if the equivalent statements in Theorem \ref{thm:4.1.s} hold then
also the equivalent statements in Theorem \ref{thm:4.1.w} hold. }
\bigskip

{\color{blue}
\begin{lemma}\label{lem:proj_comod}
Let ${\mathcal H}$ be a Hopf algebroid with constituent left bialgebroid
$\hH_L=(H,L,s_L,t_L,\gamma_L,\pi_L)$, right bialgebroid
$\hH_R=(H,R,s_R,t_R,\gamma_R, \pi_R)$, and a bijective antipode $S$. 
Assume that $H$ is a projective right comodule for the $R$-coring
$(H,\gamma_R,\pi_R)$ via $\gamma_R$.
Then $H$ is a projective left $L$-module via left multiplication by $s_L$.
\end{lemma}

\begin{proof}
By \cite[18.20(1)]{BrzWis:cor}, projectivity of $H$ as a right
$\hH_R$-comodule implies that $H$ is a projective right $R$-module via the
action 
\begin{equation}\label{eq:s_R_ac}
H \otimes R \to H, \qquad h\otimes r \mapsto h s_R(r).
\end{equation}
By bijectivity of the antipode, the right $R$-module \eqref{eq:s_R_ac} is
isomorphic to the right $R$-module $H$, with action
\begin{equation}\label{eq:t_R_ac}
H \otimes R \to H, \qquad h\otimes r \mapsto t_R(r) h.
\end{equation}
Hence also the right $R$-module \eqref{eq:t_R_ac} is projective. Furthermore,
the algebra isomorphism $\pi_R \circ s_L:L^{op}\to R$ induces a category
isomorphism $\cM_R \cong {}_L\cM$. This isomorphism takes the projective right
$R$-module \eqref{eq:t_R_ac} to the projective left $L$-module $H$, with
action
$$
L\otimes H \to H,\qquad l\otimes h \mapsto t_R \circ \pi_R \circ
s_L(l)h=s_L(l)h. 
$$
\end{proof}
}

Theorem \ref{thm:4.1.s} makes us able to answer a question which
was left open in \cite{Bohm:hgdGal}. 
{\color{blue}
Consider a Hopf algebroid $\hH$ with a
bijective antipode and a right $\hH$-comodule algebra $A$. Denote
$B:=A^{co\hH_R}=A^{co\hH_L}$, cf. \ref{prop:coinv}.  
By Appendix \ref{rem:bijantip}, one can associate to $A$ four (anti-)
isomorphic corings. Clearly, if any of them is a Galois coring (with respect
to the grouplike element determined by the unit elements in $A$ and $\hH$),
then all of them are Galois corings. In other words, the four properties that
$B \subseteq A$ is a left or right Galois extension by $\hH_R$ or $\hH_L$ are
all equivalent to each other. 
In Proposition \ref{cor:fgp} below, $\hH$-comodule algebras $A$ are
studied, such that these equivalent Galois conditions hold.}

\begin{proposition}\label{cor:fgp} Let ${\mathcal H}$ be a Hopf algebroid with
  constituent left bialgebroid 
$\hH_L=(H,L,s_L,t_L,\gamma_L,$ $\pi_L)$, right bialgebroid 
$\hH_R=(H,R,s_R,t_R,\gamma_R, \pi_R)$, and a bijective antipode $S$.
Assume that $H$ is a projective left $R$-module via $t_R$ and a
projective right comodule for the $R$-coring $(H,\gamma_R,\pi_R)$
via $\gamma_R$. (These assumptions hold e.g. if $H$ is a finitely
generated and projective 
both as a right and left $L$-module and also as a right and left 
$R$-module, 
cf. \cite[Section 4]{Bohm:hgdGal}.) 
{\color{blue}
Then $\cM^{\hH_L}\cong \cM^\hH\cong \cM^{\hH_R}$ and 
${}^{\hH_L}\cM\cong {}^\hH\cM\cong {}^{\hH_R}\cM$
as monoidal categories. Moreover, for a right $\hH$-comodule algebra $A$, such
that $B:=A^{co\hH_R}\subseteq A$ is a right $\hH_R$-Galois extension,}
the following assertions are equivalent.
\begin{itemize}
\item[$(a)$] $A$ is a faithfully flat right $B$-module.
\item[$(b)$] $B$ is a direct summand of the right $B$-module $A$.
\item[$(c)$] The functors $A\sstac{B}\bullet:{}_B\cM\to
{\color{blue} {}^{\hH}_A\cM}$
and ${}^{co{\mathcal H}}(\bullet):
{\color{blue} {}^{\hH}_A\cM}\to {}_B \cM$
  are inverse equivalences and $H\sstac{R} A$ is a flat right $A$-module.
\item[$(d)$] $A$ is a projective generator in
{\color{blue} ${}^{\hH}_A\cM$}
and $H\sstac{R} A$ is a flat right $A$-module.
\item[$(e)$] $A$ is a generator of right $B$-modules.
\item[$(f)$] $A$ is a faithfully flat left $B$-module.
\item[$(g)$] $B$ is a direct summand of the left $B$-module $A$.
\item[$(h)$] The functors $\bullet\sstac{B} A:\cM_B\to
{\color{blue} \cM^{\hH}_A}$
and $(\bullet)^{co{\mathcal H}}:
{\color{blue} \cM^{\hH}_A} \to \cM_B$ are inverse equivalences.
\item[$(i)$] $A$ is a projective generator in
{\color{blue} $\cM^{\hH}_A$}.
\item[$(j)$] $A$ is a generator of left $B$-modules.
\item[$(k)$] {\color{blue}
 The equivalent conditions in Theorem \ref{thm:4.1.s} hold.}
\end{itemize}
\end{proposition}

\begin{proof}
{\color{blue}
Since $H$ is a projective left $R$-module by assumption, it is in particular
flat. 
As a left $L$-module, $H$ is projective hence flat by Lemma
\ref{lem:proj_comod}. Thus the monoidal isomorphisms $\cM^{\hH}\cong
\cM^\hH\cong \cM^{\hH_R}$ follow by \ref{thm:pure}~3) and
\ref{thm:hgd.com.mon}. By Appendix \ref{preli:left-right}, bijectivity of the
antipode implies strict anti-monoidal isomorphisms 
${}^{\hH_R}\cM\cong \cM^{\hH_L}$, 
${}^{\hH_L}\cM\cong \cM^{\hH_R}$ and 
${}^{\hH}\cM\cong \cM^{\hH}$. 
Hence also ${}^{\hH_L}\cM\cong {}^\hH\cM\cong {}^{\hH_R}\cM$.
Note that this implies in particular ${}^{\hH_R}_A\cM\cong {}^\hH_A\cM$ and
$\cM^{\hH_R}_A \cong \cM^\hH_A$.}

$(a)\Leftrightarrow (b)$ and $(f)\Leftrightarrow (g)$ These equivalences
follow by \cite[2.11.29]{Row:rin}, as $A$ is a projective left and right
$B$-module by \cite[Proposition 4.2]{Bohm:hgdGal}.

$(b)\Rightarrow (k)$ and $(k)\Rightarrow (b)$
{\color{blue}
Note that in the current case assertion (k) is equivalent to relative
injectivity of $A$ as a right $\hH_R$-comodule, or as a right
$\hH_L$-comodule, or as a left $\hH_L$-comodule or as a left
$\hH_R$-comodule.}
Since a comodule algebra for a Hopf algebroid with a bijective antipode
determines bijective entwining structures \eqref{eq:hgdentwR} and
\eqref{eq:hgdentwL}, these implications follow by
\cite[Proposition 4.1]{BohmBrz:relchern} (which is a simple
generalization of \cite[Remark 4.2]{SchSch:genGal} to the case of
non-commutative base algebras).

$(k)\Rightarrow (b)$ and $(k)\Rightarrow (g)$ These assertions follow by
Proposition \ref{prop:cdrinj}.

$(a)\Leftrightarrow (d)\Leftrightarrow (c)$ 
{\color{blue}
Since ${}^\hH_A\cM\cong {}^{\hH_R}_A \cM$ is isomorphic to the category of left
comodules for the $A$-coring $H \sstac R A$ (cf. Appendix \ref{rem:entwbgd}),}
these equivalences follow by the
{\em Galois Coring Structure Theorem} \cite[28.19 (2)]{BrzWis:cor}.

$(f)\Leftrightarrow (i)\Leftrightarrow (h)$ Since $H$ is a projective left
$R$-module by assumption, $A\sstac{R} H$ is a projective (hence flat) left
$A$-module. Therefore also these equivalences follow by the
{\em Galois Coring Structure Theorem} \cite[28.19 (2)]{BrzWis:cor}, 
{\color{blue}
as $\cM^\hH_A\cong \cM^{\hH_R}_A\cong \cM^{A\sstac R H}$}.

$(b)\Rightarrow (e)$ and $(g)\Rightarrow (j)$ These implications are trivial.

$(e)\Rightarrow (b)$ and $(j)\Rightarrow (g)$ $A$ is a generator
of right (respectively, left) $B$-modules if and only if there
exist finite sets $\{a_i\}$ in $A$ and $\{\alpha_i\}$ in
$\mathrm{Hom}_B(A,B)$ (respectively, in ${}_B \mathrm{Hom}(A,B)$),
satisfying $\sum_i\alpha_i(a_i)=1_B$. In terms of these elements,
a right $B$-linear retraction of the inclusion $B\to A$ is given
by the map $a\mapsto \sum_i\alpha_i(a_ia)$ (respectively, a left
$B$-linear retraction of the inclusion $B\to A$ is given by the
map $a\mapsto \sum_i\alpha_i(aa_i)$).
\end{proof}

Applying Proposition \ref{cor:fgp} to the co-opposite Hopf
algebroid ${\mathcal
  H}_{cop}$, we see that the claims in  Proposition \ref{cor:fgp}
-- with the only modification that claims (h) and (i) need to be supplemented
  by the assertion that $A\sstac{R} H$ is a flat left $A$-module --
can be proven alternatively by replacing the assumptions about
the projectivity of $H$ a left $R$-module (via $t_R$) and a right
{\color{blue}
${\mathcal H_R}$}-comodule (via the coproduct
{\color{blue}
of $\hH_R$})
with the assumptions that it is a projective right $R$-module (via $s_R$) and
a projective left
{\color{blue}
${\mathcal H}_R$}-comodule (via the coproduct
{\color{blue}
of $\hH_R$}).

\section{A Schneider type theorem}\label{sec:Schneider}

This section contains the main result of the paper, Theorem
  \ref{thm:schhgd}. The starting point of our study is the following result
  \cite[Theorem 2.1]{BrzTurWri:WeakcoaGal}. Recall that
a right comodule $P$ for an $A$-coring $\cC$, which is a finitely generated
and projective right $A$-module, is a {\em Galois comodule} if the canonical
map
\begin{equation}\label{eq:can.coring}
\mathrm{can}:\mathrm{Hom}_A(P,A)\stac{S} P\to \cC,\qquad \phi\stac{S} p\mapsto
\phi(p^{[0]}) p^{[1]}
\end{equation}
is bijective, where $S: = \mathrm{End}^\cC(P)$. {Assume that $S$
is a $T$-ring (e.g. $T$ is a subalgebra of $S$).} Denote
$P^*=\mathrm{Hom}_A(P,A)$. A symmetrical (and slightly extended)
version of \cite[Theorem 2.1]{BrzTurWri:WeakcoaGal}, formulated
for right comodules, is the following.

\begin{theorem} \label{thm:Brzthm}
The canonical map \eqref{eq:can.coring} is bijective and $P^*$ is
a {$T$}-relative projective right $S$-module provided that the
following conditions hold true.
\begin{itemize}
\item[$i)$] The map $P^*\sstac{ T}S\to
\mathrm{Hom}^\cC(P,P^*\sstac{T} P)$, $\xi\sstac{T} s\mapsto \big(\
p\mapsto \xi\sstac{T} s(p)\ \big)$ is an isomorphism (of right
$S$-modules); \item[$ii)$] The {\em lifted canonical map},
\begin{equation}\label{eq:liftcan}
\widetilde{\mathrm{can}}^{T}: P^*\stac{T} P\to \cC,\qquad
\phi\stac{T} p\mapsto \phi(p^{[0]}) p^{[1]}
\end{equation}
is a
split epimorphism of right $\cC$-comodules.
\end{itemize}
\end{theorem}
Motivated by this result, in the present section we investigate
how one can use
{\color{blue}
$(\cM^\hH,{\mathbb R})$-separability of the functor ${\mathbb UV}$ in
\eqref{eq:U} }
to derive properties i) and ii) in Theorem \ref{thm:Brzthm}, for a coring
{\color{blue}
$A\ot_R H$ associated to a right comodule algebra $A$ of a Hopf algebroid
  $\hH$.}

\begin{remark}\label{rem:cond.6.1.i}
In the particular case when the right $\cC$-comodule $P$ in
Theorem \ref{thm:Brzthm} is equal to the base algebra $A$,
property i) reduces to $(A\sstac{T} A)^{co\cC}=A\sstac{T}
A^{co\cC}$. Let us investigate this condition. Note that, for an
$A$-coring $\cC$ possessing a grouplike element $g$, and any
right $T$-module $V$, $V\otimes_{T} A$ is a right $\cC$-comodule via
the comodule structure of $A$. There is an obvious map
$V\otimes_{T} A^{co\cC}\to (V\otimes_{T} A)^{co\cC}$, which is an
isomorphism in appropriate situations: e.g. if $V$ is a flat
${T}$-module, or in the situation described in Lemma
\ref{lem:coinv}. Indeed, in the last case, by applying
Lemma \ref{lem:coinv} to a right $B{:=A^{co\cC}}$-module
$V\sstac{T} \,B$, for a {right $T$}-module $V$, we conclude that
$(V\sstac{T} A)^{co\cC}=V\sstac{T} B$ whenever $A$ is a relative
injective left $\cC$-comodule.
\end{remark}
As it is explained in Appendix
\ref{rem:bijantip}, to a right
comodule algebra $A$ of a Hopf algebroid ${\mathcal H}$ (with
constituent left and right bialgebroids
$(H,L,s_L,t_L,\gamma_L,\pi_L)$ and $(H,R,s_R,t_R,\gamma_R,\pi_R)$,
and a {\em bijective} antipode $S$) one associates two isomorphic
$A$-corings, on the $k$-modules $A\sstac{R} H$ and $H\sstac{R} A$,
and two isomorphic $A^{op}$-corings, on the $k$-modules
$A\sstac{L} H$ and $H\sstac{L} A$. These $A$- and $A^{op}$-corings
are anti-isomorphic, cf. \eqref{eq:phi}. The grouplike element
$1_H$, in the $L$- and $R$-corings underlying ${\mathcal H}$,
determines grouplike elements in all associated $A$- and
$A^{op}$-corings (preserved by the coring (anti-) isomorphisms
\eqref{eq:hgdentwR}, \eqref{eq:hgdentwL} and \eqref{eq:phi}
between them). That is, $A$ (or $A^{op}$) is a right comodule in
each case. Corresponding to the four corings, there are four
canonical maps of the type \eqref{eq:can.coring}, which differ by
the respective coring (anti-) isomorphisms in {Section}
\ref{rem:bijantip}. Following Theorem \ref{thm:schhgd} is
formulated in terms of the $A$-coring $A\sstac R H$ and the
corresponding canonical map \eqref{eq:canRR}.
Certainly, all claims can be reformulated in terms of any of the other three
(anti-) isomorphic corings.

Let $\hH$ be a Hopf algebroid over base algebras $L$ and $R$.
Recall that {a} right $\hH$-comodule algebra $A$ is an $R$-ring.
Assume that the coinvariant subalgebra
{\color{blue}
$B:=A^{co\hH_R}$}
is a $T$-ring (e.g. $T$ is a $k$-subalgebra of $B$). Consider the
lifted version of the canonical map \eqref{eq:canRR}
\begin{equation}\label{eq:can.tilde}
\widetilde{\mathrm{can}}^{T}:A\stac{T} A\to A\stac{R} H,\qquad
a\stac{T} a'\mapsto aa^{\prime [0]}\stac{R} a^{\prime [1]}.
\end{equation}
Note that it is right $L$-linear with respect to the module structures
\begin{equation}\label{eq:rightLmod}
(a\stac {T} a')l=a\stac {T} \pi_R\circ t_L(l) a'\qquad
  \textrm{and}\qquad
(a\stac R h)l=a\stac R t_L(l)h,
\end{equation}
for $a\sstac {T} a'\in A\sstac {T} A$, $a\sstac R h\in A\sstac R
H$ and $l\in L$. Moreover, the lifted canonical map
\eqref{eq:can.tilde} is also left $L$-linear with respect to the
module structures
\begin{equation}\label{eq:leftLmod}
l(a\stac {T} a')= a\pi_R\circ s_L(l) \stac {T} a'\qquad
\textrm{and}\qquad l(a\stac R h)=a\stac R s_L(l)h,
\end{equation}
for $a\sstac {T} a'\in A\sstac {T} A$, $a\sstac R h\in A\sstac R
H$ and $l\in L$.

{\color{blue}
\begin{remark}\label{rem:tcan.split}
Consider a Hopf algebroid $\hH$ and a right $\hH$-comodule algebra $A$. The
lifted canonical map \eqref{eq:can.tilde} is a split epimorphism of right
$L$-modules i.e., using the notations in \eqref{eq:U}, it belongs to ${\mathcal
  E}_{{\mathbb U} {\mathbb V} {\mathbb R}}$ in various situations.

1) If the (right $L$-linear) canonical map \eqref{eq:canRR} is surjective and
$A\ot_R H$ is a projective right $L$-module.
The latter condition holds provided that $H$ is a projective right $L$-module
(via $t_L$) and $A$ is a projective right $R$-module.

2) If the (right $A$-linear) canonical map \eqref{eq:canRR} is surjective,
$A\ot_R H$ is a projective right $A$-module (e.g. $H$ is a projective right
  $R$-module via $s_R$) and ${\mathcal E}_{{\mathbb U}^\cC}\subseteq
  {\mathcal E}_{{\mathbb U} {\mathbb V} {\mathbb R} }$, where ${\mathbb
    U}^\cC$ denotes the forgetful functor $\cM^\hH_A \to \cM_A$.

The condition ${\mathcal E}_{{\mathbb U}^\cC}\subseteq
  {\mathcal E}_{{\mathbb U} {\mathbb V} {\mathbb R}}$ holds whenever
  dealing with a comodule algebra $A$ of a Hopf algebra $\hH$ over a commutative
  ring $k$. Indeed, in this case ${\mathbb V}$ is the identity functor
  $\cM^\hH$, and the functors ${\mathbb U}^\cC$, ${\mathbb R}$
  and ${\mathbb U}$ are forgetful functors. A fourth forgetful functor $\cM_A
  \to \cM_k$ makes the following diagram commutative.
$$
\xymatrix{ \cM^\hH_A\ar[r]^{{\mathbb R}}\ar[d]_{{\mathbb U}^\cC}&
\cM^\hH \ar[d]^{{\mathbb U}}\\
\cM_A\ar[r]&\cM_k }
$$
This proves that in this case ${\mathcal E}_{{\mathbb
U}^\cC}\subseteq
  {\mathcal E}_{{\mathbb U} {\mathbb V} {\mathbb R}}$, thus assumptions
  2) hold e.g. in Schneider's theorem \cite[Theorem I]{Schn:PriHomS}.
\end{remark}
}

{\color{blue}
\begin{lemma}\label{lem:psi}
Consider a right comodule algebra $A$ of a Hopf algebroid $\hH$ as a right
comodule algebra of the constituent right bialgebroid in $\hH$.
Then the following maps are morphisms in $\cM^\hH_A$.
\begin{itemize}
\item[{1)}] the entwining map in \eqref{eq:hgdentwR};
\item[{2)}] the lifted canonical map \eqref{eq:can.tilde}.
\end{itemize}
\end{lemma}

\begin{proof}
$H\ot_R A$ is an object in $\cM^\hH_A$ via the $A$-action induced by the
  multiplication in $A$ and the diagonal
$\hH_R$- and $\hH_L$-coactions.
$A\ot_R H$ is an object in $\cM^\hH_A$ via the $A$-action
$
(a\ot_R h)a' = a a^{\prime [0]} \ot_R h a^{\prime [1]},
$
(where $a\mapsto a^{[0]}\ot_R a^{[1]}$ denotes the $\hH_R$-coaction) and
  $\hH_R$- and $\hH_L$-coactions induced by the respective coproducts in $\hH$.
$A\ot_T A$ is an object in $\cM^\hH_A$ via the relative Hopf module structure
of the second factor.
It is left to the reader to check that both maps in the lemma are compatible
with these structures.
\end{proof}
}

{\begin{lemma}\label{lem:zeta}
Let ${\mathcal H}$ be a Hopf algebroid with constituent left bialgebroid
$(H,L,s_L,t_L,\gamma_L,\pi_L)$, right bialgebroid
$(H,R,s_R,t_R,\gamma_R,\pi_R)$, and a bijective antipode $S$. Let
$A$ be a right {$\hH$-}comodule algebra
{\color{blue}
such that there exists a normalized right $\hH$-comodule map $j:H\to A$}.
Assume that
{\color{blue}
$A^{co{\mathcal H}_R}$}
is a $T$-ring (e.g. $T$ is a subalgebra of
{\color{blue}
$A^{co{\mathcal H}_R}$}). Assume
furthermore that the lifted canonical map \eqref{eq:can.tilde} possesses a
right $L$-module section $\zeta_0^T$ ({with respect to the module
structures \eqref{eq:rightLmod}}).
Then \eqref{eq:can.tilde} possesses a section in $\cM^\hH_A$, given as
\begin{eqnarray*}
\zeta^T: A\stac R H &\to& A\stac T A,\nonumber\\
a\stac R h &\mapsto& \zeta_0^T(1_A\stac R h_{(1)} S^{-1}(a_{[1]})_{(1)})^{[0]}
j\left(S(\zeta_0^T(1_A\stac R h_{(1)} S^{-1}(a_{[1]})_{(1)})^{[1]})h_{(2)}
S^{-1}(a_{[1]})_{(2)}\right)a_{[0]},\nonumber\\
\end{eqnarray*}
where $a\mapsto a_{[0]}\sstac L a_{[1]}$ and $a\mapsto a^{[0]}\sstac R
a^{[1]}$ are the $\hH_L$, and $\hH_R$-coactions in $A$,
respectively, related via
{\color{blue}
\eqref{eq:hgd_coac}},
and $\gamma_L(h)=h_{(1)}\sstac L h_{(2)}$, for $h\in H$.
\end{lemma}
\begin{proof}
{\color{blue}
By Theorem \ref{thm:4.1.s},}
the forgetful functor $\cM^{\mathcal H}\to \cM_L$ is $(\cM^{\mathcal H},{\mathbb
  R})$-separable, where ${\mathbb R}$ is the forgetful functor $\cM^{\mathcal
  H}_A\to \cM^{\mathcal H}$.
By Theorem \ref{thm:rmaschke} 1) this implies that \eqref{eq:can.tilde} is a
split epimorphism in $\cM^\hH$.
Furthermore,
by Lemma \ref{lem:psi}~1),
the object $A\sstac R H\in \cM^\hH_A$ is
isomorphic to $H \sstac R A$. 
{\color{blue}
Since by definition $\cM^\hH_A$ is the category of modules for
  the monad $-\ot_R A:\cM^\hH\to \cM^\hH$, the forgetful functor $\cM^\hH_A\to
  \cM^\hH$ possesses a left adjoint $-\ot_R A:\cM^\hH \to \cM^\hH_A$ (where
  for any right $\hH$-comodule $M$, $M\ot_R A$ is a relative Hopf module via
  the $A$-action on the second factor and the diagonal coactions).
}
Hence $H \sstac R A$ (thus $A \sstac R H$) is  ${\mathcal E}_{{\mathbb
    R}}$-projective by Theorem \ref{teo 4.2.33}. 
This proves that the lifted canonical map
\eqref{eq:can.tilde} is a split epimorphism in $\cM^\hH_A$.

A section can be explicitly constructed as follows.
A right $\hH$-comodule section $\zeta_1^T$ can be constructed using arguments
in the proof of the Rafael type theorem \ref{thm:rrafael} 3). Indeed, in terms
of a right $L$-module section $\zeta_0^T$
and a natural retraction $\nu$ of $\tau_{\mathbb{R}(\bullet)}$
(where $\tau$ denotes the $\hH_L$-coaction, i.e. the unit of the
adjunction of the forgetful functor $\cM^\hH\to \cM_L$ and the
induction functor $\bullet \sstac L H:\cM_L\to \cM^\hH$), one can
put
$$
\zeta_1^T:=\nu_{A\sstac T A} \circ (\zeta_0^{T}\stac L H)\circ
(A\stac R \gamma_L).
$$
The natural retraction $\nu$ was constructed in
{\color{blue}
\eqref{eq:nu_M}}.
Thus a right $A$-module right $\hH$-comodule section of the lifted canonical
map \eqref{eq:can.tilde} is given by
$$
\zeta^T=(A\stac T \mu)\circ(\zeta_1^T\stac R A)\circ (\psi_R\stac R A)\circ
(H\stac R \eta\stac R A)\circ {\psi_R}^{-1},
$$
where $\eta:R\to A$ and $\mu:A\sstac R A \to A$ are unit and multiplication
maps in the $R$-ring $A$, respectively. The map $\zeta^T$ comes out explicitly
as in the claim.
\end{proof}

\begin{theorem} \label{thm:schhgd}
Let ${\mathcal H}$ be a Hopf algebroid with constituent left
bialgebroid $(H,L,s_L,t_L,\gamma_L,\pi_L)$, right bialgebroid
$(H,R,s_R,t_R,\gamma_R,\pi_R)$, and a bijective antipode $S$. Let
$A$ be a right {$\hH$-}comodule algebra and put $B: =
{\color{blue}
A^{co{\mathcal H}_R}}$. {Assume that $B$ is a $T$-ring (e.g. $T$ is a
  $k$-subalgebra of $B$).}
\begin{itemize}
\item[$1)$] If the lifted canonical map \eqref{eq:can.tilde} is a
split epimorphism of right $L$-modules ({with respect to the
  module structures \eqref{eq:rightLmod}}) then the following assertions are
  equivalent.
\begin{itemize}
\item[$(a)$] The canonical map
$\mathrm{can}:A\sstac B A \to A \sstac R H$ in
\eqref{eq:canRR} is bijective and $B$ is a
  direct summand of $A$ as a right $B$-module (hence $A$ is a generator of
  right $B$-modules).
\item[$(b)$]  {\color{blue}
The equivalent conditions in Theorem \ref{thm:4.1.s} hold.}
\item[$(c)$] The functor $^{{\color{blue}co\hH_R}}(\bullet):{}^\hH_A \cM\to {_B
\cM}$ is an
  equivalence, with inverse $A\sstac B\, \bullet$, and $B$ is a
  direct summand of $A$ as a right $B$-module.
\end{itemize}
Furthermore, if these equivalent statements hold then $A$ is a
${T}$-relative projective right $B$-module.
\item[$2)$] If the lifted canonical map \eqref{eq:can.tilde} is a
split epimorphism of left $L$-modules {(with respect to the
  module structures \eqref{eq:leftLmod})}
then the following assertions are equivalent.
\begin{itemize}
\item[$(a)$] The canonical map \eqref{eq:canRR} is bijective and $B$ is a
  direct summand of $A$ as a left $B$-module.
\item[$(b)$] {\color{blue}
The equivalent conditions in Theorem \ref{thm:4.1.s} hold.}
\item[$(c)$] The functor $(\bullet)^{{\color{blue}co\hH_R}}:\cM^\hH_A \to \cM_B$
is an
  equivalence, with inverse $\bullet\sstac B\, A$, and $B$ is a
  direct summand of $A$ as a left $B$-module.
\end{itemize}
Furthermore, if these equivalent statements hold then $A$ is a
${T}$-relative projective left $B$-module.
\end{itemize}
\end{theorem}
\begin{proof}
Recall from {Section} \ref{rem:bijantip} that bijectivity of the
antipode $S$ in the Hopf algebroid ${\mathcal H}$ implies that
also the entwining map \eqref{eq:hgdentwR} is bijective.

$(a)\Rightarrow (b)$
{\color{blue}
In terms of the canonical map \eqref{eq:canRR}, introduce the index notation
$h^{\{1\}}\ot_B h^{\{2\}}:=\can^{-1}(1_A\ot_R h)$ for $h\in H$ 
(implicit summation is understood). 
Using a right $B$-module retraction $p$ of the inclusion $B\to A$, a
normalized right $\hH$-comodule map is given by
$j:H \to A$, $h\mapsto p(h^{\{1\}}) h^{\{2\}}$.
}

$(b)\Rightarrow (a)$
The lifted canonical map \eqref{eq:can.tilde} is a split
epimorphism in $\cM^\hH_A$ by Lemma \ref{lem:zeta}.
{\color{blue}
Hence it is in particular
a split epimorphism of right comodules for the $A$-coring $A\ot_R
H$.}
By considerations in {Section} \ref{rem:entwcorext}, coinvariants
of the right comodules $A$ and  $A\sstac {T} A$ (latter one
defined via the second tensorand) for the $A$-coring $A\sstac R H$
coincide with the
{\color{blue}
$\hH_R$}-coinvariants in $A$ and $A\sstac {T} A$,
respectively. By
{\color{blue}
Theorem \ref{thm:4.1.w}}, $A$ (with coaction $a\mapsto S^{-1}(a_{[1]})\ot_R
a_{[0]}$) is a relative injective left comodule for $\hH_R$.
So  that, by Proposition \ref{prop:cdrinj} $(b)\Rightarrow (a)$,
$A$ is a relative injective left comodule for the $A$-coring
$A\sstac R H$. Taking Remark \ref{rem:cond.6.1.i} into account, it
follows that $(A\sstac {T} A)^{{\color{blue}co\hH_R}}=A\sstac {T} B$, hence all
assumptions in Theorem \ref{thm:Brzthm} hold. Therefore the
canonical map \eqref{eq:canRR} is bijective and $A$ is a
${T}$-relative projective right $B$-module by Theorem
\ref{thm:Brzthm}. It follows by Proposition \ref{prop:cdrinj} that
the right regular $B$-module is a direct summand in $A$.

$(b)\Rightarrow (c)$ By
{\color{blue}
Theorem \ref{thm:4.1.w} and}
Proposition \ref{prop:cdrinj}
  $(b)\Rightarrow (a)$, assertion (b) in the claim implies that $A$ is a
  relative injective right comodule for the $A$-coring
$A\sstac R H$. The $A$-coring $A\sstac R H$
  possesses a grouplike element $1_A\sstac R 1_H$, hence the
unit of the adjunction
{\color{blue}
$(A\sstac B \, \bullet:{}_B \cM\to {}^{\hH_R}_A\cM\ ,\
^{co\hH_R}(\bullet):{}^{\hH_R}_A \cM\to {}_B \cM)$}
(cf. {Section} \ref{preli:coring}) is an
isomorphism, by a left-right symmetric version of Lemma
\ref{lem:coinv} (recall that coinvariants with respect to $\hH$
and $A\sstac R H$ are the same, by arguments in {Section}
\ref{rem:entwcorext}).
{\color{blue}
Then the unit of the adjunction in \ref{prop:coinv} is a natural isomorphism
in light of \ref{prop:two_adj}. }

Let us construct the inverse of the counit,
$$
n_M:A\stac B  ^{{\color{blue}co\hH_R}} M\to M, \qquad a\stac B m\mapsto am,
$$
for $M\in {}^\hH_A \cM$. {Denote the left $\hH_R$-, and
$\hH_L$-coactions on $M$ by
  $m\mapsto m^{[-1]}\sstac R m^{[0]}$ and $m\mapsto m_{[-1]}\sstac L m_{[0]}$,
  respectively.}
The canonical map \eqref{eq:canRR} is bijective by
part (a). Consider the map
\begin{equation}\label{eq:ninv}
M \to A\stac B M, \qquad m\mapsto \can^{-1}(1_A \stac R m^{[-1]}) m^{[0]}.
\end{equation}
By Lemma \ref{lem:zeta} the lifted canonical map
\eqref{eq:can.tilde} has a section $\zeta^T$ in $\cM^\hH_A$.
The map \eqref{eq:ninv} is equal to the composite of
\begin{equation}\label{eq:nlift}
M \to A\stac {T} M, \qquad m\mapsto \zeta^{T}(1_A \stac R
m^{[-1]}) m^{[0]}
\end{equation}
and the canonical epimorphism $A\sstac {T} M\to A\sstac B M$. We
claim that the range of \eqref{eq:nlift} is in $A\sstac {T}
{}^{{\color{blue}co\hH_R}}M$.
{\color{blue}
By Theorem \ref{thm:4.1.s}},
the left $\hH_L$-coaction in $M$ has a
retraction in ${}^\hH\cM$. {The} `{\color{blue} $\hH_R$}-coinvariants part' of
this retraction yields a $k$-linear retraction of the inclusion
${}^{\color{blue}co\hH_R} M\to
{}^{\color{blue}{}^{co\hH_R}} (H\sstac L M)\cong M$.
Explicitly, in terms of a
{\color{blue} normalized right $\hH$-comodule map $j:H\to A$}
(cf. Theorem
{\color{blue}
\ref{thm:4.1.s}}
(a)), we obtain an idempotent map
$$
E_M:M\to {}^{{\color{blue}co\hH_R}} M, \qquad m\mapsto j(m^{[-1]}) m^{[0]}.
$$
Consider the right $L$-module $A$ with action $al:=\pi_R\circ
  t_L(l)a$, for $l\in L$ and $a\in A$. Take a
{\color{blue} normalized right $\hH$-comodule map $j:H\to A$}
as in part (a) of Theorem
{\color{blue}
\ref{thm:4.1.s}}
and introduce a left $B$-module map
$$
P_M:A\stac L M\to M,\qquad a\stac L m\mapsto
a^{[0]}j\big(S(a^{[1]})m^{[-1]}\big) m^{[0]}.
$$
It is well defined by the right $L$-linearity of the right
$\hH_R$-coaction in $A$ and the left $L$-linearity of the left
$\hH_R$-coaction in $M$, and module map properties of $S$ and $j$.
Making use of the relative Hopf module structure of $M$, that is
\eqref{eq:leftentw}, one checks that $E_M\circ P_M =P_M$.
This means that the range of $P_M$ is within
${}^
{\color{blue}
{co\hH_R}}M$. Since the section
$\zeta^T$ in Lemma \ref{lem:zeta} of the lifted
canonical map \eqref{eq:can.tilde} satisfies, for $m\in M$,
$$
\zeta^T(1_A\stac R m^{[-1]})m^{[0]}=\big((A\stac T P_M)\circ (\zeta_0^T\stac L
M)\big)(1_A \stac R m_{[-1]}\stac L m_{[0]}),
$$
the range of \eqref{eq:nlift} is in $A\sstac T
{}^
{\color{blue}
{co\hH_R}}M$. This implies that the range of \eqref{eq:ninv} is
in $A\sstac B {}^
{\color{blue}
{co\hH_R}}
M$. The proof is completed by showing that
the corestriction of \eqref{eq:ninv} to a map ${\widetilde
n}_M:M\to A\sstac B {}^
{\color{blue}
{co\hH_R}}
M$ yields the inverse of $n_M$.
Indeed, since $(A\sstac R\, \pi_R)\circ \can (a\sstac B\,
a')=aa'$, for $a,a'\in A$, and $\can$ is bijective, $n_M \circ
{\widetilde n}_M (m)=\pi_R(m^{[-1]}) m^{[0]}=m$, for any $m\in M$.
On the other hand, since $M$ is an object in ${}^\hH_A \cM$, it
follows by \eqref{eq:leftentw} that  ${\widetilde n}_M\circ
n_M(a\sstac B m)= \big(\can^{-1}(1_A\sstac R
S^{-1}(a_{[1]}))a_{[0]}\big) m$, for $a\sstac B m \in A\sstac B
{}^
{\color{blue}
{co\hH_R}}
 M$. The right $A$-linearity of $\can$ implies that
$\can^{-1}(1_A\sstac R S^{-1}(a_{[1]}))a_{[0]}=a\sstac B 1_A$, for
$a\in A$, which proves ${\widetilde n}_M\circ n_M(a\sstac B
m)=a\sstac B m$.

$(c)\Rightarrow (a)$
{\color{blue}
Observe that $H\sstac R A$ is an object in ${}^\hH_A \cM$, with $A$-action 
$$
a'(h\ot_R a)=h S^{-1}(a'_{[1]}) \ot_R a'_{[0]} a,
$$
where $a\mapsto a_{[0]} \ot_L a_{[1]}$ denotes the $\hH_L$-coaction on $A$,
and $\hH_L$ and $\hH_R$-coactions induced by the respective coproducts.}
The counit of the adjunction
{\color{blue}
$(A\sstac B \,\bullet:{}_B \cM \to {}^\hH_A \cM,
{}^{co\hH_R}(\bullet):{}^\hH_A \cM \to {}_B \cM$)}, evaluated at
  the object $H\sstac R A$, is the isomorphism
$$
n_{H\sstac R A}:A\stac B A\to H\stac R A,\qquad
a\stac B a'\mapsto S^{-1}(a_{[1]})\stac R a_{[0]}a'.
$$
The canonical map \eqref{eq:canRR} is a composite of isomorphisms,
$\can=\psi_R\circ n_{H\sstac R A}$, where $\psi_R$ is the
bijective entwining map \eqref{eq:hgdentwR}. This proves
bijectivity of the canonical map \eqref{eq:canRR}.

In view of
{\color{blue}
Theorem \ref{thm:4.1.s}},
part 2) follows by applying part 1) to the co-opposite Hopf
algebroid ${\mathcal H}_{cop}$ and its right comodule algebra
$A^{op}$,
with $(\hH_R)_{cop}$-coaction $a\mapsto a_{[0]}\sstac
{R^{op}} S^{-1}(a_{[1]})$ and  $(\hH_L)_{cop}$-coaction $a\mapsto a^{[0]}\sstac
  {L^{op}} S^{-1}(a^{[1]})$.
\end{proof}

Observe that ${T}$-relative projectivity and generator
properties of the $B$-module $A$ in Theorem \ref{thm:schhgd} are
as close to its faithful flatness as it is possible for an
arbitrary {base algebra $T$ of $B$}. If $A$ is a projective
$T$-module (e.g. ${T=}k$ is a field) then the equivalent
assertions in part 1) {or} 2) imply that $A$ is a projective
{right or left} $B$-module. Hence the
properties, that $B$ is a direct summand in the right or left
  $B$-module $A$,
in 1) (a) and (c) {or} 2) (a) and (c) imply that $A$ is a
faithfully flat {right or left} $B$-module, cf.
\cite[2.11.29]{Row:rin}.

If $\hH$ is a {\em coseparable} Hopf algebroid (i.e. the underlying $L$-coring
or, equivalently, the underlying $R$-coring is coseparable, cf. \cite[Theorem
3.2]{Bohm:hgdint}) then the forgetful
functors
{\color{blue}
$\cM^{\hH_R}\to \cM_R$ and $\cM^{\hH_L}\to \cM_L$}
are separable (cf. \cite[Corollary 3.6]{Brz:str}).
{\color{blue}
Moreover, also the following proposition holds.}

{\color{blue}
\begin{proposition}\label{prop:cosep}
Let $\hH$ be a Hopf algebroid whose constituent $R$-coring (equivalently, the
constituent $L$-coring) is coseparable. Then the forgetful functors
$\cM^\hH\to \cM^{\hH_R}$ and $\cM^\hH\to \cM^{\hH_L}$ are strict monoidal
isomorphisms.
\end{proposition}

\begin{proof}
The forgetful functors $\cM^\hH\to \cM^{\hH_R}$ and $\cM^\hH\to \cM^{\hH_L}$
are strict monoidal in view of 
\ref{thm:hgd.com.mon}.
For any right $\hH_R$-comodule $(M,\varrho_R)$, the equalizer \eqref{eq:R.M} is
$H\ot_L H$-pure, cf. 
\ref{ex:cosep.pure}.
Symmetrically, also the purity conditions in statement 2) in 
4\ref{thm:pure} hold.
In view of claim 3) in 
\ref{thm:pure}, this proves that the forgetful
functors $\cM^\hH\to \cM^{\hH_R}$ and $\cM^\hH\to \cM^{\hH_L}$ are isomorphisms.
\end{proof}}

{\color{blue}
Proposition \ref{prop:cosep}}
implies that if coseparability of the Hopf algebroid $\hH$ is assumed (not
only relative injectivity of the comodule algebra $A$) then the
assumption in Theorem \ref{thm:schhgd} 1) about
splitting of the lifted canonical map \eqref{eq:can.tilde} as a
right $L$-module map can be replaced by its splitting as a right $R$-module
map. Latter assumption holds in
particular if the canonical map \eqref{eq:canRR} is surjective and
$H$ is projective as a right $R$-module (via the source map of the
constituent right bialgebroid). Indeed, under this assumption,
$A\sstac R H\cong H\sstac R A$ is a projective right $A$-module.
So surjectivity of the right $A$-module map \eqref{eq:canRR}
implies that its lifted version \eqref{eq:can.tilde} is a
retraction of right $A$-modules, and hence of right $R$-modules.
By the separability of the forgetful functor $
{\color{blue} \cM^\hH\cong \cM^{\hH_R}}\to \cM_R$ it follows then that
\eqref{eq:can.tilde} is a retraction of right $\hH$-comodules and
the proof can be completed as in Theorem \ref{thm:schhgd}.
Thus we obtain the following.
\begin{theorem}
Let $\hH$ be a coseparable Hopf algebroid with constituent left bialgebroid
$(H,L,s_L,t_L,$ $ \gamma_L,\pi_L)$, right bialgebroid
$(H,R,s_R,t_R, \gamma_R,\pi_R)$, and a bijective antipode $S$. Let $A$ be a
right comodule algebra and put $B: = A^
{\color{blue}
{co\hH}_R}$. Assume that the
canonical map \eqref{eq:canRR} is surjective.
\begin{itemize}
\item[1)] If $H$ is a projective right $R$-module (via $s_R$) then
  assertions (a), (b) and (c) in Theorem \ref{thm:schhgd} 1) are
  equivalent. Furthermore, if they hold then $A$ is a $k$-relative projective
  {right} $B$-module.
\item[2)] If $H$ is a projective left $R$-module (via $t_R$) then
  assertions (a), (b) and (c) in Theorem \ref{thm:schhgd} 2) are
  equivalent. Furthermore, if they hold then $A$ is a $k$-relative projective
  {left} $B$-module.
\end{itemize}
\end{theorem}

\section{Equivariant injectivity and projectivity}
\label{sec:eqproj}

The notion of equivariant projectivity of a Hopf Galois extension
was introduced in the papers \cite{DaGrHa:Strconn} and
\cite{HaMa:Projmod}. Equivariant projectivity of a Hopf Galois
extension is a crucial property from the non-commutative geometric
point of view, as it turns out to be equivalent to the existence
of a strong connection -- a non-commutative formulation of local
triviality of a principal bundle (see \cite{Haj:Strconn}). {In the
context of Galois extensions $B\subseteq A$ by corings (or
bialgebroids or Hopf algebroids), equivariant projectivity {\em
relative} to some subalgebra of $B$ was shown to be equivalent to
the existence of more general strong connections in the paper
\cite{BohmBrz:relchern}.}

In this section we look for conditions on a Galois extension by a
Hopf algebroid, under which it obeys {(relative)} equivariant
injectivity and projectivity properties.
Recall that having a Hopf algebra $H$ over a commutative ring $k$
and a right {$H$}-comodule algebra $A$, which is a relative
injective right $H$-comodule, $A$ was shown to be a
$B(=A^{coH})$-equivariantly injective $H$-comodule in
\cite[Theorem 5.6]{SchSch:genGal}. (This result is extended to
{algebra
  extensions by Hopf algebroids} in Theorem \ref{thm:eqinj}
below.) What is more, using the proven $B$-equivariant injectivity
of a relative injective $H$-comodule algebra $A$, it was also
shown in \cite[Theorem 5.6]{SchSch:genGal} that the $B$-module $A$
is $H$-equivariantly projective if and only if it is $k$-relative
projective. If $A$ is a relative injective right comodule algebra
of a Hopf algebra $H$ with a bijective antipode, with coinvariants $B$, and
the lifted canonical map \eqref{eq:can.lifted} is a split epimorphism of
$k$-modules, then $A$ is an $H$-Galois extension of $B$
and the $B$-module $A$ is relative projective (cf. Theorem
\ref{thm:schhgd}). Hence the $B$-module $A$ is also
$H$-equivariantly projective by the quoted result in \cite[Theorem
5.6]{SchSch:genGal}. A most naive generalization of this result to
Hopf algebroid Galois extensions seems not to hold. The reason is
that -- if working with a Hopf algebroid $\hH$ over different
non-commutative base algebras $L$ and $R$ -- relative projectivity
of the $B$-module $A$ is not enough to prove its {(relative)}
$\hH$-equivariant projectivity. One needs more: {(relative)}
$L$-equivariant projectivity
(see Theorem \ref{thm:eqproj} below).
As a matter of fact, for a relative injective right comodule algebra $A$ of
a Hopf algebroid with a bijective antipode, with coinvariants $B$,
we were not able to deduce {relative} $L$-equivariant projectivity
of the $B$-module $A$ form the {splitting}
of the lifted canonical map \eqref{eq:can.tilde} as a left, or right
$L$-module map, as assumed in Theorem \ref{thm:schhgd}. We needed
a stronger assumption: splitting of the lifted canonical map
\eqref{eq:can.tilde} as an $L$-$L$ bimodule map (see Proposition
\ref{prop:Leqproj} below).

\begin{definition}\label{def:eq_proj&inj}
Let $\cD$ be an $L$-coring and $B$ a $T$-ring. A left $B$-module
and right $\cD$-comodule $V$, with left $B$-linear right
$\cD$-coaction, is called a {\em {$T$-relative}
$\cD$-equivariantly projective} left $B$-module if the left action
$B\sstac {T} V \to V$ is an epimorphism split by a left $B$-module
right $\cD$-comodule map. We call $V$ a {\em $B$-equivariantly
injective} right $\cD$-comodule if the right coaction $V \to
V\sstac L \cD$ is a monomorphism split by a left $B$-module right
$\cD$-comodule map.

Analogous notions for right $B$-modules and left $\cD$-comodules, with a right
$B$-linear left $\cD$-coaction, are defined symmetrically.
\end{definition}
Considering an algebra $L$ as a trivial $L$-coring $\mathcal{L}$,
a {$T$-relative} $\mathcal{L}$-equivariantly projective left
$B$-module $V$ is called simply {$T$-relative} $L$-equivariantly
projective. Clearly, for an $L$-coring $\cD$ and a {$T$-ring $B$},
a {$T$-relative} $\cD$-equivariantly projective left $B$-module
$V$ is necessarily {$T$-relative} $L$-equivariantly projective.

For a (left or right) bialgebroid $\mathcal{S}$ over an algebra
$L$, and a {$T$-ring $B$}, a {$T$-relative}
$\mathcal{S}$-equivariantly projective $B$-module means a
$B$-module which is {$T$-relative} equivariantly projective for
the $L$-coring underlying $\mathcal{S}$.

{\color{blue}
In the same spirit, relative $\hH$-equivariant projectivity for a Hopf
algebroid $\hH$ can be introduced:

\begin{definition}\label{def:6.1}
Consider a Hopf algebroid $\hH$ and a $T$-ring $B$. A left $B$-module and right
$\hH$-comodule $V$, such that the left $B$-action on $V$ is a right
$\hH$-comodule map, is said to be {\em $T$-relative $\hH$-equivariantly
  projective} if the action $B\ot_T V\to V$ is an epimorphism split by a left
$B$-linear and right $\hH$-colinear map.
\end{definition}
}

The following two theorems extend \cite[Theorem 5.6]{SchSch:genGal} to
non-commutative base algebras.

\begin{theorem}\label{thm:eqinj}
Let $\hH$ be a Hopf algebroid with constituent left bialgebroid
$\hH_L=(H,L,s_L,t_L,\gamma_L,\pi_L)$, right bialgebroid
$\hH_R=(H,R,s_R,t_R,\gamma_R,\pi_R)$, and antipode $S$. Let $A$ be
a right {$\hH$-}comodule algebra and denote $B:= A^
{\color{blue}
{co\hH_R}}$.
{\color{blue}
If the equivalent conditions in Theorem \ref{thm:4.1.s} hold then the
$\hH_L$-coaction on $A$ possesses a left $B$-linear and right $\hH$-colinear
retraction.}
\end{theorem}

\begin{proof}
Using a method in \cite[Lemma 4.1]{SchSch:genGal},
one constructs a left $B$-linear right $\hH$-colinear retraction $\phi$ of the
$\hH_L$-coaction $a\mapsto  a_{[0]}\sstac L a_{[1]}$ on $A$, in terms of an
$\hH$-colinear retraction $\nu$. Explicitly,
$$
\phi:A\stac L H \to A, \qquad
a\stac L h \mapsto a^{[0]}\,\,\nu\!\left(1_A\stac L S(a^{[1]})h\right),
$$
where $a\mapsto a^{[0]}\sstac R a^{[1]}$ denotes the right
$\hH_R$-coaction on $A$.
Note that since $\nu$ is an $\hH$- comodule map
it is left $R$-linear and hence the map $\phi$ is well defined.
\end{proof}

\begin{theorem}\label{thm:eqproj}
Let $\hH$ be a Hopf algebroid with constituent left bialgebroid
$\hH_L=(H,L,s_L,t_L,\gamma_L,\pi_L)$, right bialgebroid
$\hH_R=(H,R,s_R,t_R,\gamma_R,\pi_R)$, and antipode $S$. Let $A$ be
a right {$\hH$-}comodule algebra and $B:=A^
{\color{blue}
{co\hH_R}}$.
Let $T$ be an algebra such that $B$ is a $T$-ring (e.g. $T$ is some
$k$-subalgebra of $B$).
{\color{blue}
If the right $\hH_L$-coaction on $A$ possesses a left $B$-linear and right
$\hH$-colinear retraction
then the following assertions are equivalent.
\begin{itemize}
\item[{(a)}] $A$ is a $T$-relative $\hH$-equivariantly projective left
  $B$-module.
\item[{(b)}] $A$ is a $T$-relative $L$-equivariantly projective left
  $B$-module.
\end{itemize}}
\end{theorem}
\begin{proof}
If $A$ is a {$T$-relative} $\hH$-equivariantly projective left
$B$-module then it is obviously {$T$-relative} $L$-equivariantly
projective.
In order to see the converse implication, take a $B$-$L$ bimodule section
$\chi_0^{T}$ of the multiplication map
$B\sstac {T} A\to A$ and a left $B$-linear and right $\hH$-colinear
retraction $\phi$ of the right $\hH_L$-coaction
$\tau_A:A\to A\sstac L H$ in $A$.
It determines a left $B$-linear and $\hH$-colinear map
$$
\chi^{T} : = (B\stac {T} \phi)\circ (\chi_0^T\stac L
H)\circ \tau_A:A \to B \stac {T} A.
$$
It follows by the left $B$-linearity of $\phi$ that $\chi^{T}$ is
a section of the multiplication map $B\sstac {T} A\to A$.
\end{proof}

The message of Theorem \ref{thm:eqinj} and Theorem
\ref{thm:eqproj} is to look for situations, in which the
{$T$-relative} $L$-equivariant projectivity condition in Theorem
\ref{thm:eqproj}
(b) holds, for a right comodule algebra
{\color{blue}
of a Hopf algebroid, obeying the conditions in Theorem \ref{thm:4.1.s}}.

It is discussed in Appendix \ref{preli:left-right} that if the antipode of
a Hopf algebroid $\hH$ {(over base algebras $L$ and $R$)} is
bijective then a right comodule algebra $A$ has a canonical left
$\hH^{op}$-comodule algebra structure, with coactions
$a\mapsto S^{-1}(a_{[1]})\sstac R a_{[0]}$
$a\mapsto S^{-1}(a^{[1]})\sstac L a^{[0]}$.
Recall that this left $\hH_L$-coaction
corresponds to the left $L$-module
structure of $A$, which is related to its right $R$-module
structure via
\begin{equation}\label{eq:leftL_mod}
la=a\pi_R\circ s_L(l), \qquad \textrm{for } l\in L, \ a\in A.
\end{equation}
Since the right actions in $A$ by $R$ and $B:=A^
{\color{blue}
{co\hH_R}}$ commute
(cf. Section \ref{preli:bgd.comodule}), $A$ is an $L$-$B$
bimodule via the left $L$-action \eqref{eq:leftL_mod} and the
obvious right $B$-action. The following proposition {concerns}
equivariant projectivity of this $L$-$B$ bimodule $A$.
\begin{proposition}\label{prop:Leqproj}
Let $\hH$ be a Hopf algebroid with constituent left bialgebroid
$\hH_L=(H,L,s_L,t_L,\gamma_L,$ $\pi_L)$, right bialgebroid
$\hH_R=(H,R,s_R,t_R,\gamma_R,\pi_R)$, and a bijective antipode
$S$. Let $A$ be a right $\hH$-comodule algebra.
Set $B:= A^
{\color{blue}
{co\hH_R}}$, {and assume that $B$ is a $T$-ring (e.g. $T$ is a
  $k$-subalgebra of $B$)}. Assume that the lifted canonical map
\eqref{eq:can.tilde} is a split epimorphism of $L$-$L$ bimodules
(with respect to the module structures
\eqref{eq:rightLmod} and \eqref{eq:leftLmod}).
If
{\color{blue}
the equivalent conditions in Theorem \ref{thm:4.1.s} hold}
then $A$ is a $T$-relative $L$-equivariantly projective right $B$-module.
\end{proposition}
\begin{proof}
Let $\zeta_0^T$ be an $L$-$L$ bimodule section of the lifted
canonical map \eqref{eq:can.tilde}. By Lemma \ref{lem:zeta}
the map \eqref{eq:can.tilde} is split by the right $A$-module and right
$\hH$-comodule map $\zeta^T$, explicitly given in Lemma \ref{lem:zeta}. From
the proof of implication $(b)\Rightarrow (a)$ in Theorem \ref{thm:schhgd} we
have that $(A\sstac {T} A)^
{\color{blue}
{co\hH_R}}
=A\sstac {T} B$. Hence taking
the `{\color{blue}$\hH_R$}-coinvariants part' of $\zeta^{T}$, we obtain a right
$B$-module section of the multiplication map $A\sstac {T} B\to A$,
$$
\chi_0^{T}:A\to A\stac {T} B,\qquad a\mapsto \zeta_0^{T}(1_A\stac
R S^{-1}(a_{[1]})_{(1)})^{[0]} j\left(S(\zeta_0^{T}(1_A\stac R
S^{-1}(a_{[1]})_{(1)})^{[1]}) S^{-1}(a_{[1]})_{(2)}\right)a_{[0]}.
$$
Consider the left $L$-module structure \eqref{eq:leftL_mod} of
$A$. The right $\hH_L$-coaction in $A$ is left $L$-linear in the
sense that $(a \pi_R\circ s_L(l))_{[0]}\sstac L (a \pi_R\circ
s_L(l))_{[1]}= a_{[0]}\sstac L a_{[1]} s_R\circ \pi_R\circ
s_L(l)$, for $l\in L$ and $a\in A$. The antipode satisfies
$S^{-1}(h s_R\circ \pi_R\circ s_L(l))= t_R\circ \pi_R\circ s_L(l)
S^{-1}(h)=s_L(l) S^{-1}(h)$, for $l\in L$ and $h\in H$. The
coproduct $\gamma_L$ is left $L$-linear, i.e. $(s_L(l)
h)_{(1)}\sstac L (s_L(l) h)_{(2)} =s_L(l) h_{(1)}\sstac L
h_{(2)}$, for $l\in L$ and $h\in H$. The map $\zeta_0^{T}$ is left
$L$-linear by assumption, with respect to the left $L$-module
structures in \eqref{eq:leftLmod}. The right $\hH_R$-coaction in
$A\sstac {T} A$ is given via the second factor, so it is obviously
left $L$-linear with respect to the left $L$-module structure in
\eqref{eq:leftLmod}. All these considerations together verify the
left $L$-linearity of $\chi_0^{T}$ with respect to the left
$L$-module structure \eqref{eq:leftL_mod} in $A$. Hence
$\chi_0^{T}$ is an $L$-$B$ bilinear section of the multiplication
map $A\sstac {T} B\to A$, which proves {$T$-relative}
$L$-equivariant projectivity of the right $B$-module $A$.
\end{proof}

Let $A$ be a right comodule algebra of a Hopf algebroid $\hH$ with a
bijective antipode $S$.
{\color{blue}
Recall that in this case the conditions (a)-(c) in Theorem \ref{thm:4.1.s} are
equivalent also the the analogous conditions for the right
$\hH_{cop}$-comodule algebra $A^{op}$, with $(\hH_R)_{cop}$-coaction $a\mapsto
a_{[0]}\ot_{R^{op}} S^{-1}(a_{[1]})$ and $(\hH_L)_{cop}$-coaction $a\mapsto
a^{[0]}\ot_{L^{op}} S^{-1}(a^{[1]})$.
}
Hence Proposition \ref{prop:Leqproj}
can be applied to the right $\hH_{cop}$- comodule algebra $A^{op}$.
It yields a result about the equivariant projectivity of $A$
as a $B$-$L$ bimodule, with obvious left $B$-action, and right $L$-action
related to the left $R$-action via
$$
al=\pi_R\circ t_L(l) a, \qquad \textrm{for } l\in L, \ a\in A.
$$
\begin{corollary}\label{cor:eqproj_left}
In the setting of Proposition \ref{prop:Leqproj}, $A$ is a
{$T$-relative}  $L$-equivariantly projective left $B$-module.
\end{corollary}
The following corollary is the main result of this section. It
formulates sufficient conditions on a Galois extension $B\subseteq
A$ by a Hopf algebroid $\hH$ with a bijective antipode, under
which $A$ is a {$T$-relative} $\hH$-equivariantly projective left
and right $B$-module, for an algebra $T$ such that $B$ is a
  $T$-ring.

\begin{corollary}\label{cor:eqproj}
Let $\hH$ be a Hopf algebroid with constituent left bialgebroid
$\hH_L=(H,L,s_L,t_L,\gamma_L,$ $\pi_L)$, right bialgebroid
$\hH_R=(H,R,s_R,t_R,\gamma_R,\pi_R)$, and a bijective antipode
$S$. Let $A$ be a right {$\hH$-}comodule algebra.
Denote $B: =A^
{\color{blue}
{co\hH_R}}$
{and assume that $B$ is a $T$-ring (e.g. $T$ is a subalgebra of $B$)}. Assume
that the lifted canonical map
\eqref{eq:can.tilde} is a split epimorphism of $L$-$L$ bimodules (with respect
to the module structures \eqref{eq:rightLmod} and \eqref{eq:leftLmod}).
{\color{blue}
If the equivalent conditions in Theorem \ref{thm:4.1.s} hold then $A$ is a
$T$-relative $\hH$-equivariantly projective left and right
$B$-module. Moreover, in this case $B\subseteq A$ is a right $\hH_R$-Galois
extension.}
\end{corollary}

\begin{proof}
The Galois property, i.e. bijectivity of the canonical map \eqref{eq:canRR},
 follows by virtue of Theorem \ref{thm:schhgd} $(b)\Rightarrow (a)$.

{\color{blue}
The right $\hH_L$-coaction on $A$ possesses a left $B$-module right
$\hH$-comodule retraction and the left $\hH_L$-coaction on $A$ possesses a
right $B$-module left $\hH$-comodule retraction},
by Theorem \ref{thm:eqinj}, and its application to the co-opposite Hopf
  algebroid $\hH_{cop}$ and the right {$\hH_{cop}$-}comodule algebra
  $A^{op}$,
respectively. By Proposition \ref{prop:Leqproj}, $A$ is
  a {$T$-relative} $L$-equivariantly projective right $B$-module. By
  Corollary \ref{cor:eqproj_left}, $A$ is a {$T$-relative}
  $L$-equivariantly projective left $B$-module. Hence $A$ is a {$T$-relative}
  $\hH$-equivariantly projective left and right $B$-module by
  Theorem \ref{thm:eqproj}, and its application to the co-opposite Hopf
  algebroid $\hH_{cop}$ and the right {$\hH_{cop}$-}comodule algebra
  $A^{op}$,
respectively.
\end{proof}

By \cite[Theorem 3.7]{BohmBrz:relchern} we conclude that there
exists a strong {$T$-}connection for an extension $B\subseteq A$
as in Corollary \ref{cor:eqproj} whenever $T$ is a
$k$-subalgebra of $B$. {In \cite[Theorem 5.14]{BohmBrz:relchern}
conditions are formulated for the independence of the
corresponding relative Chern-Galois character of the choice of a
strong $T$-connection. Note that in the case when in Corollary
\ref{cor:eqproj} the $k$-algebra $T$ is equal to $k$, these
conditions reduce to the assumption that} $A$ is a locally
projective $k$-module.

{\begin{example} {\em Cleft} extensions by Hopf algebroids were
introduced in \cite[Definition
  3.5]{BohmBrz:cleft}, as follows. Let $\hH$ be a Hopf algebroid with
constituent left bialgebroid $\hH_L=(H,L,s_L,t_L,\gamma_L,$ $\pi_L)$, right
bialgebroid $\hH_R=(H,R,s_R,t_R,\gamma_R,\pi_R)$, and antipode $S$. Let $A$ be
a right $\hH$-comodule algebra with coinvariants $B:=A^
{\color{blue}
{co\hH_R}}$. Denote the
unit map of the corresponding $R$-ring $A$ by $\eta_R:R\to A$. The algebra
extension $B\subseteq A$ is called {\em $\hH$-cleft} provided that the
following conditions hold.
\begin{itemize}
\item[a)] $A$ is an $L$-ring (with some unit map $\eta_L:L\to A$) and $B$ is
  an $L$-subring of $A$.
\item[b)] There exist morphisms $j\in {}_L\mathrm{Hom}^\hH(H,A)$ and
  ${\widetilde j}\in {}_R\mathrm{Hom}_L(H,A)$, satisfying
$$
\mu\circ({\widetilde j}\stac L j)\circ \gamma_L=\eta_R\circ \pi_R \qquad
\textrm{and}\qquad
\mu\circ(j\stac R {\widetilde j})\circ \gamma_R=\eta_L\circ \pi_L,
$$
where $\mu$ denotes the multiplication in $A$, both as an $L$-ring and as an
$R$-ring. The bimodule structures in $H$ are given by
$$
lhr:=s_L(l) h s_R(r)\qquad  \textrm{and}\qquad
rhl=t_L(l)ht_R(r),\qquad \textrm{for }l\in L, \ r\in R,\ h\in H.
$$
The bimodule structures in $A$ are given by
$$
lar:=\eta_L(l) a \eta_R(r)\qquad  \textrm{and}\qquad
ral=\eta_R(r)a\eta_L(l),\qquad \textrm{for }l\in L, \ r\in R,\ a\in A.
$$
\end{itemize}

In an $\hH$-cleft extension $B\subseteq A$
the map ${\widetilde j}(1_H)j(-):H\to A$ is right $\hH$-colinear and
normalized. Hence
{\color{blue}
the equivalent conditions in Theorem \ref{thm:4.1.s} hold}.
By definition $B$ is an
$L$-ring. The lifted canonical map
$$
{\widetilde{\mathrm{can}}}^L:A\stac L A\to A\stac R H,\qquad a\stac
L a'\mapsto aa^{\prime[0]}\stac R a^{\prime[1]}
$$
possesses an $L$-$L$ bilinear section (with respect to the module structures
\eqref{eq:rightLmod} and \eqref{eq:leftLmod}):
\begin{equation}\label{eq:zeta_L}
\zeta_0^L:A\stac R H \to A\stac L A,\qquad a\stac R h  \mapsto a{\widetilde
  j}(h_{(1)})\stac L j(h_{(2)}).
\end{equation}
The map \eqref{eq:zeta_L} is well defined by the module map properties of $j$
and ${\widetilde j}$. It is left $L$-linear by the identity ${\widetilde
  j}(t_R(r)h)= {\widetilde j}(h)\eta_R(r)$, for $r\in R$ and $h\in H$, see
\cite[Lemma 3.7]{BohmBrz:cleft}. Right $L$-linearity of
\eqref{eq:zeta_L} follows by the left $R$-linearity of (the right
$\hH$-comodule map) $j$, i.e. $j(s_R(r)h)= \eta_R(r) j(h)$, for
$r\in R$ and $h\in H$.
In view of Corollary \ref{cor:eqproj}, all these considerations
together imply that a cleft extension $B\subseteq A$ by a Hopf
algebroid $\hH$ with a bijective antipode is an $\hH_R$-Galois
extension which is an $L$-relative $\hH$-equivariantly projective
left and right $B$-module, cf. \cite[Lemma 5.1]{BohmBrz:cleft}.
\end{example}

\appendix
\section{Coring extensions, entwining structures and Hopf algebroids}
\label{sec:prelims}

\begin{claim}\label{preli:A.ring}
For a $k$-algebra $A$, an {\bf $A$-ring} $T$ means an algebra (or monoid) in
the monoidal category of $A$-$A$ bimodules. More explicitly, it consists of an
$A$-$A$ bimodule $T$, equipped with a bilinear associative product
$\mu:T\sstac A T\to T$ and a bilinear unit map $\eta:A\to T$. An $A$-ring $T$
is equivalent to a $k$-algebra $T$ and a $k$-algebra map $\eta:A\to T$.

For an $A$-ring $(T,\mu,\eta)$, the opposite means the
$A^{op}$-ring
  $T^{op}$, with $A^{op}$-$A^{op}$ bimodule structure
$$
A^{op}\otimes T\otimes A^{op}\to T,\qquad a\otimes t\otimes a'\mapsto a'ta,
$$
product $t\sstac {A^{op}} t'\mapsto t' t$ and unit $\eta:A^{op}\to T^{op}$.

An $A$-ring $T$ determines a monad $\bullet \sstac A T$ on
$\cM_A$. By right $T$-modules we mean algebras for this monad.
This notion coincides with that of right modules for the
$k$-algebra $T$. Left modules are defined symmetrically.
\end{claim}

\begin{claim}\label{preli:coring}
A {\bf coring} over an algebra $A$ means a coalgebra (or comonoid) in
the monoidal category of $A$-$A$ bimodules. More explicitly, it consists of an
$A$-$A$ bimodule $\cC$, equipped with a bilinear coassociative coproduct
$\Delta:\cC\to \cC\sstac A \cC$ and a bilinear counit map $\epsilon:\cC\to A$.
This extends the notion of a coalgebra.

For an $A$-coring $(\cC,\Delta,\epsilon)$, the co-opposite means the
$A^{op}$-coring $\cC_{cop}$, with $A^{op}$-$A^{op}$ bimodule $\cC$,
$$
A^{op}\otimes \cC\otimes A^{op}\to \cC,\qquad
a\otimes c \otimes a'\mapsto a'ca,
$$
{coproduct} $\Delta^{op}:c\mapsto c^{(2)}\sstac {A^{op}} c^{(1)}$
and counit $\epsilon$.

An $A$-coring $\cC$ determines a comonad $\bullet \sstac A \cC$ on $\cM_A$. By
right $\cC$-comodules we mean coalgebras for this comonad. That is, a right
$\cC$-comodule is a right $A$-module $M$, equipped with a right $A$-linear
coassociative and counital coaction $\varrho^M$.
Right $\cC$-comodule maps are right $A$-linear maps which are compatible with
the coactions.
Left $\cC$-comodules are defined symmetrically.

The forgetful functor $\cM^\cC\to \cM_A$ possesses a right
adjoint, the functor $\bullet\sstac A \cC:\cM_A\to \cM^\cC$. The
unit of the adjunction is given by the right coaction
$\varrho^M:M\to M\sstac A \cC$, for $M\in \cM^\cC$, and the counit
of the adjunction is given in terms of the counit $\epsilon$ of
$\cC$ as $N\sstac A \epsilon:N\sstac A \cC\to N$, for $N\in
\cM_A$, cf. \cite[18.13 (2)]{BrzWis:cor}.

A right comodule $M$ for an $A$-coring $\cC$ is called {\bf relative
injective} if any $\cC$-comodule map of domain $M$, which is a
split monomorphism of $A$-modules, is a split monomorphism
of $\cC$-comodules, too. By \cite[18.18]{BrzWis:cor},
$M$ is a relative injective $\cC$-comodule if and only if the coaction
$\varrho^M$ is a section of $\cC$-comodules.

By \cite[Lemma 5.1]{Brz:str}, the right regular $A$-module extends to a
comodule for an $A$-coring $\cC$ if and only if there exists a {\bf grouplike
element} $g$ in $\cC$ (meaning that $\Delta(g)=g\sstac A g$ and
$\epsilon(g)=1_A$). A bijective correspondence between grouplike elements $g$
in $\cC$ and right $\cC$-coactions $\varrho^A$ in $A$ is given by $g\mapsto
(\varrho^A: a\mapsto ga)$. A similar equivalence holds between grouplike
elements and left $\cC$-comodule structures in $A$.

For an $A$-coring $\cC$ with a grouplike element $g$, the {\bf coinvariants}
with respect to $g$ of a right $\cC$-comodule $M$ are defined as the elements
of the set
$$
M^{co\cC}=\{\ m\in M\ |\ \varrho^M(m)=m\stac A g\ \}.
$$
Coinvariants of left $\cC$-comodules are defined symmetrically. In
particular, the coinvariants of $A$, both as a right
$\cC$-comodule and a left $\cC$-comodule, are the elements of the
subalgebra
$$
B=\{\ b\in A\ |\ gb=bg\ \}.
$$
A grouplike element $g$ in $\cC$ determines an adjunction 
$(\bullet\sstac B A,(\bullet)^{co\cC})$, 
between the categories $\cM_B$ and $\cM^\cC$. The unit
and counit are given by the maps
\begin{eqnarray}
&u_N:N\to (N\stac B A)^{co\cC},\qquad &x\mapsto x\stac B
  1_A, \quad
\textrm{and} \label{eq:co_un}\\
& n_M: M^{co\cC}\stac B A\to M,\qquad &y\stac B a \mapsto ya,
\label{eq:co_coun}
\end{eqnarray}
respectively, for any $N\in \cM_B$ and $M\in \cM^\cC$, cf.
\cite[28.8]{BrzWis:cor}. There is a symmetrical adjunction between
the categories ${}_B \cM$ and ${}^\cC\cM$.

An $A$-coring $\cC$ with a grouplike element $g$ is called a {\bf Galois
coring} if the canonical map
\begin{equation}\label{eq:cor_can}
\mathrm{can}:\ A\stac B A \to \cC,\qquad a\otimes a'\mapsto aga'
\end{equation}
is bijective. For more information about corings we refer to the
{monograph} \cite{BrzWis:cor}.
\end{claim}

\begin{claim}\label{preli:cor.ext}
Let $\cD$ be a coring over a base $k$-algebra $L$ and $\cC$ a
coring over a $k$-algebra $A$. Assume that $\cC$ is a
$\cC$-$\cD$ bicomodule with the left regular $\cC$-coaction
$\Delta_\cC$ and some right $\cD$-coaction $\tau_\cC$. By
definition \cite[22.1]{BrzWis:cor}, this means that $\tau_\cC$ is
left $A$-linear (hence $\cC\otimes_{A} \cC$ is also a right
$\cD$-comodule with coaction $\cC\otimes_{A}\tau_\cC$) and the
coproduct $\Delta_\cC$ is right $\cD$-colinear. Equivalently, the
coproduct $\Delta_\cC$ is right $L$-linear (hence
$\cC\otimes_{L}\cD$ is a left $\cC$-comodule with coaction
$\Delta_\cC\otimes_{L}\cD$) and the $\cD$-coaction $\tau_\cC$ is
left $\cC$-colinear. {In this case, following
\cite[Definition 2.1]{Brz:Corext}, we say that}
$\cD$ is a {\bf right extension} of $\cC$.
{\color{blue}
For a right extension $\cD$ of $\cC$, assume that the equalizer
\begin{equation}\label{eq:eq}
\xymatrix{
M \ar[rr]^-{\varrho^M} &&
M\ot_A \cC \ar@<2pt>[rr]^-{\varrho^M\ot_A \cC}\ar@<-2pt>[rr]_-{M\ot_A
  \Delta_\cC}&&
M\ot_A \cC\ot_A \cC
}
\end{equation}
in $\cM_L$ is $\cD\ot_L \cD$-pure, i.e. it is preserved by the functor $-\ot_L
{\mathcal D} \ot_L \cD:\cM_L\to \cM_L$, for any right $\cC$-comodule
$M$. If this condition holds, we say that $\cD$ is a {\bf pure}
right coring extension of $\cC$.
By \cite[22.3]{BrzWis:cor}
and its Erratum, for pure coring extension $\cD$ of $\cC$, 
there is a functor ${\mathbb R}:= \bullet \Box_\cC \cC:\cM^\cC\to \cM^\cD$,
given by a cotensor product, 
that renders diagram \eqref{fig:Brzfunc} commutative (up to the
  natural isomorphism $M\cong M\Box_\cC\, \cC$, for $M\in {\mathfrak M}^\cC$). 
The explicit form of the functor ${\mathbb R}$ is computed in \cite[Theorem
  2.6]{Brz:Corext} (mind the missing purity condition in the journal
version):} Using the right $\cD$-coaction $\tau_\cC:c\mapsto
c_{[0]}\otimes_{L}c_{[1]}$, for $c\in \cC$ (note our convention to use
character $\tau$ for $\cD$-coactions and lower indices of the Sweedler type
to denote components of the coproduct and coactions of $\cD$), any right
$\cC$-comodule $M$ is equipped with a right $\cD$-comodule structure with
right $L$-action 
\begin{equation}
ml: = m^{[0]}\epsilon_\cC({m^{[1]}}l),
\qquad \textrm{for } m\in M \textrm{ and } l\in L,
\end{equation}
and $\cD$-coaction
\begin{equation}\label{eq:D_coac}
\tau_M\ :\ M\to M\stac{L} \cD, \qquad m\mapsto m_{[0]}\stac{L} m_{[1]}: =
m^{[0]}\epsilon_\cC({m^{[1]}}_{[0]})\stac{L} {m^{[1]}}_{[1]},
\qquad \textrm{for } m\in M,
\end{equation}
where $\varrho^M:m\mapsto m^{[0]}\otimes_{A} m^{[1]}$ denotes the
$\cC$-coaction on $M$ (note our convention to use character
$\varrho$ for $\cC$-coactions and upper indices of the Sweedler
type to denote components of the coproduct and coactions of $\cC$). 
With this definition any right $\cC$-comodule
map is $\cD$-colinear. In particular, a right $\cC$-coaction, being
$\cC$-colinear by coassociativity, is $\cD$-colinear.
\end{claim}

{\color{blue}
\begin{claim}\label{ex:cosep.pure}
{\bf Any right coring extension of a coseparable $A$-coring $\cC$ is
pure.} Indeed, \eqref{eq:eq} is a split equalizer in $\cM_A$ (split by the
right $A$-module map $M\ot_A \cC\ot_A \varepsilon_\cC$). By
separability of the functor $\cM^\cC \to \cM_A$, it is a split equalizer also
in $\cM^\cC$.
If $\cD$ is an $L$-coring that is a right extension of $\cC$, then taking
cotensor products with the $\cC$-$\cD$ bicomodule $\cC$ defines a functor
$-\Box_{\cC} {\cC}:\cM^\cC \to \cM_L$, equipping any right
$\cC$-comodule $M\cong M\Box_\cC \cC$ with a right $L$-action.
By right $L$-linearity of any $\cC$-comodule map, splitting of the equalizer
\eqref{eq:eq} in $\cM^\cC$ implies that it splits in also in $\cM_L$. Hence
the purity condition holds.
\end{claim}}

\begin{claim}\label{rem:entwcorext}
An {\bf entwining structure} over a ({not necessarily
commutative}) algebra $L$ consists of an $L$-ring $A$, with
multiplication $\mu$ and unit $\eta$, an $L$-coring $\cD$, with
comultiplication $\Delta$ and counit $\epsilon$, and an $L$-$L$
bilinear map $\psi: \cD\sstac L A \to A \sstac L\cD$, satisfying
the following compatibility conditions.
\begin{eqnarray*}
&\psi\circ (\cD \stac L \mu)=(\mu \stac L \cD)\circ (A\stac L \psi)\circ
  (\psi   \stac L A)\qquad
&\psi\circ (\cD \stac L \eta)=\eta \stac L \cD\\
&(A\stac L \Delta)\circ \psi =(\psi\stac L \cD)\circ (\cD \stac L \psi)\circ
  (\Delta \stac L A)\qquad
&(A\stac L \epsilon)\circ \psi=\epsilon \stac L A.
\end{eqnarray*}
In complete analogy with \cite[Proposition 2.2]{Brz:str}, $A\sstac L\cD$ is an
$A$-coring. Its bimodule structure is given by
$$
a_1(a\stac L d)a_2=a_1 a \psi(d\stac L a_2), \qquad \textrm{for }a_1,a_2\in A,
\ a\stac L d\in A\stac L \cD.
$$
The coproduct is equal to $A\sstac L \Delta:A\sstac L \cD\to A\sstac L
\cD\sstac L \cD\cong (A\sstac L \cD)\sstac A (A\sstac L \cD)$ and the counit is
$A\sstac L \epsilon:A\sstac L \cD\to A$. Via the canonical isomorphism
$M\sstac A A \sstac L \cD\cong M\sstac L\cD$, for any right $A$-module $M$,
right comodules for the $A$-coring $A\sstac L \cD$ are identified with {\bf
entwined modules}. A right-right entwined module means a right $A$-module and
right $\cD$-comodule $M$, with coaction $\tau_M:m\mapsto m_{[0]}\sstac L
m_{[1]}$, such that
$$
m\eta(l)=ml
\quad \textrm{and}\quad
\tau_M(ma)=m_{[0]}\psi(m_{[1]}\stac L a), \qquad \textrm{for }m\in M,
\ l\in L,\ a\in A.
$$
Morphisms of entwined modules are $A$-linear and $\cD$-colinear maps. The
category of right-right entwined modules is denoted by $\cM^\cD_A(\psi)$.

Entwining structures $(A,\cD,\psi)$ over an algebra $L$ provide
examples of coring extensions. Namely, the associated $A$-coring
$\cC :=(A\sstac{L}\cD,A\sstac{L} \Delta,A\sstac{L}\epsilon)$ is a
right $\cD$-comodule with coaction $\tau_\cC :=
A\sstac{L}\Delta:A\sstac{L}\cD\to A\sstac{L}\cD\sstac{L}\cD$. By
the coassociativity of the coproduct {$\Delta$} in $\cD$,
$\tau_\cC$ is left $\cC$-colinear. This means that the $L$-coring
$\cD$ is right extension of the $A$-coring $\cC$. 

{\color{blue}
\label{ex:entw.pure}
Note that any coring extension arising from an $L$-entwining structure
$(A,\cD,\psi)$ is pure. Use again that \eqref{eq:eq} is a split equalizer in
$\cM_A$. Thus the existence of a
forgetful functor $\cM_A \to \cM_L$ implies that \eqref{eq:eq} is a split
equalizer in $\cM_L$, hence it is preserved by any functor of domain $\cM_L$.}
In this
situation the functor ${\mathbb R}$ in Figure \eqref{fig:Brzfunc}
can be identified with the forgetful functor $\cM^\cC\cong
\cM_A^\cD(\psi)\to \cM^\cD$.

Let $(A,\cD,\psi)$ be an entwining structure over an algebra $L$
and $\cC :=A\sstac{L} \cD$ the associated $A$-coring. If $e$ is a
grouplike element in $\cD$ then $1_A\sstac{L} e$ is a grouplike
element in $\cC$. In this case $A$ is a right $\cC$-comodule hence
a right-right entwined module. The $\cD$-coaction in $A$ comes out
as
\begin{equation}\label{eq:A.entw.mod}
A\to A\stac L\cD,\qquad a\mapsto \psi(e\stac L a).
\end{equation}
The coinvariants of a right $\cC$-comodule (i.e. entwined module)
$M$ with respect to $1_A\sstac{L} e$ can be identified with
$\mathrm{Hom}^\cC(A,M)$, and the coinvariants of $M$ as a right
$\cD$-comodule with respect to $e$ can be identified with
$\mathrm{Hom}^\cD(L,{\mathbb R}(M))$ (cf.
\cite[28.4]{BrzWis:cor}). Since in this case the forgetful functor
${\mathbb R}:\cM^\cC\cong \cM^\cD_A(\psi)\to \cM^\cD$ possesses a
left adjoint, $\bullet\sstac{L} A$, it follows that
$\mathrm{Hom}^\cD(L,{\mathbb R}(M)) \cong \mathrm{Hom}^\cC(A,M)$.
That is to say, the coinvariants of a right $\cC$-comodule (i.e.
entwined module) with respect to $1_A\sstac{L} e$ are the same as its
coinvariants as a right $\cD$-comodule with respect to $e$.

If the entwining map $\psi$ is bijective then it induces an
$A$-coring structure in $\cD\sstac L A$. Its left comodules are
identified with left $A$-modules and left $\cD$-comodules,
satisfying a compatibility condition with $\psi$. If there exists
a grouplike element $e$ in $\cD$ then the corresponding left
$\cD$-coaction in $A$ is given by
\begin{equation}\label{eq:entw.mod.A}
A\to \cD\stac L A,\qquad a\mapsto \psi^{-1}(a\stac L e).
\end{equation}
\end{claim}

\begin{claim}\label{preli:bgd}
The notion of a {\bf bialgebroid} over an algebra $L$ was introduced by
Takeuchi in \cite{Tak:bgd} under the original name {\bf$\times_L$-bialgebra}.
Takeuchi's definition was shown by Brzezi\'nski and Militaru in
\cite{BrzMil:bgd} to be equivalent to the structure introduced in
\cite{Lu:hgd}. As a
$k$-bialgebra consists of compatible algebra and coalgebra structures on the
same $k$-module, an $L$-bialgebroid comprises compatible $L\sstac k
L^{op}$-ring and $L$-coring structures. More explicitly, a {\bf left
bialgebroid} is given
by the data $(H,L,s_L,t_L,\gamma_L,\pi_L)$. Here $H$ and $L$ are $k$-algebras
and $s_L:L\to H$ and $t_L:L^{op}\to H$ are algebra maps, called the source and
target maps, respectively. The map
$$
L\stac k L^{op}\to H,\qquad l\stac k l'\mapsto s_L(l) t_L(l')
$$
is required to be an algebra map, equipping $H$ with the structure of an
$L\sstac k L^{op}$-ring. The $L$-$L$ bimodule $H$, with actions
\begin{equation}\label{eq:bgd.left}
lhl'=s_L(l) t_L(l') h,\qquad \textrm{for } l,l'\in L,\ h\in H,
\end{equation}
is required to be an $L$-coring with coproduct $\gamma_L$ and
counit $\pi_L$. For the coproduct we use a Sweedler type index
notation with {\em lower} indexes, $\gamma_L(h)= h_{(1)}\sstac L
h_{(2)}$, for $h\in H$, where implicit summation is understood.
The compatibility axioms between the $L\sstac k L^{op}$-ring and
$L$-coring structures are the following. Consider the subset of
the $L$-module tensor square of the bimodule \eqref{eq:bgd.left},
the so called Takeuchi product
$$
H\times_L H =\{\ \sum_i h_i\stac L h'_i\in H \stac L H\ |\ \forall l\in L\
\sum_i h_it_L(l)\stac L h'_i=\sum_i h_i\stac L h'_is_L(l)\ \}.
$$
Note that $H\times_L H$ is an $L\sstac k L^{op}$-ring, with factorwise
multiplication and unit map
$$
L\stac k L^{op}\to H\times_L H, \qquad
l\stac k l'\mapsto s_L(l)\stac L t_L(l').
$$
The first bialgebroid axiom asserts that the coproduct
{corestricts to} a map of $L\sstac k L^{op}$-rings $H\to H\times_L
H$. The requirement, that the range of the coproduct lies within
$H\times_L H$, is referred to as the {\bf Takeuchi axiom}. Further
axioms require the counit to preserve the unit and satisfy
$$
\pi_L(h s_L\circ \pi_L(h'))=\pi_L(hh')=\pi_L(h t_L\circ \pi_L(h')),
$$
for all $h,h'\in H$.

The $L$-$L$ bimodule \eqref{eq:bgd.left} is defined in terms of left
multiplication by the source and target maps. Symmetrically, one defines {\bf
right bialgebroids} by interchanging the roles of left and right
multiplications. Explicitly, a right bialgebroid is given by the data
$(H,R,s_R,t_R,\gamma_R,\pi_R)$, where $H$ and $R$ are $k$-algebras and
$s_R:R\to H$ and $t_R:R^{op}\to H$ are algebra maps, called the source and
target maps, respectively. $H$ is required to be an $R\sstac k R^{op}$-ring
with unit
$$
R\stac k R^{op}\to H, \qquad r\stac k r'\mapsto s_R(r)t_R(r'),
$$
and an $R$-coring, with bimodule structure
\begin{equation}\label{eq:bgd.right}
rhr'=hs_R(r') t_R(r),\qquad \textrm{for } r,r'\in R, \ h\in H,
\end{equation}
coproduct $\gamma_R$ and counit $\pi_R$. For the coproduct we use a Sweedler
type index notation with {\em upper} indices, $\gamma_R(h)= h^{(1)} \sstac R
h^{(2)}$, for $h\in H$, where implicit summation is understood. The
coproduct is required to be a map of $R\sstac k R^{op}$-rings from $H$ to the
Takeuchi product
$$
H\times_R H =\{\ \sum_i h_i\stac R h'_i\in H \stac R H\ |\ \forall r\in R\
\sum_i s_R(r)h_i\stac R h'_i=\sum_i h_i\stac R t_R(r)h'_i\ \},
$$
where the $R$-module tensor product is taken with respect to the bimodule
structure \eqref{eq:bgd.right}. The counit is defined to preserve the unit and
satisfy
$$
\pi_R(s_R\circ \pi_R(h) h')=\pi_R(hh')=\pi_R(t_R\circ \pi_R(h) h'),
$$
for all $h,h'\in H$. For more details we refer to \cite{KadSzl:D2}.

The co-opposite $(H,L^{op},t_L,s_L,\gamma_L^{op},\pi_L)$ of a left bialgebroid
$(H,L,s_L,t_L,\gamma_L,\pi_L)$ is a left bialgebroid too.
The opposite $(H^{op},L,t_L,s_L,\gamma_L,\pi_L)$ is a right bialgebroid.
\end{claim}

\begin{claim}\label{preli:bgd.comodule}
A {\bf right comodule} of a right bialgebroid
$\hH_R=(H,R,s_R,t_R,\gamma_R,\pi_R)$ means a right comodule of the
$R$-coring $(H,\gamma_R,\pi_R)$ (with bimodule structure
\eqref{eq:bgd.right}). The category of right $\hH_R$-comodules is
denoted by $\cM^{\hH_R}$. In Section \ref{preli:coring}, a right
  $\hH_R$-comodule was defined to be in particular a
right $R$-module. Using the bialgebroid structure of $\hH_R$ (not
its coring structure alone), one can introduce also a left
$R$-module structure in a right $\hH_R$-comodule $M$,
\begin{equation}\label{eq:leftRmodM}
rm: = m^{[0]}\pi_R\big(s_R(r)m^{[1]}\big),\qquad \textrm{for } m\in M,\
r\in R.
\end{equation}
This makes $M$ an $R$-$R$ bimodule such that the (so called
Takeuchi) identity
$$
rm^{[0]}\stac R m^{[1]}=m^{[0]}\stac R t_R(r) m^{[1]}
$$
holds, for all $m\in M$ and $r\in R$. Any $\hH_R$-comodule map is
$R$-$R$ bilinear. This amounts to saying that there is a forgetful
functor $\cM^{\hH_R}\to {}_R \cM_R$. It was observed in
\cite[Proposition 5.6]{Scha:bgdnc} that the forgetful functor
$\cM^{\hH_R}\to {}_R \cM_R$ is strict monoidal. That is,
$\cM^{\hH_R}$ is a monoidal category via the $R$-module tensor
product. The coaction in the product $M\sstac R N$ of two right
$\hH_R$-comodules $M$ and $N$ is
$$
(m\stac R n)^{[0]}\stac R (m\stac R n)^{[1]}=
(m^{[0]}\stac R n^{[0]})\stac R m^{[1]} n^{[1]},\qquad
\textrm{for } m\stac R n\in M\stac R N.
$$
The monoidal unit is $R$ with coaction given by the source map
$s_R$. A {\bf right $\hH_R$-comodule algebra} is an algebra in the
monoidal category $\cM^{\hH_R}$ (hence it is in particular an
$R$-ring). Explicitly, it means an {$R$-ring}  and right
$\hH_R$-comodule $A$ whose coaction $\varrho^A$ satisfies
$$
\varrho^A(1_A)=1_A\stac R 1_H,\qquad
\varrho^A(aa')=a^{[0]}a^{\prime [0]}\stac R a^{[1]}a^{\prime [1]},\qquad
\textrm{for }a,a'\in A.
$$
The $R$-coring $(H,\gamma_R,\pi_R)$ underlying $\hH_R$ possesses a
grouplike element $1_H$. Coinvariants of a right $\hH_R$-comodule
are meant always with respect to the distinguished grouplike
element $1_H$. By the $R$-$R$ bilinearity of the coaction
$\varrho^A$ in a right $\hH_R$-comodule algebra $A$, for any
element $r$ in $R$ and any coinvariant $b$ in $A$, the unit map
$\eta:R\to A$ satisfies
$$
\varrho^A(b\eta(r))=b\stac R s_R(r)=\varrho^A(\eta(r)b).
$$
Hence the elements $b\in A^{co\hH_R}$ and $\eta(r)$, for $r\in R$, commute in
$A$.

{\bf Left comodules} of a right bialgebroid
$\hH_R=(H,R,s_R,t_R,\gamma_R,\pi_R)$
(i.e. of the $R$-coring $(H,\gamma_R,\pi_R)$) are treated symmetrically. Their
category is denoted by ${}^{\hH_R}\cM$. A left $\hH_R$-comodule $M$ (which is
a priori a left $R$-module) can be equipped with an $R$-$R$ bimodule structure
with right $R$-action
$$
mr=\pi_R\big(s_R(r) m^{[-1]}\big)m^{[0]}\qquad \textrm{for } m\in M,\
r\in R.
$$
The forgetful functor ${}^{\hH_R}\cM\to {}_{R^{op}}\cM_{R^{op}}$ is strict
monoidal. For two left $\hH_R$-comodules $M$ and $N$, the left and right
$R$-actions and the left $\hH_R$-coaction in the product $M\sstac {R^{op}} N$
take the form
$$
r(m\stac {R^{op}} n)r'=mr'\stac {R^{op}} rn,\qquad
(m\stac {R^{op}} n)^{[-1]}\stac R (m\stac {R^{op}} n)^{[0]}=
m^{[-1]} n^{[-1]}
\stac R (m^{[0]}\stac {R^{op}} n^{[0]}),
$$
for $r,r'\in R$ and $m\sstac {R^{op}} n \in M\sstac {R^{op}} N$.
The monoidal unit is $R^{op}$, with coaction given by the target
map $t_R$. A {\bf left $\hH_R$-comodule algebra} is defined as an
algebra in the monoidal category ${}^{\hH_R}\cM$. It is in
particular an $R^{op}$-ring. Explicitly, a left $\hH_R$-comodule
algebra is an {$R^{op}$-ring} and  left $\hH_R$-comodule $A$,
whose coaction ${}^A\varrho$ satisfies
$$
{}^A\varrho(1_A)=1_H\stac R 1_A,\qquad
{}^A\varrho(aa')=a^{[-1]}a^{\prime [-1]}\stac R a^{[0]}a^{\prime [0]},\qquad
\textrm{for }a,a'\in A.
$$
Coinvariants of left $\hH_R$-comodules are meant always with respect to the
distinguished grouplike element $1_H$.

Comodules of left bialgebroids can be described symmetrically.
{For a right bialgebroid $\hH_R$, the categories
${}^{(\hH_R)_{cop}}\cM$ and $\cM^{\hH_R}$ are monoidally
isomorphic. The categories $\cM^{(\hH_R)^{op}}$ and $\cM^{\hH_R}$
are anti-monoidally isomorphic.}
\end{claim}

\begin{claim}\label{rem:entwbgd}
Let ${\mathcal H}_R=(H,R,s_R,t_R,\gamma_R,\pi_R)$ be a right bialgebroid and
$A$ a right $\hH_R$-comodule algebra {with coaction $a\mapsto
  a^{[0]}\sstac R a^{[1]}$}. Recall from {Section}
  \ref{preli:bgd.comodule} that $A$ possesses an $R$-ring structure.
The $R$-ring $A$, the $R$-coring
$(H,\gamma_R, \pi_R)$ and the $R$-$R$ bimodule map
\begin{equation}\label{eq:hgdentwR}
\psi_R:H\stac{R} A\to A\stac{R} H,\qquad h\stac{R} a\mapsto a^{[0]}\stac{R}
ha^{[1]}
\end{equation}
form an entwining structure over $R$. This implies that $A\otimes_R H$ is an
$A$-coring, with bimodule structure
$$
a_1(a\stac{R} h)a_2=a_1a{a_2}^{[0]}\stac{R} h{a_2}^{[1]},\qquad \textrm{for }
a_1,a_2\in A\textrm{ and }\ a\stac{R} h\in A\stac{R} H,
$$
coproduct $A\sstac{R} \gamma_R:A\sstac{R}H\to A\sstac{R}
H\sstac{R} H\cong (A\sstac{R} H)\sstac{A}(A\sstac{R} H)$ and
counit $A\sstac{R}\pi_R:A\sstac{R} H\to A$. Via the canonical
isomorphism $M\sstac{A}(A\sstac{R}H)\cong M\sstac{R}H$, {for $M\in
\cM_A$,} right comodules for the $A$-coring $(A\sstac{R} H,
A\sstac{R}\gamma_R, A\sstac{R}\pi_R)$ can be identified with
right-right entwined modules for the entwining structure
$(A,H,\psi_R)$. Such entwined modules are also called {\bf
right-right $(A,{\mathcal H}_R)$-relative Hopf modules}. They can
be described equivalently as right modules for the algebra $A$ in
the category of right $\hH_R$-comodules. That is, right
$A$-modules and right $\hH_R$-comodules $M$, such that the
$A$-action is $\hH_R$-colinear, in the sense that the
compatibility condition
$$
(ma)^{[0]}\stac{R} (ma)^{[1]}=m^{[0]} a^{[0]}\stac{R} m^{[1]} a^{[1]}
$$
holds, for $m\in M$, $a\in A$. The category of right-right
$(A,{\mathcal H}_R)$-relative Hopf modules will be denoted by
$\cM^{{\mathcal H}_R}_A$. As it is explained in {Section}
\ref{rem:entwcorext}, in the $R$-entwining structure
$(A,H,\psi_R)$ the R-coring $(H,\gamma_R,\pi_R)$ is a right
extension of the $A$-coring $(A\sstac{R} H, A\sstac{R}\gamma_R,
A\sstac{R}\pi_R)$. For this coring extension, the functor
${\mathbb R}$ on Figure \eqref{fig:Brzfunc} can be identified with
the forgetful functor $\cM^{{\mathcal H}_R}_A\to \cM^{{\mathcal
H}_R}$.

A right $\hH_R$-comodule algebra $A$ is called an {\bf
$\hH_R$-Galois extension} of its coinvariant subalgebra $B$ if the
associated $A$-coring $(A\sstac R H,A\sstac R \gamma_R,A\sstac R
\pi_R)$, with grouplike element $1_A\sstac R 1_H$, is a Galois
coring, i.e. the canonical map
\begin{equation}\label{eq:canRR}
\can:A\stac B A \to A\stac R H,\qquad a\stac B a'\mapsto a a^{\prime [0]}\stac
R a^{\prime [1]}
\end{equation}
is bijective.

For a left comodule algebra $A'$ for a left bialgebroid ${\mathcal H}_L$ one
defines left-left $(A',{\mathcal H}_L)$-relative Hopf modules in a symmetrical
way.
\end{claim}

\begin{claim}\label{preli:hgd}
A {\bf Hopf algebroid} \cite{BohmSzl:hgdax}, \cite{Bohm:hgdint} is
a triple $\hH=(\hH_L,\hH_R,S)$. It consists of a left bialgebroid
$\hH_L=(H,L,s_L,t_L,\gamma_L,\pi_L)$ and a right bialgebroid
$\hH_R=(H,R,s_R,t_R,\gamma_R,\pi_R)$ on the {\em same} total
algebra $H$. They are subject to the following compatibility
axioms
\begin{eqnarray}\label{eq:hgd.ax.base.iso}
&s_R\circ \pi_R\circ t_L=t_L\qquad &s_L\circ \pi_L\circ t_R=t_R\\
&t_R\circ \pi_R\circ s_L=s_L\qquad &t_L\circ \pi_L\circ s_R=s_R
\nonumber
\end{eqnarray}
and
\begin{equation}\label{eq:hgd.ax.cor.ext}
(\gamma_R \stac L H)\circ \gamma_L=(H\stac R \gamma_L)\circ \gamma_R,\qquad
(\gamma_L \stac R H)\circ \gamma_R=(H\stac L \gamma_R)\circ \gamma_L.
\end{equation}
The $k$-linear map $S:H\to H$ is called the {\bf antipode}. It is required to
be $R$-$L$ bilinear in the sense that
$$
S\big(t_L(l)ht_R(r)\big)=s_R(r) S(h) s_L(l),\qquad \textrm{for }l\in L,\ r \in
R,\ h\in H.
$$
The antipode axioms read as
$$
\mu\circ(H\stac R S)\circ \gamma_R=s_L\circ \pi_L,\qquad
\mu\circ(S\stac L H)\circ \gamma_L=s_R\circ \pi_R,
$$
where $\mu$ denotes the multiplication both in the $L$-ring $s_L:L\to H$ and
the $R$-ring $s_R:R\to H$.

In a Hopf algebroid there are two bialgebroid (hence two coring)
structures present. Throughout this paper we insist on using upper
indices of the Sweedler type to denote components of the
coproduct and coactions of the right bialgebroid $\hH_R$, and
lower indices in the case of the left bialgebroid $\hH_L$.

{Similarly to the case of Hopf algebras, the antipode of a Hopf
  algebroid $\hH=(\hH_L,\hH_R,S)$ is an anti-algebra map on the total algebra
  $H$. That is,
$$
S(1_H)=1_H\quad \textrm{and}\quad S(hh')=S(h')S(h),\qquad \textrm{for }h,h'\in
H.
$$
It is also an anti-coring map $\hH_L\to \hH_R$ and $\hH_R \to \hH_L$. That is,
\begin{eqnarray*}
&\pi_R\circ S =\pi_R\circ s_L\circ \pi_L\quad \textrm{and}\quad
&S(h)^{(1)}\stac R S(h)^{(2)}=S(h_{(2)})\stac R S(h_{(1)}),\\
&\pi_L\circ S =\pi_L\circ s_R\circ \pi_R\quad \textrm{and}\quad
&S(h)_{(1)}\stac L S(h)_{(2)}=S(h^{(2)})\stac L S(h^{(1)}),
\qquad \textrm{for }h\in H.
\end{eqnarray*}}

For a Hopf algebroid $\hH=(\hH_L,\hH_R,S)$, also the opposite-co-opposite
$\hH^{op}_{cop}=((\hH_R)^{op}_{cop},(\hH_L)^{op}_{cop},$ $S)$ is a Hopf
algebroid. If the antipode $S$ is bijective then so are the opposite
$\hH^{op}=((\hH_R)^{op},(\hH_L)^{op},$ $S^{-1})$ and the co-opposite
$\hH_{cop}=((\hH_L)_{cop},(\hH_R)_{cop},S^{-1})$, too.
\end{claim}

{\color{blue}
\begin{claim}\label{preli:hgd.comodule}
The following definition was proposed in \cite[Definition 3.2]{Bohm:hgdGal}
and \cite[Section 2.2]{BalSzl:fin.Gal}.

A {\bf right comodule} of a Hopf algebroid $\hH$ is a right
$L$-module as well as a right $R$-module $M$, together with a right
coaction $\varrho_R: M \to M\ot_R H$ of the constituent right bialgebroid
$\hH_R$ and a right coaction $\varrho_L:M\to M\ot_L H$ of the constituent
left bialgebroid $\hH_L$,
such that $\varrho_R$ is an $\hH_L$-comodule map and $\varrho_L$ is an
$\hH_R$-comodule map. Explicitly,
\begin{equation}\label{eq:hgd_coac}
(M\ot_R \gamma_L)\circ \varrho_R = (\varrho_R \ot_L H) \circ \varrho_L
\qquad \textrm{and}\qquad
(M\ot_L \gamma_R)\circ \varrho_L = (\varrho_L \ot_R H) \circ \varrho_R.
\end{equation}
Morphisms of $\hH$-comodules are $\hH_R$-comodule maps as well as
$\hH_L$-comodule maps. The category of right $\hH$-comodules is denoted by
$\cM^\hH$.

Note that any right $\hH$-comodule is a right $R\ot L$-module.

The category ${}^\hH\cM$ of left $\hH$-comodules is defined symmetrically.
\end{claim}

\begin{claim}\label{prop:hgd_coinv}
Since a comodule $M$ of a Hopf algebroid $\hH$ is a comodule of both constituent
bialgebroids $\hH_L$ and $\hH_R$, we can consider the coinvariants $M^{co
  \hH_R}$ and $M^{co \hH_L}$ in the sense of Appendix A.6.
By \cite[Corrigendum]{BohmBrz:cleft}, for any $\hH$-comodule $M$, $M^{co
  \hH_R}\subseteq M^{co \hH_L}$. If the antipode of $\hH$ is bijective then an
equality holds.
\end{claim}

\begin{claim}\label{prop:ff_adj}
For any Hopf algebroid $\hH$ the following hold.

1) The forgetful functor $\cM^\hH \to \cM_L$ possesses a
  right adjoint $-\ot_L H$.

2) The forgetful functor $\cM^\hH \to \cM_R$ possesses a
  right adjoint $-\ot_R H$.

\begin{proof}
1) The unit of the adjunction is given by the $\hH_L$-coaction $M \to M\ot_L
  H$, for any right $\hH$-comodule $M$. It is an $\hH$-comodule map by
  definition. Counit is given by $N\ot_L \pi_L:N\ot_L H \to N$, for any right
  $L$-module $N$.
Part 2) is proven symmetrically.
\end{proof}
\end{claim}

\begin{claim}\label{thm:pure}\cite[Corrigendum]{BohmBrz:cleft}
Consider a Hopf algebroid $\hH$. Denote by $F_R$ and $F_L$ the
forgetful functors $\cM^{\hH_R} \to \cM_k$ and $\cM^{\hH_L} \to \cM_k$,
respectively.

1) If the equalizer
\begin{equation}\label{eq:R.M}
\xymatrix{
M \ar[rr]^-{\varrho_R}&&
M\ot_R H \ar@<2pt>[rr]^-{\varrho_R\ot_R H}\ar@<-2pt>[rr]_-{M\ot_R \gamma_R}&&
M\ot_R H \ot_R H
}
\end{equation}
in $\cM_L$ is $H\ot_L H$-pure, i.e. it is preserved by the functor $-\ot_L
H\ot_L H:\cM_L\to \cM_L$, for any right $\hH_R$-comodule $(M,\varrho_R)$, then
there exists a functor $U:\cM^{\hH_R}\to \cM^{\hH_L}$, such that $F_L \circ U
=F_R$. Moreover, in this case the forgetful
functor $G_R:\cM^\hH\to \cM^{\hH_R}$ is fully faithful.

2) If the equalizer
$$
\xymatrix{
N \ar[rr]^-{\varrho_L}&&
N\ot_L H \ar@<2pt>[rr]^-{\varrho_L\ot_L H}\ar@<-2pt>[rr]_-{N\ot_L \gamma_L}&&
N\ot_L H \ot_L H
}
$$
in $\cM_R$ is $H\ot_R H$-pure, i.e. it is preserved by the functor $-\ot_R
H\ot_R H:\cM_R\to \cM_R$, for any right $\hH_L$-comodule $(N,\varrho_L)$, then
there exists a functor $V:\cM^{\hH_L}\to \cM^{\hH_R}$, such that $F_R \circ V
=F_L$.
Moreover, in this case the forgetful
functor $G_L:\cM^\hH\to \cM^{\hH_L}$ is fully faithful.

3) If both purity assumptions in parts 1) and 2) hold, then the forgetful
functors $G_R:\cM^\hH\to \cM^{\hH_R}$ and $G_L:\cM^\hH\to \cM^{\hH_L}$
are isomorphisms. Moreover, $G_L \circ G_R^{-1}=U$ and $G_R \circ G_L^{-1}=V$,
hence $U$ and $V$ are inverse isomorphisms.
\end{claim}

\begin{claim}\label{thm:hgd.com.mon}\cite[Corrigendum]{BohmBrz:cleft}
For any Hopf algebroid $\hH$,
the category $\cM^{\mathcal H}$ of right $\hH$-comodules is
monoidal. Moreover, 
the following diagram is commutative
and all occurring forgetful functors are strict monoidal. 
$$
\xymatrix{
\cM^{\mathcal H} \ar[r]\ar[d]&
\cM^{\hH_R} \ar[d]\\
\cM^{\hH_L} \ar[r]&
{}_R \cM_R.
}
$$
\label{def:Hopf_mod}
In light of this observation, the following definition can be made.

A {\bf right comodule algebra} of a Hopf algebroid $\hH$ is an algebra
in the monoidal category $\cM^\hH$. Right/left modules of a right
$\hH$-comodule algebra $A$ {\em in} $\cM^\hH$ are termed {\bf
(right-right/left-right) relative Hopf modules}. Their categories are denoted by
$\cM^\hH_A$ and ${}_A \cM^\hH$, respectively.
\end{claim}

\begin{claim}\label{prop:coinv}
For a right comodule algebra $A$ of a Hopf algebroid $\hH$, denote $B:=A^{co
  \hH_R}$. Then there is an adjunction
$$
-\ot_B A:\cM_B \to \cM^\hH_A\qquad (-)^{co \hH_R}:\cM^\hH_A \to \cM_B.
$$

\begin{proof}
For any right $B$-module $N$, the unit of the adjunction is given by
$$
N \to (N\ot_B A)^{co \hH_R},\qquad n\mapsto n \ot_B 1_A.
$$
For any relative Hopf module $M\in \cM^\hH_A$, counit is given by
$$
M^{co \hH_R}\ot_B A \to M,\qquad m\ot_B a \mapsto ma.
$$
Obviously, it is a right $A$-module map. In light of 
\ref{prop:hgd_coinv}, it
is also a morphism of $\hH$-comodules. Verification of the adjunction
relations is a routine computation.
\end{proof}
\end{claim}

\begin{claim}\label{prop:two_adj}
For a Hopf algebroid $\hH$ and a right $\hH$-comodule algebra $A$, denote
$B:=A^{co\hH_R}$.

1) The functor $-\ot_B A:\cM_B \to \cM^{\hH_R}_A$ is fully faithful if and only
if the functor $-\ot_B A:\cM_B \to \cM^\hH_A$ is fully faithful.

2) If the functor $-\ot_B A:\cM_B \to \cM^{\hH_R}_A$ is an equivalence then
also the functor $-\ot_B A:\cM_B \to \cM^\hH_A$ is an equivalence.

\begin{proof}
Consider the adjunction in Appendix \ref{prop:coinv} and the adjunction
\begin{equation}\label{eq:bgd_adj}
-\ot_B A:\cM_B \to \cM^{\hH_R}_A\qquad (-)^{co\hH_R}:\cM^{\hH_R}_A \to \cM_B,
\end{equation}
cf. \eqref{eq:co_un}-\eqref{eq:co_coun}.
Both statements follow by noticing that the units of the two adjunctions
coincide and counit of the adjunction in Appendix
\ref{prop:coinv} is equal to the restriction of the counit of the adjunction
\eqref{eq:bgd_adj} to the objects of $\cM^\hH_A$.
\end{proof}
\end{claim}}

\begin{claim}\label{rem:entwhgd}
Let ${\mathcal H}$ be a Hopf algebroid with
constituent left bialgebroid ${\mathcal H}_L=(H,L,s_L,t_L,\gamma_L,\pi_L)$,
right bialgebroid  ${\mathcal H}_R=(H,R,s_R,t_R,\gamma_R,\pi_R)$, and antipode
$S$, and let $A$ be a right ${\mathcal H}$-comodule algebra.
This means in particular that $A$ is a right comodule algebra for the
right $R$-bialgebroid ${\mathcal H}_R$, with coaction $a\mapsto a^{[0]}\sstac
R a^{[1]}$.
What is more, since $A$ is a right comodule algebra for the left bialgebroid
$\hH_L$ as well, with coaction $a\mapsto a_{[0]}\sstac L a_{[1]}$,
related to the $\hH_R$-coaction as in
{\color{blue}
\eqref{eq:hgd_coac}}, the opposite algebra $A^{op}$ is a right
comodule algebra for the right $L$-bialgebroid $({\mathcal
H}_L)^{op}$. Hence in addition to the
$R$-entwining structure \eqref{eq:hgdentwR}, $A$ determines also an
$L$-entwining structure. It consists of the $L$-ring $A^{op}$
(with unit, expressed in terms of the unit $\eta$ of the $R$-ring $A$ as
$\eta\circ \pi_R\circ t_L$), the $L$-coring $(H,\gamma_L,\pi_L)$, and the
entwining map
\begin{equation}\label{eq:hgdentwL}
\psi_L:H\stac{L} A\to A\stac{L} H,\qquad h\stac{L} a\mapsto a_{[0]}\stac{L}
a_{[1]}h.
\end{equation}
Therefore there is an associated $A^{op}$-coring structure on $A\sstac{L} H$.

Note that the entwining map \eqref{eq:hgdentwL} is bijective with
inverse $\psi_L^{-1}(a\sstac L h)=S(a^{[1]})h \sstac L a^{[0]}$.
Hence $H\sstac{L} A$ has a unique $A^{op}$-coring structure such
that \eqref{eq:hgdentwL} is an isomorphism of corings. Clearly, by
the existence of grouplike elements, $A^{op}$ is a left comodule
for the $A^{op}$-corings $A\sstac{L} H\cong H\sstac{L} A$.
\end{claim}

{\color{blue}
\begin{claim}\label{preli:left-right}
The antipode $S$ of a Hopf algebroid induces strict anti-monoidal functors
${}^{\hH_R}\cM \to \cM^{\hH_L}$, 
${}^{\hH_L}\cM \to \cM^{\hH_R}$ and
${}^\hH\cM\to \cM^\hH$:
Let $M$ be a left ${\mathcal H}_R$-comodule with
coaction $m\mapsto m^{[-1]}\sstac R m^{[0]}$.
Then $M$ has a right $\hH_L$-comodule structure with right $L$-action
$ml:=\pi_R \circ t_L(l)m$, for $l\in L$ and $m\in M$,
and coaction
\begin{equation}\label{eq:MleftH_Lcom}
m\mapsto m^{[0]}\stac L S(m^{[-1]}).
\end{equation}
If $M$ is a left ${\mathcal H}_L$-comodule with coaction $m\mapsto
m_{[-1]}\sstac L m_{[0]}$, then $M$ has
a right $\hH_R$-comodule structure, with right
$R$-action $mr:=\pi_L\circ t_R(r)m$, for $r\in R$ and $m\in M$,
and coaction
\begin{equation}\label{eq:MleftH_Rcom}
m\mapsto m_{[0]}\stac R S(m_{[-1]}).
\end{equation}
If $M$ is a left $\hH$-comodule 
then the $\hH_R$-coaction \eqref{eq:MleftH_Rcom} and the $\hH_L$-coaction
\eqref{eq:MleftH_Lcom} are checked to constitute a right $\hH$-comodule
structure on $M$. 

Clearly, if $S$ is bijective, then all these functors are isomorphisms.
Therefore, $A$ is a right $\hH$-comodule algebra if and only if the opposite
algebra $A^{op}$ possesses a left $\hH$-comodule algebra structure.
For a right comodule algebra $A$ of a Hopf algebroid $\hH$
with a bijective antipode, left/right
$A^{op}$-modules in ${}^\hH\cM$ are called {\bf (right-left/left-left)
  relative Hopf modules}.  Their categories are denoted by
${}^\hH \cM_A$ and ${}_A ^\hH \cM$, respectively.

In particular, left-left relative Hopf modules are left $A$-modules and left
$\hH$-comodules, subject to the compatibility conditions
\begin{eqnarray}\label{eq:leftentw}
&&(am)^{[-1]}\stac R (am)^{[0]}=m^{[-1]} S^{-1}(a_{[1]})\stac R a_{[0]} m^{[0]}
\qquad \textrm{and}\\
&&{\color{blue}
(am)_{[-1]}\stac L (am)_{[0]}=m_{[-1]} S^{-1}(a^{[1]})\stac L a^{[0]} m_{[0]}}
\qquad \textrm{for } a\in A, \ m\in M.\nonumber
\end{eqnarray}
\end{claim}}

\begin{claim}\label{rem:bijantip}
Let ${\mathcal H}$ be a Hopf algebroid with constituent left
bialgebroid ${\mathcal H}_L=(H,L,s_L,t_L,\gamma_L,\pi_L)$, right
bialgebroid  ${\mathcal H}_R=(H,R,s_R,t_R,\gamma_R,\pi_R)$, and a
{\em bijective} antipode $S$. Let $A$ be a right ${\mathcal
H}$-comodule algebra with $\hH_R$-coaction $a\mapsto a^{[0]}\sstac
R a^{[1]}$ and $\hH_L$-coaction $a\mapsto a_{[0]}\sstac L
a_{[1]}$, related via \eqref{eq:hgd_coac}.
The two isomorphic $A^{op}$-corings $A\sstac L H\cong H\sstac L A$
in Section \ref{rem:entwhgd} are anti-isomorphic to the
$A$-corings $A\sstac R H \cong H\sstac R A$. An anti-isomorphism
is given by a bijection in \cite[Lemma 3.3]{Bohm:hgdGal},
\begin{equation}\label{eq:phi}
A\stac R H \to A\stac L H,\qquad a\stac R h \mapsto a_{[0]}\stac L
a_{[1]}S(h),
\end{equation}
and an isomorphism $H\sstac R A \to A\sstac R H$ is given by the entwining map
\eqref{eq:hgdentwR}, with inverse
\begin{equation}\label{eq:psi_R_inv}
a\stac R h \mapsto h S^{-1}(a_{[1]})\stac R a_{[0]},
\end{equation}
cf. \cite[Lemma 4.1]{Bohm:hgdGal}.
By the existence of grouplike elements
(given by the units in $A$ and $H$), $A$ is a left comodule for
all these corings.
\end{claim}

\section*{Acknowledgements}
This paper was written while A. Ardizzoni and C. Menini were
members of G.N.S.A.G.A. with partial financial support from
M.I.U.R.. Work of G. B\"ohm is supported by the Hungarian
Scientific Research Fund OTKA T043159 and the Bolyai J\'anos
Fellowship. {Her stay, as a visiting professor at University of
Ferrara, was supported by I.N.D.A.M..} She would like to express
her gratitude to the members of the Department of Mathematics at
University of Ferrara for a very warm hospitality.


\begin{thebibliography}{99}
\bibitem[Ar]{Ar2} A. Ardizzoni, \textit{Separable Functors and formal
  smoothness}, Journal of K-Theory \textbf{1} (2008), 535--582.

\bibitem[BaSz]{BalSzl:fin.Gal} I. B\'alint and K. Szlach\'anyi, {\em Finitary
  Galois  extensions over noncommutative bases}, J. Algebra {\bf 296} (2006),
  520--560.

\bibitem[B\"o1]{Bohm:hgdint} G.\ B\"ohm, {\em Integral theory for Hopf
  algebroids}, Alg. Rep. Theory \textbf{8} (2005), 563--599.
{\color{blue}
{\em Corrigendum}, to be published. See also  
\href{http://arxiv.org/abs/math/0403195}{\tt
http://arxiv.org/abs/math/0403195}, to be replaced by version 4.}

\bibitem[B\"o2]{Bohm:hgdGal} G.\ B\"ohm, {\em Galois theory for Hopf
  algebroids}, Ann. Univ. Ferrara - Sez. VII - Sc. Mat., Vol {\textbf{LI}}
  (2005), 233--262.
{\color{blue}
A corrected version is available at 
\href{http://arxiv.org/abs/math/0409513}{\tt
  http://arxiv.org/abs/math/0409513v3}.}

\bibitem[BB1]{BohmBrz:relchern} G. B\"ohm and T. Brzezi\'nski, \emph{Strong
  connections and the relative Chern-Galois character for corings},
  Int. Math. Res. Not. 2005, \textbf{42} (2005), 2579--2625.

\bibitem[BB2]{BohmBrz:cleft} G. B\"ohm and T. Brzezi\'nski, \emph{Cleft
  extensions of Hopf algebroids},
  Appl. Categorical Structures \textbf{14} (2006), 431--469.
{\color{blue}
{\em Corrigendum}, to be published.
See also 
\href{http://arxiv.org/abs/math/0510253}{\tt
  http://arxiv.org/abs/math/0510253v2}.}
 
\bibitem[BB3]{BohmBrz:pre_tor} G. B\"ohm and T. Brzezi\'nski,
  \emph{Pre-torsors and equivalences},  J. Algebra {\bf 317} (2007), 544--580.
{\em Corrigendum}, J. Algebra {\bf 319} (2008), 1339 --1340.

\bibitem[BSz]{BohmSzl:hgdax} G. B\"ohm and K. Szlach\'anyi, \emph{Hopf
  algebroids with bijective antipodes: axioms, integrals and duals}, J. Algebra
  \textbf{274} (2004), 708--750.

\bibitem[Bro]{Br} K. S. Brown, \emph{Cohomology of groups.} Graduate Texts in
Mathematics, \textbf{87}. Springer-Verlag, New York-Berlin, 1982.

\bibitem[Brz1]{Brz:str} T.\ Brzezi\'nski,
{\it The structure of corings. Induction functors, Maschke-type
theorem, and Frobenius and Galois-type properties,} Alg.\ Rep.\
Theory \textbf{5} (2002), 389--410.

\bibitem[Brz2]{Brz:Galcom} T. Brzezi\'nski, \emph{Galois comodules}, J. Algebra
  \textbf{290} (2005), 503--537.

\bibitem[Brz3]{Brz:Corext} T.\ Brzezi\'nski, {\em A note on coring extensions},
Ann. Univ. Ferrara - Sez. VII - Sc. Mat., Vol {\textbf{LI}}
(2005), 15--27.
{\color{blue}
A corrected version is available at 
\href{http://arxiv.org/abs/math/0410020}{\tt
  http://arxiv.org/abs/math/0410020v3}.}

\bibitem[BM]{BrzMil:bgd} T.\ Brzezi\'nski and G. Militaru,
\emph{Bialgebroids, $\times_R$-bialgebras and duality}, J. Algebra
\textbf{251} (2002), 279--294.

\bibitem[BTW]{BrzTurWri:WeakcoaGal} T. Brzezi\'nski, R.B Turner and
  A.P. Wrightson, \emph{The structure of weak coalgebra-Galois extensions},
{Comm. Algebra {\bf 34} (2006), 1489--1519.}

\bibitem[BW]{BrzWis:cor} T.\ Brzezi\'nski and R.\ Wisbauer, {\it Corings and
 Comodules}. Cambridge University Press, Cambridge, 2003.
 Erratum:
\href{http://www-maths.swan.ac.uk/staff/tb/corinerr.pdf}{\tt
 http://www-maths.swan.ac.uk/staff/tb/corinerr.pdf}. 

\bibitem[CM]{CaeMil:SecK} S. Caenepeel and G. Militaru, \emph{Maschke functors,
    semisimple functors and separable functors of the second
    kind. Applications,} J. Pure Appl. Algebra \textbf{178} (2003), 131--157.

\bibitem[CIMZ]{CaeIonMilZhu:SepDoiH} S. Caenepeel, B. Ion, G. Militaru and
  S. Zhu, \emph{Separable functors applied to Doi Hopf modules. Applications},
  Adv. Math. \textbf{145} (1999), 239--290.

\bibitem[CMZ]{CMZ}  S. Caenepeel, G. Militaru and Shenglin Zhu, \emph{Frobenius
Separable Functors for Generalized Module Categories and Nonlinear Equations}.
LNM \textbf{1787} (2002), Springer-Verlag, Berlin - New York.

\bibitem[DGH]{DaGrHa:Strconn} L. Dabrowski, H. Grosse and P.M. Hajac,
  \emph{Strong connections and Chern-Connes pairing in the Hopf-Galois
  theory}, Comm. Math. Phys. \textbf{220} (2001), 301--331.

\bibitem[Doi]{Doi:totint} Y. Doi, \emph{Algebras with total integrals},
  Comm. Algebra {\bf 13} (1985), 2137--2159.

\bibitem[H] {Haj:Strconn} P.M. Hajac, \emph{Strong connections on
    quantum principal bundles}, Comm. Math. Phys. \textbf{182} (1996),
    579--617.

\bibitem[HM]{HaMa:Projmod} P.M. Hajac and S. Majid, \emph{Projective
    module  description of the q-monopole}, Comm. Math. Phys. \textbf{206}
    (1999), 247--264.

\bibitem[HS]{HS}  P.J. Hilton and U. Stambach, \emph{A course in
Homological algebra}. Graduate Text in Mathematics \textbf{4},
Springer, New York, 1971.

\bibitem[KSz]{KadSzl:D2} L. Kadison and K. Szlach\'anyi, \emph{Bialgebroid
  actions on depth two extensions and duality}, Adv. Math. \textbf{179}
  (2003), 75--121.

\bibitem[KT]{KreTak:HopfGal} H. F. Kreimer and M. Takeuchi, {\em Hopf algebras
and Galois extensions of an algebra}, Indiana Univ. Math. J. {\bf
30} (1981), 675--692.

\bibitem[Lu]{Lu:hgd} J.-H. Lu, \emph{Hopf algebbroids and quantum groupoids},
  Int. J. Math. \textbf{7} (1996), 47--70.

\bibitem[MM1]{MenMil:intYD} C. Menini and G. Militaru, \emph{Integrals, quantum
    Galois extensions and the affineness criterion for quantum
    Yetter-Drinfel'd modules},  J. Algebra \textbf{247} (2002), 467--508.

\bibitem[MM2]{MenMil:affDK} C. Menini and G. Militaru, \emph{The affineness
    criterion for Doi-Koppinen modules}, in: "Hopf algebras in non-commutative
    geometry and physics", S. Caenepeel and F. Van Oystaeyen (eds.), Lecture
    Notes in Pure and Applied Mathematics \textbf{239} pp 215--228, Marcel
    Dekker, New York, 2004.

\bibitem[NVV]{NdO}  C. N\u{a}st\u{a}sescu, M. Van den Bergh and F.
Van Oystaeyen, \emph{Separable functors applied to graded rings},
J.Algebra \textbf{123} (1989), 397--413.

\bibitem[Raf]{Rafael}  M. D. Rafael, \emph{Separable Functors Revisited, }Comm.
Algebra \textbf{18} (1990), 1445--1459.

\bibitem[Row]{Row:rin}
L.H.\ Rowen,  {\em Ring theory. {V}ol.\ {I}}, Academic Press,
Boston (1988).

\bibitem[Scha1]{Scha:HGal&biGal} P. Schauenburg, \emph{Hopf-Galois and
    bi-Galois extensions}, in: Galois theory, Hopf algebras, and semiabelian
    categories, Fields Inst. Commun. \textbf{43}, AMS 2004, pp. 469--515.

\bibitem[Scha2]{Scha:bgdnc} P. Schauenburg, \emph{Bialgebras over
    noncommutative rings and a structure theorem for Hopf bimodules},
    Appl. Categorical Structures \textbf{6} (1998), 193--22.

\bibitem[SS]{SchSch:genGal} P. Schauenburg and H.-J. Schneider, \emph{On
    generalized Hopf Galois
    extensions}, Journal of Pure and Applied Algebra \textbf{202} (2005),
    168--194.

\bibitem[Schn]{Schn:PriHomS} H.-J. Schneider, \emph{Principal homogeneous
    spaces for arbitrary Hopf algebras}, Israel J. Math. \textbf{72} (1990),
    167--195.

\bibitem[Ta]{Tak:bgd} M. Takeuchi, \emph{Groups of algebras over $A\otimes
  {\overline A}$}, J. Math. Soc. Japan \textbf{29} (1977), 459--492.

\bibitem[We]{Weibel}  C. Weibel, \emph{An introduction to homological algebra.
} Cambridge Studies in Advanced Mathematics \textbf{38}, Cambridge University
Press, 1994.
\end{thebibliography}
\end{document}